\documentclass[final,onefignum,onetabnum]{siamonline171218}

\usepackage{lipsum}
\usepackage{amsfonts}
\usepackage{graphicx}
\usepackage{epstopdf}
\usepackage{algorithmic}
\usepackage{amssymb}
\usepackage{mathpazo}
\usepackage{mathrsfs}
\usepackage{cancel}
\ifpdf
  \DeclareGraphicsExtensions{.eps,.pdf,.png,.jpg}
\else
  \DeclareGraphicsExtensions{.eps}
\fi
\newcommand{\dd}{\mathrm{d}}
\newcommand{\ee}{\mathrm{e}}
\newcommand{\ii}{\mathrm{i}}
\renewcommand{\varepsilon}{\epsilon}

\usepackage[shortlabels]{enumitem}
\setlist[enumerate]{leftmargin=.5in}
\setlist[itemize]{leftmargin=.5in}

\newsiamremark{remark}{Remark}
\newsiamremark{hypothesis}{Hypothesis}
\crefname{hypothesis}{Hypothesis}{Hypotheses}
\newsiamthm{claim}{Claim}
\newsiamthm{example}{Example}

\headers{Benjamin-Ono Small Dispersion Limit}{E. Blackstone, L. Gassot, P. G\'erard, and P. D. Miller}

\title{The {B}enjamin-{O}no equation in the zero-dispersion limit for rational initial data:  generation of dispersive shock waves
}

\author{Elliot Blackstone\thanks{Department of Mathematics, University of Michigan, Ann Arbor, MI 
  (\email{eblackst@umich.edu}).}
\and Louise Gassot\thanks{CNRS and Department of Mathematics, University of Rennes, Rennes, France
	(\email{louise.gassot@cnrs.fr}).}
 \and
 Patrick G\'erard\thanks{Laboratoire de Mathématiques d'Orsay,  Université Paris-Saclay, Orsay, France (\email{patrick.gerard@universite-paris-saclay.fr}).}
\and Peter D. Miller\thanks{Department of Mathematics, University of Michigan, Ann Arbor, MI 
  (\email{millerpd@umich.edu}).}
}

\usepackage{amsopn}

\makeatletter
\newcommand*{\addFileDependency}[1]{
  \typeout{(#1)}
  \@addtofilelist{#1}
  \IfFileExists{#1}{}{\typeout{No file #1.}}
}
\makeatother

\newcommand*{\myexternaldocument}[1]{%
    \externaldocument{#1}%
    \addFileDependency{#1.tex}%
    \addFileDependency{#1.aux}%
}

\ifpdf
\hypersetup{
  pdftitle={Strong Zero-Dispersion Limit for Benjamin-Ono},
  pdfauthor={Elliot Blackstone and Louise Gassot and Patrick Gérard and Peter D. Miller}
}
\fi

\myexternaldocument{ex_supplement}

\usepackage{comment}
\usepackage{tikz}
\usetikzlibrary{arrows}
\usetikzlibrary{decorations.pathmorphing}
\usetikzlibrary{decorations.markings}
\usetikzlibrary{patterns}
\usetikzlibrary{automata}
\usetikzlibrary{positioning}
\usepackage{tikz-cd}

\tikzset{->-/.style={decoration={markings,mark=at position #1 with {\arrow[thick]{>}}},postaction={decorate}}}
	
\tikzset{-<-/.style={decoration={markings,mark=at position #1 with{\arrow[thick]{<}}},postaction={decorate}}}

\tikzset{
  frac arrow/.style={postaction={decorate,decoration={
        markings,
        mark=at position #1 with {\arrow[]{Stealth[round]}}
      }}},
}

\tikzset{
  mid arrow/.style={postaction={decorate,decoration={
        markings,
        mark=at position .5 with {\arrow[#1]{Stealth[round]}}
      }}},
}

\tikzset{
  midp arrow/.style={postaction={decorate,decoration={
        markings,
        mark=at position .6 with {\arrow[#1]{Stealth[round]}}
      }}},
}

\tikzset{
  midpp arrow/.style={postaction={decorate,decoration={
        markings,
        mark=at position .8 with {\arrow[#1]{Stealth[round]}}
      }}},
}
\tikzset{
  midm arrow/.style={postaction={decorate,decoration={
        markings,
        mark=at position .4 with {\arrow[#1]{Stealth[round]}}
      }}},
}

\tikzset{
  midmm arrow/.style={postaction={decorate,decoration={
        markings,
        mark=at position .3 with {\arrow[#1]{Stealth[round]}}
      }}},
}

\newcommand\spiral{}
\def\spiral[#1](#2)(#3:#4:#5){
\pgfmathsetmacro{\domain}{pi*#3/180+#4*2*pi}
\draw [#1,
       shift={(#2)},
       domain=0:\domain,
       variable=\t,
       smooth,
       samples=int(\domain/0.08)] plot ({\t r}: {#5*\t/\domain})
}

\newcommand\bonusspiral{}
\def\bonusspiral[#1](#2)(#3:#4)(#5:#6)[#7]{
\pgfmathsetmacro{\domain}{#4+#7*360}
\pgfmathsetmacro{\growth}{180*(#6-#5)/(pi*(\domain-#3))}
\draw [#1,
       shift={(#2)},
       domain=#3*pi/180:\domain*pi/180,
       variable=\t,
       smooth,
       samples=int(\domain/5)] plot ({\t r}: {#5+\growth*\t-\growth*#3*pi/180})
}

\DeclareRobustCommand{\sourcenode}{\tikz[baseline=-3pt]{\node at (0,0) [thick,black,circle,fill=lightgray,draw=black,inner sep=2pt]{\tiny $+$}}}

\DeclareRobustCommand{\sinknode}{\tikz[baseline=-3pt]{\node at (0,0) [thick,black,circle,fill=lightgray,draw=black,inner sep=2pt]{\tiny $-$}}}

\DeclareRobustCommand{\inftynode}{
    \tikz[baseline=-3pt]{\node at (0,0) [thick,black,circle,fill=lightgray,draw=black,inner sep=2pt]{\tiny $\infty$}}}

\DeclareRobustCommand{\saddlenode}{\tikz[baseline=-3pt]{\node at (0,0) [circle,fill=black,minimum size=8pt,inner sep=0pt]{}}}

\begin{document}

\maketitle

\begin{abstract}
  The leading-order asymptotic behavior of the solution of the Cauchy initial-value problem for the Benjamin-Ono equation in $L^2(\mathbb{R})$ is obtained explicitly for generic rational initial data  $u_0$.  
  An explicit asymptotic wave profile $u^\mathrm{ZD}(t,x;\epsilon)$ is given, in terms of the branches of the multivalued solution of the inviscid Burgers equation with initial data $u_0$, such that the solution $u(t,x;\epsilon)$ of the Benjamin-Ono equation with dispersion parameter $\epsilon>0$ and initial data $u_0$ satisfies $u(t,x;\epsilon)-u^\mathrm{ZD}(t,x;\epsilon)\to 0$ in the locally uniform sense as $\epsilon\to 0$, provided a discriminant inequality holds implying that certain caustic curves in the $(t,x)$-plane are avoided.  In some cases this convergence implies strong $L^2(\mathbb{R})$ convergence.  The asymptotic profile $u^\mathrm{ZD}(t,x;\epsilon)$ is consistent with the modulated multi-phase wave solutions described by Dobrokhotov and Krichever.
\end{abstract}

\begin{keywords}
  Benjamin-Ono equation, zero-dispersion limit, steepest descent integration.
\end{keywords}

\begin{AMS}
35C20, 35Q51, 41A60.
\end{AMS}

\section{Introduction}
\subsection{Background:  weak dispersion and dispersive shock waves}
Weakly dispersive perturbations of the initial-value problem for the inviscid Burgers equation 
\begin{equation}
u^\mathrm{B}_t+2u^\mathrm{B}u^\mathrm{B}_x=0,\quad u^\mathrm{B}(0,x)=u_0(x),
\label{eq:inviscid-Burgers}
\end{equation}
have been studied for many years.  Their theory is complementary to that of weakly diffusive perturbations such as the well-known viscous Burgers equation $u_t+2uu_x-\epsilon u_{xx}=0$, $\epsilon>0$, that are understood in some generality to lead to viscous shock waves converging as $\epsilon\to 0$ in a suitable weak topology to weak solutions of~\eqref{eq:inviscid-Burgers} satisfying appropriate entropy conditions.  See~\cite[Section 8]{Bressan13} for an accessible review.  When the perturbing terms are dispersive rather than diffusive, the perturbed system reacts to exploding gradients, a phenomenon occurring generally in solutions of~\eqref{eq:inviscid-Burgers} in finite time, by generating rapid oscillations rather than a sharp shock profile modeled as $\epsilon\to 0$ by a traveling jump discontinuity.  The field of rapid oscillations that resolve the gradient catastrophe in the solution of~\eqref{eq:inviscid-Burgers} is called a \emph{dispersive shock wave}.  

At a formal level, the oscillations making up the dispersive shock wave can be studied using the wave modulation theory of Whitham~\cite{Whitham}.  Whitham's theory starts with a family of exact periodic traveling wave solutions of the dispersive model indexed by constant parameters such as the wavenumber, frequency, and amplitude of the wave.  Subject to an ansatz of a solution fitting this form but with constant parameters replaced by slowly-varying functions instead, the theory predicts a system of \emph{modulation equations} governing the wave parameters.  There are several approaches to deriving the modulation equations going back to the early work of Whitham; for instance one can average the densities and fluxes in a sufficient number of local conservation laws~\cite{Whitham-ConservationLaws}, or equivalently one can average the density of a Lagrangian~\cite{Whitham-Variational} and obtain the modulation equations from an averaged variational principle.  This theory requires rather few assumptions about the dispersive perturbation of~\eqref{eq:inviscid-Burgers} and applies also to weakly-dispersive perturbations of more complicated quasilinear systems.  See~\cite{El-Hoefer16} and~\cite{BenzoniMietkaRodrigues21} for an up-to-date review.

While there is no specific reason to doubt its validity, the Whitham theory has only been rigorously established in isolated cases, usually with proofs relying on complete integrability.  The most famous completely integrable dispersive perturbation of the inviscid Burgers initial-value problem~\eqref{eq:inviscid-Burgers} is that of the Korteweg-de Vries equation 
\begin{equation}
u_t+2uu_x +\epsilon^2u_{xxx}=0.
\label{eq:KdV}
\end{equation}
In some sense, it was the formation of dispersive shock waves that led to the discovery of complete integrability for certain nonlinear dispersive wave equations.  Indeed, in the famous paper of Zabusky and Kruskal~\cite{ZabuskyKruskal}, the Korteweg-de Vries initial-value problem was solved numerically with sinusoidal initial data and a small value of $\epsilon$, and the individual peaks of the dispersive shock wave were identified for the first time as ``solitons''.  There were also hints that the Whitham theory for the Korteweg-de Vries equation had special properties, going back to the nontrivial observation of Whitham~\cite{Whitham-ConservationLaws} that the $3\times 3$ quasilinear first-order system of modulation equations was not only hyperbolic but also admitted a complete set of independent Riemann invariants.  Also, the periodic traveling wave solutions on which the Whitham theory is based were generalized in the Korteweg-de Vries setting to multi-phase wave solutions, for which modulation equations were also obtained and shown to have the same property for reasons related to complete integrability~\cite{FlaschkaFM80}.  The rigorous mathematical theory of the zero-dispersion limit $\epsilon\to 0$ of the Korteweg-de Vries equation began with the work of Lax and Levermore~\cite{LL83}, who established the existence of an $\epsilon$-independent limit $\bar{u}(t,x)$ of the solution $u(t,x;\epsilon)$ in the weak $L^2(\mathbb{R})$ sense for each $t\ge 0$.  The Lax-Levermore theory describes $\bar{u}(t,x)$ as the solution of a limiting variational problem for measures that are weak-$\ast$ limits of discrete measures; subsequently Venakides~\cite{Venakides90} showed how the retention of the discrete nature of the measure allows one to capture the profile of the $\mathcal{O}(\epsilon)$-wavelength oscillations of the dispersive shock wave.  Venakides' theory was not able to fully resolve the phase of the dispersive shock wave profile, but that was eventually achieved using the Riemann-Hilbert formulation of the inverse-scattering problem for the Korteweg-de Vries equation and applying the Deift-Zhou steepest descent technique~\cite{DeiftVZ97}. Further developments towards universality properties near wave caustics were achieved by Claeys and Grava in~\cite{ClaeysG09,ClaeysG10a, ClaeysG10b}.

\subsection{The Benjamin-Ono equation}
Another weakly-dispersive correction to the inviscid Burgers equation is the Benjamin-Ono equation
\begin{equation}
u_t + (u^2)_x = \epsilon \partial_x|D_x|u,\quad x\in\mathbb{R},\quad t\ge 0,\quad \epsilon>0,
\label{eq:BO}
\end{equation}
where the operator $|D_x|$ can be defined via a Fourier multiplier as $\widehat{|D_x|u}(\xi)=|\xi|\widehat{u}(\xi)$.
Like the Korteweg-de Vries equation, equation~\eqref{eq:BO} can be derived as an asymptotic model for long water waves of small amplitude, although for~\eqref{eq:BO} these are internal waves on a free surface between two density-stratified fluids with the lower layer being infinitely deep.

Like the Korteweg-de Vries equation, the Benjamin-Ono equation is understood to be a completely integrable system, although many aspects of the corresponding theory are quite different in character.  Here we restrict the discussion to the aspects that are germane to the zero-dispersion limit $\epsilon\to 0$ for~\eqref{eq:BO} posed with given initial data $u(0,x;\epsilon)=u_0(x)$ independent of $\epsilon$.  

Firstly, the Benjamin-Ono equation admits a family of periodic traveling wave solutions first obtained by Benjamin~\cite{Benjamin67} and Ono~\cite{Ono75}.  Later, all smooth periodic traveling wave solutions for the Benjamin-Ono equation in the form~\eqref{eq:BO} were systematically classified by Amick and Toland~\cite{AmickToland1991}, implying the completeness of the earlier results of~\cite{Benjamin67,Ono75}. These waves correspond to the initial data $u_0(x)=M U_r(M\epsilon^{-1}x+\alpha)+a$ where $M>0$ (taking $M\in\epsilon \mathbb{N}$ gives the periodic solutions with spatial period $2\pi$), $\alpha\in\mathbb{T}:=\mathbb{R}\pmod{2\pi}$, $a\in \mathbb{R}$ and for some a parameter $0<r<1$,
\begin{equation}\label{eq:Ur}
U_r(\theta):=\frac{1-r^2}{1+r^2-2r\cos(\theta)}.
\end{equation}
The traveling wave solution itself is then $u(t,x;\epsilon)=u_0(x-c_r t)$ where the 
velocity is $c_r=M(1+r^2)/(1-r^2)$, the maximum amplitude is $a+M(1+r)/(1-r)$ and the peak-to-trough variation is $4Mr/(1-r^2)$.  These results did not require integrability, but the generalization to multi-phase wave solutions of Benjamin-Ono was achieved in a paper of Dobrokhotov and Krichever~\cite{DobrokhotovK91} making full use of the Lax pair structure of~\eqref{eq:BO}.

\begin{theorem}[Multi-phase solutions of Dobrokhotov, Krichever~\cite{DobrokhotovK91}]\label{thm:J-phase}
The $J$-phase solutions are characterized by real constant parameters $R_0<R_1<R_2<\cdots<R_{2J}$, and nonzero complex parameters $\gamma_1,\dots,\gamma_J$ satisfying
\begin{equation}
    |\gamma_j|^2 = -\frac{(R_{2j}-R_0)\displaystyle\mathop{\prod_{k=1}^J}_{k\neq j}(R_{2j-1}-R_{2k-1})\cdot\mathop{\prod_{i=1}^J}_{i\neq j}(R_{2j}-R_{2i})}{(R_{2j-1}-R_0)\displaystyle\mathop{\prod_{k=1}^J(R_{2j}-R_{2k-1})\cdot \prod_{i=1}^J(R_{2j-1}-R_{2i})}},\quad j=1,\dots,J.
\label{eq:mod-gamma-exact}
\end{equation}
They are given by
\begin{equation}
    u(t,x;\epsilon)
    =	R_0+\sum_{j=1}^J(R_{2j-1}-R_{2j})-2\epsilon\mathrm{Im}\left(\frac{\partial}{\partial x}\log(\det(\mathbf{M}(t,x;\epsilon)))\right)
\label{eq:DK-exact}
\end{equation}
in which $\mathbf{M}(t,x;\epsilon)$ is a $J\times J$ matrix with elements
\begin{equation}
    M_{jk}:=\gamma_j \ee^{\ii\theta^\mathrm{L}_j(t,x)/\epsilon} \delta_{jk}+\frac{1}{R_{2j-1}-R_{2k}},\quad 1\le j,k\le J
\label{eq:DK-exact-M}
\end{equation}
and $\theta^\mathrm{L}_j(t,x)$ is a linear phase given by
\begin{equation}
    \theta^\mathrm{L}_j(t,x):=(R_{2j-1}-R_{2j})x - (R_{2j-1}^2-R_{2j}^2)t,\quad j=1,\dots,J.
    \label{eq:linear-phase}
\end{equation}
\end{theorem}
In the same paper~\cite{DobrokhotovK91}, the Benjamin-Ono analogue of the multiphase modulation theory developed for Korteweg-de Vries in~\cite{FlaschkaFM80} was advanced; in this theory the parameters $R_j$ are replaced with functions of $(t,x)$ that are all required to satisfy individually the same equation, namely~\eqref{eq:inviscid-Burgers}.  It is worth pointing out that the exact multi-phase solutions for Benjamin-Ono are rational trigonometric functions and hence are far simpler than the corresponding Korteweg-de Vries solutions which are built instead from Riemann theta functions of hyperelliptic curves of genus $J$.  Similarly, the multi-phase modulation theory for Benjamin-Ono is far simpler than that of the Korteweg-de Vries equation which takes the form of a \emph{coupled} system of $2J+1$ equations in Riemann invariant form~\cite{FlaschkaFM80}.  In \cite{Moll20}, a correspondence was established between $J$-phase solutions and finite-gap potentials for which only the first $J$ Birkhoff coordinates may be nonzero. See~\cite[Section 7]{GerardKappeler2019} for a formula similar to~\eqref{eq:DK-exact} by using a $J\times J$ matrix defined from the first $J$ Birkhoff coordinates of the initial data. 

\begin{figure}
\begin{center}
    \includegraphics[width=0.48\linewidth]{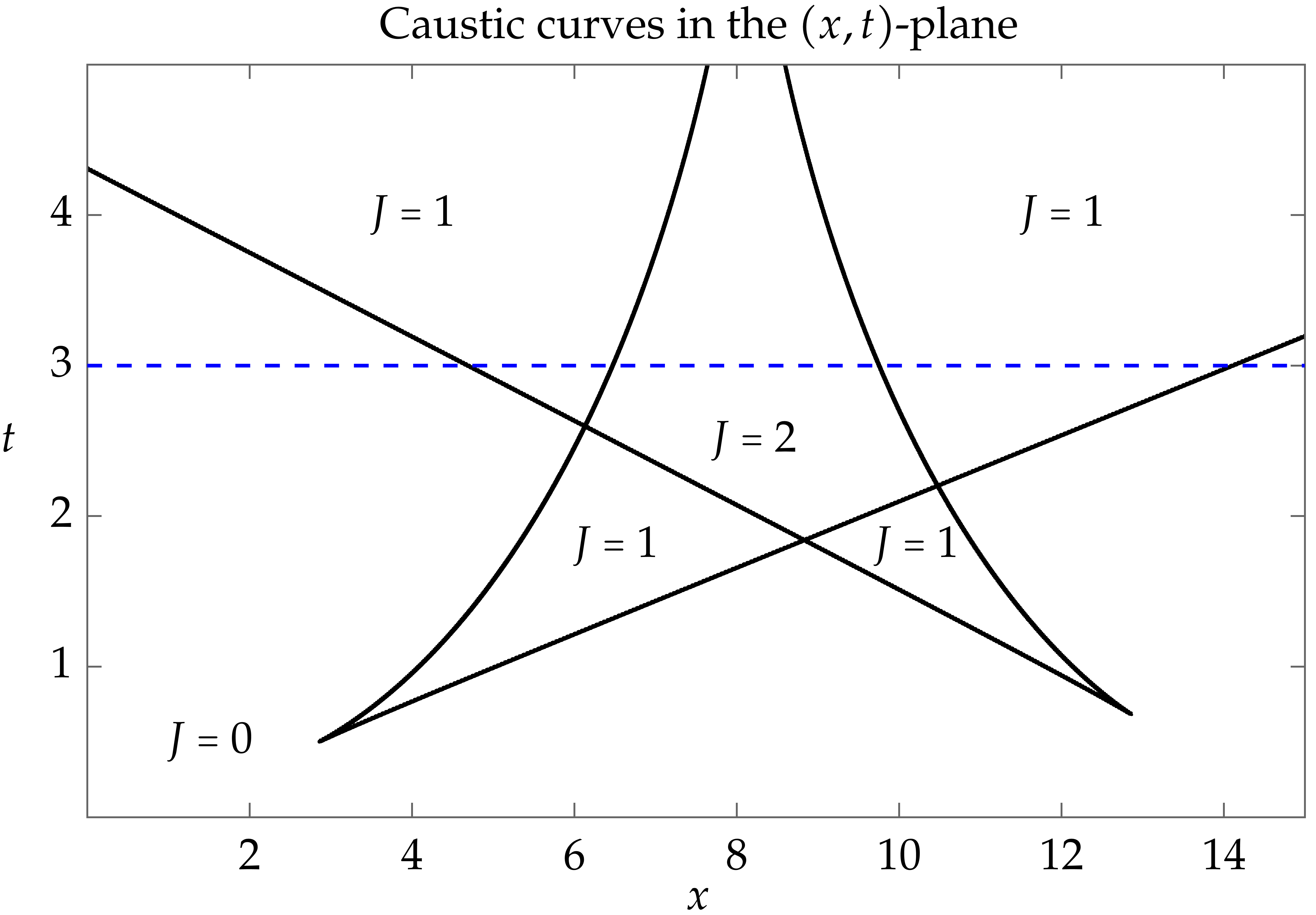}\hfill%
    \includegraphics[width=0.48\linewidth]{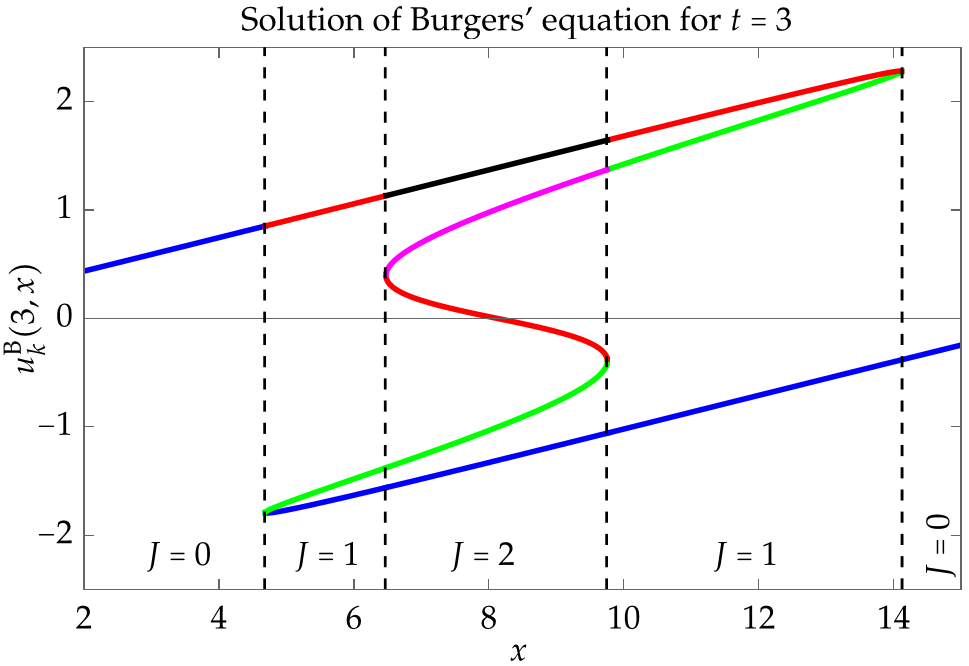}
\end{center}
\caption{The plots in this figure correspond to the initial condition $u_0(x)$ of the form \eqref{eq:rationalIC} for $N=2$ with $c_1=1-\ii$, $c_2=1+\ii/\sqrt{2}$, $p_1=\ii$, $p_2=16+\ii$.  Left:  the caustic curves dividing the plane into regions where the solution of Burgers' equation is $2J+1=2J(t,x)+1$-valued, as shown. Right: the evolution of $u_0(x)$ under Burgers' equation at $t=3$.  The blue, green, red, purple, black curves are $u^{\mathrm{B}}_k(3,x)$, $k=0,\ldots,4$, respectively.  The dashed black lines in the right pane are placed at the $x$ values where the dashed blue line intersects the caustic curves in the left pane.}
\label{fig: Burgers}
\end{figure}

A close analogue of the method of Lax and Levermore was developed for the Benjamin-Ono equation in~\cite{MillerX10}.  The weak zero-dispersion limit identified in~\cite{MillerX10} is constructed as follows.  Firstly, one solves the initial-value problem~\eqref{eq:inviscid-Burgers} by the method of characteristics, which produces a solution that one can write in implicit form as $u^\mathrm{B}=u_0(x-2u^\mathrm{B}t)$.  This solution is classical for small $t$, but for larger $t$ the implicit form undergoes a bifurcation that generates multiple solutions.  For generic coordinates $(t,x)\in\mathbb{R}^2$, there is a nonnegative integer $J=J(t,x)$ such that there are $2J+1$ distinct values (``branches'' or ``sheets'') of $u^\mathrm{B}(t,x)$ denoted $u_k^\mathrm{B}(t,x)$ that can be ordered as 
\begin{equation}
    u^\mathrm{B}_0(t,x)<u^\mathrm{B}_1(t,x)<\cdots<u^\mathrm{B}_{2J}(t,x).
\label{eq:uB-branches-ordered}
\end{equation}
Each sheet defines a locally smooth solution of the same equation~\eqref{eq:inviscid-Burgers}, see Figure~\ref{fig: Burgers} or~\cite[Fig.~1]{BGM23} for a simpler example. Then the result of~\cite{MillerX10} is that as $\epsilon\to 0$, the solution of the initial-value problem for the Benjamin-Ono equation~\eqref{eq:BO} with suitable initial data approximating $u_0$ in the strong $L^2(\mathbb{R})$ sense converges weakly in $L^2(\mathbb{R})$ for each fixed $t\ge 0$ to $\bar{u}(t,x)$ defined explicitly by
\begin{equation}
    \bar{u}(t,x):=\sum_{k=0}^{2J(t,x)}(-1)^ku_k^\mathrm{B}(t,x)
    \label{eq:ubar}
\end{equation}
in the limit $\epsilon\to 0$.  This result is far simpler to state than the corresponding result for the Korteweg-de Vries equation~\cite{LL83} which requires the solution of a quadratic variational problem for each $(t,x)\in\mathbb{R}^2$.  More recently, the same formula~\eqref{eq:ubar} was proven to yield the weak zero-dispersion limit for Benjamin-Ono posed with periodic initial data~\cite{Gassot23}, without resorting to an approximation of the initial data~\cite{Gassot23b}.  

An attempt to strengthen the zero-dispersion asymptotics for the Benjamin-Ono equation was given in~\cite{BGM23}, following along the lines of what Venakides accomplished for the Korteweg-de Vries equation~\cite{Venakides90}.  However, to date there has not been complete success in fully capturing the dispersive shock waves generated from smooth initial data for the Benjamin-Ono equation as $\epsilon\to 0$.  

Some recent progress for the Benjamin-Ono equation~\eqref{eq:BO} came in the paper~\cite{Gerard22}, which presented for the first time a closed-form expression for the solution with general initial data $u_0\in L^2(\mathbb{R})\cap L^{\infty}(\mathbb{R})$. That formula was quickly used~\cite{Gerard23} to prove that the weak convergence result of~\cite{MillerX10} holds more generally, and the proof was also much more efficient. The inversion formula for the Benjamin-Ono equation was then extended to initial data  $u_0\in L^2(\mathbb{R})$, and a new proof of the weak convergence result of the solutions was given at low regularity in~\cite{Chen24}. A similar strategy was implemented for the Calogero-Moser derivative NLS equation on the line, where in~\cite{Badreddine24zero} the author  establishes a weak convergence result in the zero-dispersion limit from an explicit inversion formula~\cite{KillipLaurensVisan23}. 

However, strong zero-dispersion asymptotics remained out of reach.  The latest development is that the solution formula of~\cite{Gerard22} was effectively implemented in~\cite{RationalBO} for rational initial data of the form
\begin{equation}
u_0(x):=\sum_{n=1}^N\left(\frac{c_n}{x-p_n}+\frac{{c_n^*}}{x-{p_n^*}}\right),
\quad \operatorname{Im}(p_n)>0, \forall n.
\label{eq:rationalIC}
\end{equation}
Note that $u_0\in L^2(\mathbb{R})$.  The solution $u(t,x;\epsilon)$ of~\eqref{eq:BO} for such initial data is expressed in~\cite{RationalBO} as a ratio of $(N+1)\times (N+1)$ determinants whose entries are certain contour integrals depending on $(t,x;\epsilon)$.  Using these determinantal formul\ae\ it is now possible to prove a rigorous result about the oscillatory profile of dispersive shock waves generated by the Benjamin-Ono equation.

\subsection{Main result: strong zero-dispersion asymptotics}\label{sec: main result}
In solving~\eqref{eq:inviscid-Burgers} by the method of characteristics, one first seeks for given $(t,x)\in\mathbb{R}^2$ the intercepts $y$ of characteristic lines through $(t,x)$ by solving the equation
\begin{equation}\label{eq:critical-points}
y-x+2t u_0(y)=0.
\end{equation}
For rational $u_0$ of the form~\eqref{eq:rationalIC}, this is equivalent to a polynomial equation of degree $2N+1$ for $y$ with real coefficients that are polynomial in $(t,x)$.  Let $\mathcal{C}\subset\mathbb{R}^2$ denote the \emph{caustic locus} of points $(t,x)$ for which there is a real root $y\in\mathbb{R}$ of this polynomial having multiplicity greater than $1$, this set is of Lebesgue measure zero, see~\cite{Gerard23}.  For given $(t,x)\in\mathbb{R}^2\setminus \mathcal{C}$, the equation~\eqref{eq:critical-points}  has a positive odd number $2J(t,x)+1$ of simple real solutions labeled as
\begin{equation}
 y_0(t,x)>y_1(t,x)>\dots > y_{2J(t,x)}(t,x).
\end{equation}
These are all analytic functions of $(t,x)$ on each maximal connected component of $\mathbb{R}^2\setminus \mathcal{C}$, and $J(t,x)\ge 0$ is an integer-valued function that is constant on the indicated component.  The corresponding branches of the (multivalued, if $J(t,x)>0$) solution of~\eqref{eq:inviscid-Burgers} with initial data $u_0$ are then given by 
\begin{equation}\label{eq:uB-yj}
u_k^{\mathrm{B}}(t,x)=u_0(y_k(t,x)),\quad k=0,\dots,2J(t,x),
\end{equation}
and they are also the analytic functions of $(t,x)$ satisfying the strict inequalities~\eqref{eq:uB-branches-ordered}.  Combining this with~\eqref{eq:critical-points} shows that $y_k(t,x)=x-2tu^\mathrm{B}_k(t,x)$ holds for $k=0,\dots,2J(t,x)$.

We now define a function $h:\mathbb{R}\to\mathbb{R}$ by
\begin{equation}
h(y)=h(y;t,x):=\frac{1}{4t}(y-x)^2+\sum_{n=1}^N\left[c_n\log(
y-p_n)+c_n^*\log(y-p_n^*)\right],
\label{eq:h-def}
\end{equation}
where the complex logarithms denote the principal branches.  Note that by using~\eqref{eq:rationalIC} we see that
\begin{equation}
    h'(y)=\frac{y-x}{2t}+u_0(y),
 \label{eq:hprime}
\end{equation}
so comparing with~\eqref{eq:critical-points}, the real critical points of $h(y)$ are precisely the intercepts of characteristic lines through $(t,x)$.
For those familiar with the theory of entropy solutions of the inviscid Burgers equation~\eqref{eq:inviscid-Burgers}, the quantity $h(y)$ is precisely the objective function that one minimizes over $y$ in the Lax-Oleinik formula~\cite{Evans} to obtain the correct weak solution.  However, an important point is that for the Benjamin-Ono equation as a dispersive rather than viscous regularization of the inviscid Burgers equation~\eqref{eq:inviscid-Burgers}, \emph{all} of the real critical points of $h(y)$ will play an equal role rather than just the minimizer.  

Now, we define a function $u^\mathrm{ZD}(t,x;\epsilon)$ for $(t,x)\in\mathbb{R}^2\setminus \mathcal{C}$ and $\epsilon>0$ in terms of the above data. To define $u^\mathrm{ZD}(t,x;\epsilon)$ when $J(t,x)>0$, we introduce some nonlinear phases, nonlinear phase corrections and modulation parameters.

First, we use~\eqref{eq:h-def} to define nonlinear phases by 
\begin{equation}
    \theta_j(t,x):=h(y_{2j-1}(t,x))-h(y_{2j}(t,x)), \quad j=1,\dots,J(t,x).
    \label{eq:theta-j}
\end{equation}
Then, we define nonlinear phase corrections. For this, we first note that the rational function defined by 
\begin{equation}\label{NEWfactorize-bis}
g(y):=\frac{y+2tu_0(y)-x}{\displaystyle{\prod_{k=0}^{2J(t,x)}(y+2tu_k^\mathrm{B}(t,x)-x)}}
\end{equation}
is real-valued and analytic for $y\in\mathbb{R}$, that $g(y)>0$ holds on $\mathbb{R}$, and also that $g(y)=y^{-2J(t,x)}(1+\mathcal{O}(y^{-1}))$ as $y\to\infty$.  Indeed, the denominator simply consists of the factors from the numerator that vanish on the real line.  From these properties of $g(y)$ it is easy to check that 
\begin{equation}\label{Phi}
\Phi (y):=-J\frac{\pi}{2}+\frac{1}{2\pi}\int_0^{+\infty}\ln\left(\frac{g(y-s)}{g(y+s)}\right)\, \frac{\dd s}{s},\quad y\in\mathbb{R}
\end{equation}
is a well-defined real-valued function of $y\in\mathbb{R}$. We define phase corrections with the help of the function $\Phi(y)$ defined in~\eqref{Phi} by 
\begin{equation}
    \varphi_j(t,x):=\frac{\pi}{2}+\Phi(y_{2j-1}(t,x))-\Phi(y_{2j}(t,x)),\quad j=1,\dots,J(t,x).
    \label{eq:varphi}
\end{equation}
Finally, we define 
\begin{equation}
R_j(t,x)
	:=u^B_j(t,x), \quad j=0,\dots,J(t,x)
\label{eq:V}
\end{equation}
and introduce complex functions $\gamma_j(t,x)$, $j=1,\dots,J(t,x)$ by
\begin{equation}
    \gamma_j(t,x):=|\gamma_j(t,x)|\ee^{\ii\varphi_j(t,x)},\quad j=1,\dots,J(t,x),
\label{eq:gamma}
\end{equation}
where $|\gamma_j(t,x)|^2$ is defined by the condition~\eqref{eq:mod-gamma-exact} from Theorem~\ref{thm:J-phase} and $\varphi_j$ was introduced in~\eqref{eq:varphi}.

\begin{figure}
    \centering
    \includegraphics[width=0.7\linewidth]{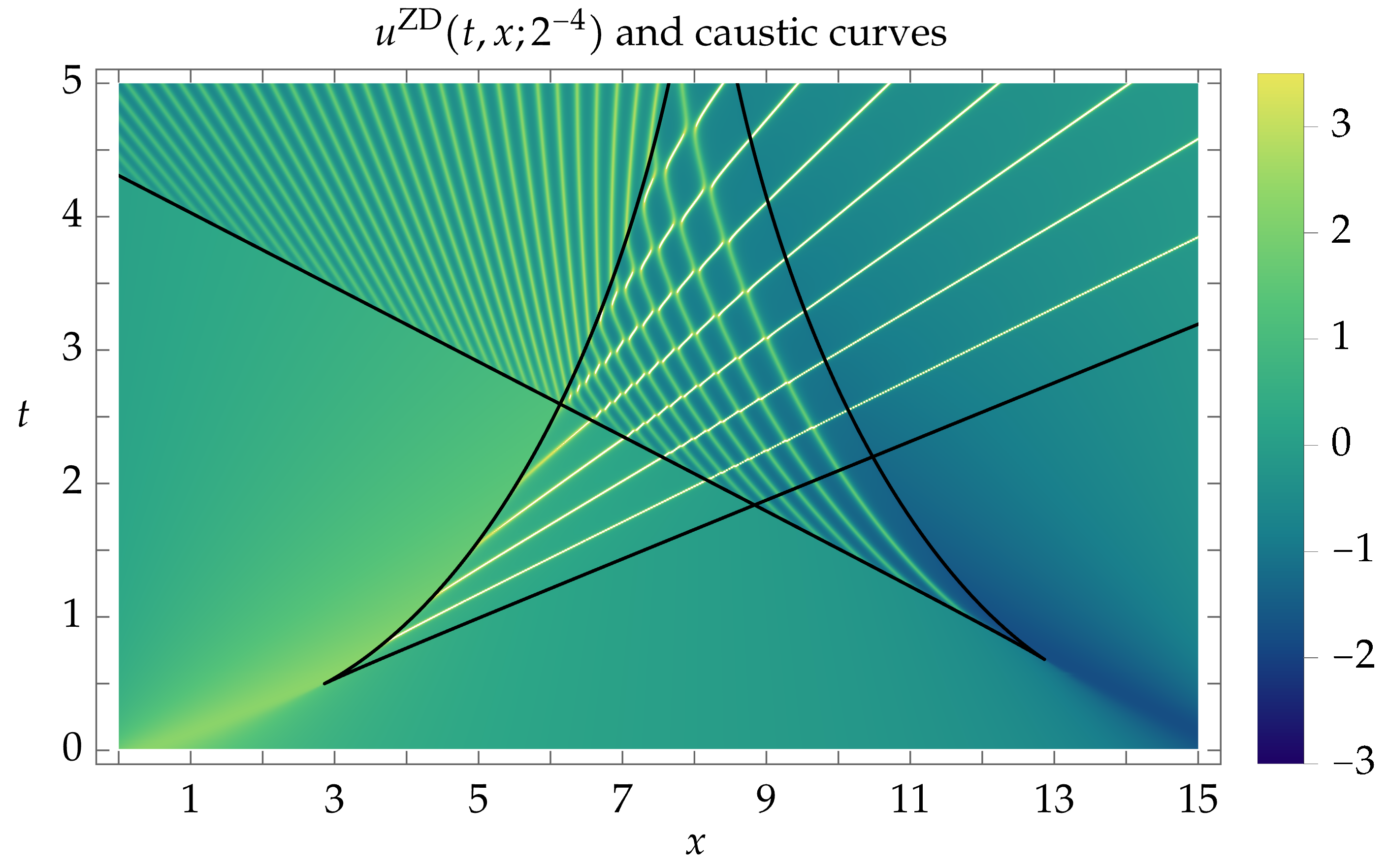}
    \caption{The zero-dispersion profile $u^{\mathrm{ZD}}(t,x;2^{-4})$ corresponding to the same $u_0(x)$ as in Figure \ref{fig: Burgers}.  The caustic curves are displayed in black.}
    \label{fig:uZD}
\end{figure}

\begin{definition}[Zero-dispersion approximation]\label{def:zero}
If $J(t,x)=0$ then we set:
\begin{equation}
    u^\mathrm{ZD}(t,x;\epsilon):=u_0^\mathrm{B}(t,x).
\end{equation}  
If $J(t,x)>0$, we define $u^\mathrm{ZD}(t,x;\epsilon)
$ as a $J$-phase profile similarly as in Theorem~\ref{thm:J-phase} with $\theta_j$, $R_j$ and $\gamma_j$ given in~\eqref{eq:theta-j}, \eqref{eq:V} and~\eqref{eq:gamma} respectively:
\begin{equation}\label{eq:uZD}
    u^\mathrm{ZD}(t,x;\epsilon)
    	:=u_0^\mathrm{B}(t,x)+\sum_{j=1}^{J(t,x)}(u_{2j-1}^\mathrm{B}(t,x)-u_{2j}^\mathrm{B}(t,x))-2\epsilon\mathrm{Im}\left(\frac{\partial}{\partial x}\log(\det(\mathbf{M}(t,x;\epsilon)))\right),
\end{equation}
\begin{equation}
    M_{jk}(t,x;\epsilon):=\gamma_j(t,x) \ee^{\ii\theta_j(t,x)/\epsilon}\delta_{jk} +\frac{1}{u^\mathrm{B}_{2j-1}(t,x)-u^\mathrm{B}_{2k}(t,x)},\quad 1\le j,k\le J(t,x).
    \label{eq:Mjk}
\end{equation}
\end{definition}
See Figure \ref{fig:uZD} for a plot of $u^{\mathrm{ZD}}(t,x,\epsilon)$.

\begin{remark}
If $J(t,x)=1$, this formula simplifies as follows:
\begin{equation}
    u^\mathrm{ZD}(t,x;\epsilon)=u_0^\mathrm{B}(t,x)+(u_2^\mathrm{B}(t,x)-u_1^\mathrm{B}(t,x))U_{r(t,x)}\left(\frac{\theta_1(t,x)}{\epsilon}+\varphi_1(t,x)\right),
    \label{eq:uZD-J1}
\end{equation}
where $U_r(\theta)$ is given by~\eqref{eq:Ur} and the parameter $r(t,x)$ is given by 
\begin{equation}\label{eq:r}
0<r(t,x):=\sqrt{\frac{u^\mathrm{B}_1(t,x)-u^\mathrm{B}_0(t,x)}{u^\mathrm{B}_2(t,x)-u^\mathrm{B}_0(t,x)}}<1.
\end{equation}
\end{remark}
\begin{remark}[Local wavenumber and frequency]
Note that according to~\eqref{eq:wavenumbers} and~\eqref{eq:frequencies} below, the local wavenumber $\kappa_j(t,x)$ and frequency $\omega_j(t,x)$ for the nonlinear phase are given respectively by 
\begin{equation}
\begin{split}
    \kappa_j(t,x):=&\frac{\partial\theta_j}{\partial x}(t,x)=u^\mathrm{B}_{2j-1}(t,x)-u_{2j}^\mathrm{B}(t,x), \\
    \omega_j(t,x):=&-\frac{\partial\theta_j}{\partial t}(t,x)=u_{2j-1}^\mathrm{B}(t,x)^2-u_{2j}^\mathrm{B}(t,x)^2,
\end{split}
\end{equation}
which match the corresponding expressions~\eqref{eq:linear-phase} for the derivatives of the linear phase $\theta_j^\mathrm{L}(t,x)$ replacing $U_k$ with $u_k^\mathrm{B}(t,x)$.
With this identification, the formul\ae\ \eqref{eq:Mjk} are seen to match~\eqref{eq:DK-exact-M}, if only we generalize from the linear phase $\theta_j^\mathrm{L}(t,x)$ to the nonlinear phase $\theta_j(t,x)$ whose derivatives reproduce the coefficients of $x$ and $t$ locally.
\label{rem:local-wavenumber}
\end{remark}

Now we may formulate our main result.  Note that for rational $u_0$, the caustic locus $\mathcal{C}\subset\mathbb{R}^2$ is a subset of the \emph{discriminant locus} $\mathcal{D}\subset\mathbb{R}^2$ consisting of points $(t,x)$ for which the characteristic equation $x=2tu_0(y)+y$ has a solution, real or complex, of multiplicity greater than $1$.  Like $\mathcal{C}$, $\mathcal{D}$ also has Lebesgue measure zero.  In fact, for each fixed $t\in\mathbb{R}$, $\mathcal{D}$ consists of at most $4N$ values of $x\in\mathbb{R}$.

\begin{theorem}[Strong zero-dispersion asymptotics]\label{thm:u-app}
Assume that $u_0\in L^2(\mathbb{R})$ is a rational initial condition of the form \eqref{eq:rationalIC}, and suppose that $t>0$ and that $(t,x)\in\mathbb{R}^2\setminus \mathcal{D}$. Then, 
\begin{equation}\label{eq:u-app}
u(t,x;\epsilon)=u^\mathrm{ZD}(t,x;\varepsilon)+\mathcal{O}({\varepsilon})
\end{equation}
where the convergence is pointwise for such $(t,x)$ and also uniform on compact subsets of 
$\mathbb{R}^2\setminus \mathcal{D}$.
\end{theorem}

\begin{corollary}[Convergence in $L^2$]\label{cor:CV-L2}
Let $u_0\in L^2(\mathbb{R})$ be a rational initial condition of the form \eqref{eq:rationalIC}.  If $t\ge 0$ is such that  $J(t,x)\leq 1$ for every generic $x$ where $J(t,x)$ is well-defined, then as $\varepsilon\to 0$, 
\begin{equation}
u(t,\diamond;\epsilon)=u^\mathrm{ZD}(t,\diamond;\epsilon)+o_{L^2}(1).
\label{eq:L2-convergence}
\end{equation}
\end{corollary}

\begin{remark}
It is reasonable to expect that \eqref{eq:L2-convergence} holds even at times $t\ge 0$ for which $J(t,x)$ has any finite upper bound.  However, the proof of Corollary~\ref{cor:CV-L2} involves an averaging argument (see Lemma~\ref{lem:averaging} below) that would need to be upgraded to a multiphase averaging argument if $J(t,x)\ge 2$.  Because such arguments are difficult to explain (Diophantine resonance conditions must be dealt with carefully), to keep the paper as simple as possible we leave the details to future work.
\end{remark}

\begin{remark}
Note that the zero-dispersion profile $u^{\mathrm{ZD}}(t,x;\epsilon)$ still makes sense for non-rational initial data, however, the proof of Theorem~\ref{thm:u-app} is open without the rational assumption. Indeed, using a rational approximation of the initial data is not sufficient, because known well-posedness results in the literature (see, for instance \cite{KillipLaurensVisan23-BO}) do not provide estimates on the modulus of continuity of the data-to-solution map that are uniform in the limit $\epsilon\downarrow 0$.  One may also notice that the approximation formula itself is very sensitive to small errors.  For instance, looking at the simpler $J(t,x)=1$ version of the formula given in \eqref{eq:uZD-J1}, one can see that even if the modulation variables $u_j^\mathrm{B}(t,x)$, $j=0,1,2$, $\theta_1(t,x)$, and $\varphi_1(t,x)$ are perturbed by terms of size proportional to $\epsilon$, one can still obtain errors not decaying to zero as $\epsilon\to 0$ because of the factor $\epsilon^{-1}$ in the phase.  This shows that, for example, a spatial translation of $u_0(\diamond)$ by a shift of order $O(\epsilon)$ (which is an $O(\epsilon^{1/2})$ perturbation in $L^2(\mathbb{R})$) will lead to a leading-order phase shift in the approximation that will ruin the accuracy of the approximation.
\end{remark}

\subsection{Plan of the paper}

We first recall from~\cite{RationalBO} the representation of the exact solution of the Cauchy initial-value problem for the Benjamin-Ono equation~\eqref{eq:BO} with rational initial data of the form~\eqref{eq:rationalIC} in Section~\ref{sec:intro-inversion}. 
Then, in Section~\ref{section:N1}, we explain the general idea of proof in the case of one pole $N=1$, so that $J=0$ or $J=1$. 

In general, the analysis of the solution formula in the small-dispersion limit amounts to the (classical) steepest-descent analysis of a family of contour integrals, and for this approach to be effective it is first necessary to specify suitable contours of integration.  In Section~\ref{section:LevelSets}, we explain how to analyze the exponent function appearing in contour integrals in the solution formula by introducing the notion of extended level curves. 
Using this analysis, the specification of contours of integration is given in Section~\ref{section:StokesGraphs}. The reader willing to accept the result stated at the beginning of Section~\ref{section:StokesGraphs} (Proposition~\ref{prop:W}) on first reading can skip ahead to Section~\ref{sec:steepest-descent}, wherein the steepest descent method is applied and the resulting asymptotic formula is simplified, proving Theorem~\ref{thm:u-app}. We finally deduce Corollary~\ref{cor:CV-L2} from Theorem~\ref{thm:u-app} in Section~\ref{section:CV-L2}. Section \ref{sec:numerics} is a short demonstration of the validity of Theorem~\ref{thm:u-app} via numerics.

\subsection*{Acknowledgements}
E. Blackstone was partially supported by the National Science Fou\-ndation under grant DMS-1812625. L. Gassot was supported by the France 2030 framework program, the Centre Henri Lebesgue ANR-11-LABX-0020-01, and the ANR project HEAD--ANR-24-CE40-3260.
P. G\'erard was partially supported by the French Agence Nationale de la Recherche under the ANR project ISAAC--ANR-23--CE40-0015-01. P. D. Miller was partially supported by the National Science Foundation under grant DMS-2204896.  The authors would like to thank Matthew Mitchell for useful discussions.

\section{The solution of Benjamin-Ono with rational initial data}
\label{sec:intro-inversion}
Here we recall the result of~\cite{RationalBO}, in which a general solution formula for the Cauchy problem of the Benjamin-Ono equation~\eqref{eq:BO} with initial data in $L^2(\mathbb{R})$ was applied to rational initial data of the form~\eqref{eq:rationalIC}.

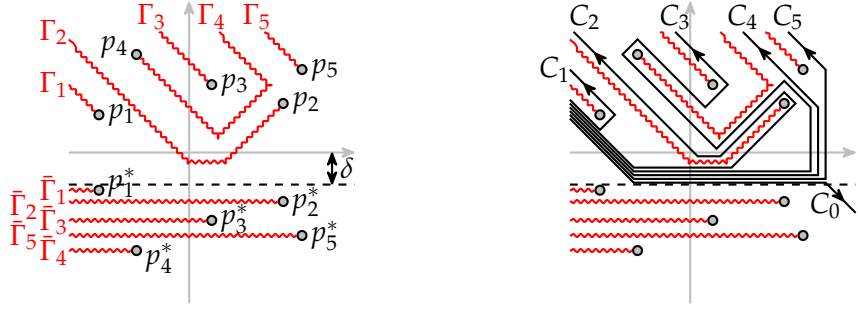
\begin{figure}
\hfill\begin{tikzpicture}
\node (p1) at (-1.2cm,0.5cm) [thick,black,circle,fill=lightgray,draw=black,inner sep=1.2pt]{};
\node (p2) at (1.25cm,0.65cm) [thick,black,circle,fill=lightgray,draw=black,inner sep=1.2pt]{}; 
\node (p3) at (0.3cm,0.9cm) [thick,black,circle,fill=lightgray,draw=black,inner sep=1.2pt]{}; 
\node (p4) at (-0.7cm,1.3cm) [thick,black,circle,fill=lightgray,draw=black,inner sep=1.2pt]{}; 
\node (p5) at (1.5cm,1.1cm) [thick,black,circle,fill=lightgray,draw=black,inner sep=1.2pt]{}; 
\node at (-0.9cm,0.5cm) {$p_1$};
\node at (1.55cm,0.65cm) {$p_2$};
\node at (0.6cm,0.9cm) {$p_3$};
\node at (-1.0cm,1.4cm) {$p_4$};
\node at (1.8cm,1.1cm) {$p_5$};
\node (p1s) at (-1.2cm,-0.5cm) [thick,black,circle,fill=lightgray,draw=black,inner sep=1.2pt]{};
\node (p2s) at (1.25cm,-0.65cm) [thick,black,circle,fill=lightgray,draw=black,inner sep=1.2pt]{}; 
\node (p3s) at (0.3cm,-0.9cm) [thick,black,circle,fill=lightgray,draw=black,inner sep=1.2pt]{}; 
\node (p4s) at (-0.7cm,-1.3cm) [thick,black,circle,fill=lightgray,draw=black,inner sep=1.2pt]{}; 
\node (p5s) at (1.5cm,-1.1cm) [thick,black,circle,fill=lightgray,draw=black,inner sep=1.2pt]{}; 
\node at (-0.9cm,-0.4cm) {$p_1^*$};
\node at (1.55cm,-0.65cm) {$p_2^*$};
\node at (0.6cm,-0.875cm) {$p_3^*$};
\node at (-0.4cm,-1.4cm) {$p_4^*$};
\node at (1.8cm,-1.1cm) {$p_5^*$};
\draw [->,>={Stealth[round]},thick,lightgray] (-1.6cm,0cm)--(2.2cm,0cm) {};
\draw [->,>={Stealth[round]},thick,lightgray] (0cm,-2cm) -- (0cm,2cm) {};
\draw [thick,red,decorate,
decoration={snake,amplitude=.2mm,segment length=1mm,post length=0mm}] (p1) -- (-1.6cm,0.9cm);
\draw [thick,red,decorate,
decoration={snake,amplitude=.2mm,segment length=1mm,post length=0mm}] (p5) -- (1cm,1.6cm);
\draw [thick,red,decorate,
decoration={snake,amplitude=.2mm,segment length=1mm,post length=0mm}] (p2) -- (0.5cm,-0.1cm) -- (0.0cm,-0.1cm) -- (-1.6cm,1.5cm);
\draw [thick,red,decorate,
decoration={snake,amplitude=.2mm,segment length=1mm,post length=0mm}] (p4) -- (0.4cm,0.2cm) -- (1.1cm,0.9cm) -- (0.4cm,1.6cm);
\draw [thick,red,decorate,
decoration={snake,amplitude=.2mm,segment length=1mm,post length=0mm}] (p3) -- (-0.4cm,1.6cm);
\draw [thick,red,decorate,
decoration={snake,amplitude=.2mm,segment length=1mm,post length=0mm}] (p1s) -- (-1.6cm,-0.5cm);
\draw [thick,red,decorate,
decoration={snake,amplitude=.2mm,segment length=1mm,post length=0mm}] (p2s) -- (-1.6cm,-0.65cm);
\draw [thick,red,decorate,
decoration={snake,amplitude=.2mm,segment length=1mm,post length=0mm}] (p3s) -- (-1.6cm,-0.9cm);
\draw [thick,red,decorate,
decoration={snake,amplitude=.2mm,segment length=1mm,post length=0mm}] (p4s) -- (-1.6cm,-1.3cm);
\draw [thick,red,decorate,
decoration={snake,amplitude=.2mm,segment length=1mm,post length=0mm}] (p5s) -- (-1.6cm,-1.1cm);

\node[red] at (-1.8cm,0.9cm) {$\Gamma_1$};
\node[red] at (-1.8cm,1.6cm) {$\Gamma_2$};
\node[red] at (-0.5cm,1.8cm) {$\Gamma_3$};
\node[red] at (0.3cm,1.8cm) {$\Gamma_4$};
\node[red] at (0.9cm,1.8cm) {$\Gamma_5$};

\node[red] at (-1.8cm,-0.5cm) {$\bar{\Gamma}_1$};
\node[red] at (-2.2cm,-0.7cm) {$\bar{\Gamma}_2$};
\node[red] at (-1.8cm,-0.9cm) {$\bar{\Gamma}_3$};
\node[red] at (-2.2cm,-1.1cm) {$\bar{\Gamma}_5$};
\node[red] at (-1.8cm,-1.3cm) {$\bar{\Gamma}_4$};

\draw [thick,black,dashed] (-1.6cm,-0.425cm) -- (2.2cm,-0.425cm);
\draw [thick,black,<->,>={Stealth[round]}] (1.9cm,0cm) -- (1.9cm,-0.425cm);
\node at (2.1cm,-0.2cm) {$\delta$};
\end{tikzpicture} \hfill%
\begin{tikzpicture}
\node (p1) at (-1.2cm,0.5cm) [thick,black,circle,fill=lightgray,draw=black,inner sep=1.2pt]{};
\node (p2) at (1.25cm,0.65cm) [thick,black,circle,fill=lightgray,draw=black,inner sep=1.2pt]{}; 
\node (p3) at (0.3cm,0.9cm) [thick,black,circle,fill=lightgray,draw=black,inner sep=1.2pt]{}; 
\node (p4) at (-0.7cm,1.3cm) [thick,black,circle,fill=lightgray,draw=black,inner sep=1.2pt]{}; 
\node (p5) at (1.5cm,1.1cm) [thick,black,circle,fill=lightgray,draw=black,inner sep=1.2pt]{}; 
\node (p1s) at (-1.2cm,-0.5cm) [thick,black,circle,fill=lightgray,draw=black,inner sep=1.2pt]{};
\node (p2s) at (1.25cm,-0.65cm) [thick,black,circle,fill=lightgray,draw=black,inner sep=1.2pt]{}; 
\node (p3s) at (0.3cm,-0.9cm) [thick,black,circle,fill=lightgray,draw=black,inner sep=1.2pt]{}; 
\node (p4s) at (-0.7cm,-1.3cm) [thick,black,circle,fill=lightgray,draw=black,inner sep=1.2pt]{}; 
\node (p5s) at (1.5cm,-1.1cm) [thick,black,circle,fill=lightgray,draw=black,inner sep=1.2pt]{}; 
\draw [->,>={Stealth[round]},thick,lightgray] (-1.6cm,0cm)--(2.2cm,0cm) {};
\draw [->,>={Stealth[round]},thick,lightgray] (0cm,-2cm) -- (0cm,2cm) {};
\draw [thick,red,decorate,
decoration={snake,amplitude=.2mm,segment length=1mm,post length=0mm}] (p1) -- (-1.6cm,0.9cm);
\draw [thick,red,decorate,
decoration={snake,amplitude=.2mm,segment length=1mm,post length=0mm}] (p5) -- (1cm,1.6cm);
\draw [thick,red,decorate,
decoration={snake,amplitude=.2mm,segment length=1mm,post length=0mm}] (p2) -- (0.5cm,-0.1cm) -- (0.0cm,-0.1cm) -- (-1.6cm,1.5cm);
\draw [thick,red,decorate,
decoration={snake,amplitude=.2mm,segment length=1mm,post length=0mm}] (p4) -- (0.4cm,0.2cm) -- (1.1cm,0.9cm) -- (0.4cm,1.6cm);
\draw [thick,red,decorate,
decoration={snake,amplitude=.2mm,segment length=1mm,post length=0mm}] (p3) -- (-0.4cm,1.6cm);
\draw [thick,red,decorate,
decoration={snake,amplitude=.2mm,segment length=1mm,post length=0mm}] (p1s) -- (-1.6cm,-0.5cm);
\draw [thick,red,decorate,
decoration={snake,amplitude=.2mm,segment length=1mm,post length=0mm}] (p2s) -- (-1.6cm,-0.65cm);
\draw [thick,red,decorate,
decoration={snake,amplitude=.2mm,segment length=1mm,post length=0mm}] (p3s) -- (-1.6cm,-0.9cm);
\draw [thick,red,decorate,
decoration={snake,amplitude=.2mm,segment length=1mm,post length=0mm}] (p4s) -- (-1.6cm,-1.3cm);
\draw [thick,red,decorate,
decoration={snake,amplitude=.2mm,segment length=1mm,post length=0mm}] (p5s) -- (-1.6cm,-1.1cm);

\draw [thick,postaction={frac arrow=0.85}] (-1.6cm,0.7cm) -- (-1.2cm,0.3cm) -- (-1.0cm,0.5cm) -- (-1.6cm,1.1cm); 
\draw [thick,postaction={frac arrow=0.95}] (-1.6cm,0.65cm) -- (-0.75cm,-0.2cm) -- (0.55cm,-0.2cm) -- (1.4cm,0.65cm) -- (1.25cm,0.8cm) -- (0.4cm,-0.05cm) -- (0.1cm,-0.05cm) -- (-1.55cm,1.6cm); 
\draw [thick,postaction={frac arrow=0.97}] (-1.6cm,0.6cm) -- (-0.75cm,-0.25cm) -- 
(1.6cm,-0.25cm) -- (1.6cm,0.6cm) -- (1.25cm,0.95cm) -- (0.35cm,0.05cm) -- (-0.9cm,1.3cm) -- (-0.7cm,1.5cm) -- (0.2cm,0.6cm) -- (0.5cm,0.9cm) -- (-0.2cm,1.6cm); 
\draw [thick,postaction={frac arrow=0.95}] (-1.6cm,0.55cm) -- (-0.75cm,-0.3cm) -- (1.7cm,-0.3cm) -- (1.7cm,0.6cm) -- (0.7cm,1.6cm); 
\draw [thick,postaction={frac arrow=0.95}] (-1.6cm,0.5cm) -- (-0.75cm,-0.35cm) -- (1.8cm,-0.35cm) -- (1.8cm,1.1cm) -- (1.3cm,1.6cm); 
\draw [thick,postaction={frac arrow=0.95}] (-1.6cm,0.45cm) -- (-0.75cm,-0.4cm) -- (1.8cm,-0.4cm) -- (2.2cm,-0.8cm); 
\node at (-1.8cm,1.1cm) {$C_1$};
\node at (-1.4cm,1.8cm) {$C_2$};
\node at (-0.2cm,1.8cm) {$C_3$};
\node at (0.7cm,1.8cm) {$C_4$};
\node at (1.3cm,1.8cm) {$C_5$};
\node at (1.8cm,-0.7cm) {$C_0$};
\draw [thick,black,dashed] (-1.6cm,-0.425cm) -- (2.2cm,-0.425cm);
\end{tikzpicture}\hfill%
    \caption{Left:  admissible branch cuts of $h(z)$ in the $z$-plane for a rational initial condition with $N=5$.  Right:  corresponding contours $C_1,C_2,C_3,C_4,C_5$ and $C_0$ in the general case that $c_n\not\in\ii\epsilon \mathbb{N}$ for each $n$.}
\label{fig:BO-BasicContours}
\end{figure}

We start with the real-valued Lax-Oleinik objective function $h(y)$ defined for $y\in\mathbb{R}$ by~\eqref{eq:h-def}. It will be necessary to analytically continue $y\mapsto h(y)$ to a maximal domain that is generally more complicated than implied by using the principal branch of the logarithm in~\eqref{eq:h-def} and taking $y=z$ to be complex.  To this end, we consider branch cuts $\{\Gamma_n,\bar{\Gamma}_n\}_{n=1}^N$ with the following properties.
\begin{definition}[Branch cuts of $h$]
\label{def:B}
The branch cuts $\Gamma_1,\dots,\Gamma_N$ are pairwise disjoint piecewise-smooth curves each emanating from exactly one of the poles $\{p_n\}_{n=1}^N$ and tending to $z=\infty$ in the other direction asymptotic to the ray $\arg(z)=3\pi/4$, all of which lie in a half-plane $\mathrm{Im}(z)>-\delta$ for some $\delta>0$ sufficiently small (in particular, we assume $\delta<\min_n\{\mathrm{Im}(p_n)\}$.
The branch cuts $\bar{\Gamma}_1,\dots,\bar{\Gamma}_N$ are straight horizontal rays each emanating from exactly one of the conjugate poles $\{p_n^*\}_{n=1}^N$ and extending to $z=\infty$ in the left half-plane.
\end{definition}
Given the branch cuts, we index the poles $\{p_n\}_{n=1}^N$ such that in the vicinity of $z=\infty$ in the upper half-plane, the branch cuts $\Gamma_1,\dots,\Gamma_N$ are ordered left-to-right.  See Figure~\ref{fig:BO-BasicContours}, left panel.
We hence obtain a well-defined function $z\mapsto h(z)$ by analytic continuation from large negative real values of $z$ where $h(z)$ is given by~\eqref{eq:h-def} to the domain $z\in\mathbb{C}\setminus (\Gamma_1\cup\cdots\cup \Gamma_N\cup \bar{\Gamma}_1\cup\cdots\cup\bar{\Gamma}_N)$.

We next define some relevant integration contours in the $z$-plane.  In the terminology of~\cite{RationalBO}, the assumption that the rational initial condition is non-special (see Definition~\ref{def:nonsingular-IC}) means that all indices $n=1,\dots,N$ are non-exceptional, so we can provide a simplified definition.

\begin{definition}\label{def:C}
For $n=1,\dots,N$, $C_n$ denotes a contour originating and terminating at $z=\infty$ in the direction $\arg(z)={3\pi/4}$, lying in the domain of analyticity of $h(z)$, and encircling with counterclockwise orientation precisely the branch cuts of $h(z)$ emanating from each of the points $z=p_m$, $1\le m\le n$.  In the special case that $c_n\in \ii \epsilon\mathbb{N}$, $C_n$ terminates instead at $z=p_n$.  $C_0$ denotes a path in the domain of analyticity of $h(z)$ originating at $z=\infty$ in the direction $\arg(z)={3\pi/4}$  to the {left} of all  branch cuts $\Gamma_1,\dots,\Gamma_N$ of $h(z)$ and terminating at $z=\infty$ in the direction $\arg(z)=-{\pi/4}$.  
\end{definition}
See the right-hand panel of Figure~\ref{fig:BO-BasicContours}.

Given $\epsilon>0$ and $t>0$, we next define $(N+1)\times (N+1)$ matrices $\mathbf{A}(t,x;\epsilon)$ and $\mathbf{B}(t,x;\epsilon)$ as follows.  For $1\leq j\leq N+1$ and $ 2\leq k\leq N+1$, set
\begin{align}
A_{j1}(t,x;\epsilon):=\int_{C_{j-1}}u_0(z)\ee^{-\ii h(z)/\epsilon}\,\dd z, 
	&\quad
A_{jk}(t,x;\epsilon):=\int_{C_{j-1}}\frac{\ee^{-\ii h(z)/\epsilon}\,\dd z}{z-p_{k-1}},
\label{eq:Amatrix}
\\
B_{j1}(t,x;\epsilon):=\int_{C_{j-1}}\ee^{-\ii h(z)/\epsilon}\,\dd z,
	&\quad
B_{jk}(t,x;\epsilon):=\int_{C_{j-1}}\frac{\ee^{-\ii h(z)/\epsilon}\,\dd z}{z-p_{k-1}}=A_{jk}(t,x;\epsilon).
\label{eq:Bmatrix}
\end{align}

Then the main result of~\cite{RationalBO} is 
\begin{theorem}[Solution of Benjamin-Ono for rational initial data]\label{thm:inversion-formula}
For each $\epsilon>0$, the exact solution of the Cauchy initial-value problem for the Benjamin-Ono  equation~\eqref{eq:BO} with rational initial condition $u(x,0;\epsilon)=u_0(x)$ of the form \eqref{eq:rationalIC} is 
\begin{equation}
    u(t,x;\epsilon)=2\mathrm{Re}\left( \frac{\det(\mathbf{A}(t,x;\epsilon))}{\det(\mathbf{B}(t,x;\epsilon))}\right),\quad t>0.
\label{eq:lambda-formula}
\end{equation}
\end{theorem}
Note that as shown in~\cite{RationalBO} the denominator cannot vanish for any $(t,x)$ with $t>0$, and that the formula returns the same result regardless of how the branch cuts and contours are chosen as long as they are consistent with Definitions~\ref{def:B} and~\ref{def:C} respectively.

\begin{remark}[Change of contours]
The contours $W_1=C_1$, $W_2=C_2-C_1$, $W_3=C_3-C_2$,\dots, $W_N=C_N-C_{N-1}$ (the notation meaning concatenation of the contours reversing orientation of the second one) each enclose just one of the branch cuts of $h$.  Conversely, given contours $W_n$, $n=1,\dots,N$, with $W_n$ enclosing only the branch cut of $h$ emanating from $z=p_n$ once in the counterclockwise sense (terminating at $z=p_n$ in the exceptional case that $c_n\in\ii\epsilon\mathbb{N}$), we can take $C_n=W_1+W_2+\cdots + W_n$ (concatenation) as the contours in Definition~\ref{def:C}.  We make the identification $W_0=C_0$. Since replacing the contours $\{C_{n}\}_{n=0}^N$ in~\eqref{eq:Amatrix}--\eqref{eq:Bmatrix} with $\{W_{n}\}_{n=0}^N$ amounts to multiplication of the matrices $\mathbf{A}(t,x;\epsilon)$ and $\mathbf{B}(t,x;\epsilon)$ on the left by the same invertible matrix (which happens to have unit determinant), the solution formula~\eqref{eq:lambda-formula} in Theorem~\ref{thm:inversion-formula} is unchanged if we use $\{W_{n}\}_{n=0}^N$ instead of $\{C_{n}\}_{n=0}^N$ to define the integrals, and assuming non-special initial data $u_0$, we will do so going forward.  
    \label{rem:change-contours}
\end{remark}

\begin{remark}[Choice of new contours]
By the generalized Cauchy integral theorem, each contour $W_n$ can be taken to lie exactly along opposite sides of the branch cut $\Gamma_n$ (along one side in the exceptional case that $c_n\in\ii\epsilon\mathbb{N}$) except possibly near the pole $p_n$.  Thus away from the pole $p_n$, the arcs on opposite sides of $\Gamma_n$ may be taken to be exactly the same curves with opposite orientation and with the exponent $h(z)$ taking values differing by a constant.  By a similar deformation argument, several pairwise disjoint branch cuts $\Gamma_n$ ``wrapped by'' corresponding contours $W_n$ may be taken to lie in part against each other in the $z$-plane.
    \label{rem:move-branch-cuts}
\end{remark}

\section{The case of one pole}\label{section:N1}  

In this section, we give a general idea of proof of Theorem~\ref{thm:u-app} in the case $N=1$ and $p_1=\ii$, assuming also for convenience that $u_0$ and $(t,x)$ are in nonspecial configuration according to Definition~\ref{def:nonsingular-IC} (see Section~\ref{sec:typical-behavior} below).
Note that once we have chosen $N=1$, there is no loss of generality in assuming that $p_1=\ii$ (by real translation and positive scaling in the complex $z$-plane, compensating for the scaling by adding a real constant to $h(z)$).

Let $(t,x)\in\mathbb{R}^2$. According to Theorem~\ref{thm:inversion-formula} and Remark~\ref{rem:change-contours}, one can write
\begin{equation}\label{eq:Pi-N1}
u (t,x;\epsilon)
	=2\operatorname{Re}\left(\frac{\det(\mathbf{A}(t,x;\epsilon))}{\det(\mathbf{B}(t,x;\epsilon))}\right),
\end{equation}
where
\begin{equation}\label{eq:AB-N1}
\begin{split}
    \mathbf{A}(t,x;\epsilon)
	&=\begin{bmatrix}\displaystyle 
\int_{W_{0}}\ee^{-\ii h(z)/\varepsilon}u_0(z)\dd z & 	\displaystyle \int_{W_{0}}\ee^{-\ii h(z)/\varepsilon}\frac{\dd z}{z-p_1} \\
\displaystyle \int_{W_1}\ee^{-\ii h(z)/\varepsilon}u_0(z)\dd z & 	\displaystyle \int_{W_1}\ee^{-\ii h(z)/\varepsilon}\frac{\dd z}{z-p_1} \\
\end{bmatrix}, \\
\mathbf{B}(t,x;\epsilon)
	&=\begin{bmatrix}
\displaystyle \int_{W_{0}}\ee^{-\ii h(z)/\varepsilon}\dd z 
	& 	\displaystyle \int_{W_{0}}\ee^{-\ii h(z)/\varepsilon}\frac{\dd z}{z-p_1} \\
	\displaystyle \int_{W_1}\ee^{-\ii h(z)/\varepsilon}\dd z
	 & \displaystyle 	\int_{W_1}\ee^{-\ii h(z)/\varepsilon}\frac{\dd z}{z-p_1} \\
\end{bmatrix}.
\end{split}
\end{equation}
We write $2c_1=a-\ii b$ where $(a,b)\in\mathbb{R}^2\setminus\{(0,0)\}$, so that
\begin{equation}
u_0(x)=\frac{ax+b}{x^2+1}.
\end{equation}
Equation~\eqref{eq:critical-points} has three roots $y_0,y_1,y_2\in\mathbb{C}$:
\begin{equation}
(y-x)(y^2+1)+2t(ay+b)
	=(y-y_0)(y-y_1)(y-y_2).
	\label{eq:critical-points-N1}
\end{equation}
When $J=0$, then $y_0$ is real and the other two are a complex-conjugate pair denoted $y_1=z_1\in\mathbb{C}_+$ and $y_2=z_1^*$, and when $J=1$, the roots $y_0>y_1>y_2$ are all real.

\subsection{Determination of suitable branch cuts and integration contours}\label{sec:N1-ChoiceOfContours}

The first step is to carefully choose the branch cut of $h(z)=h(z;t,x)$, denoted $\Gamma_1$. Then, we choose the contours $W_0$, $W_1$ in $\mathbb{C}\setminus\Gamma_1$ consistent with Definition~\ref{def:C} and Remark~\ref{rem:change-contours} such that the integrals along each contour are dominated by contributions from different critical points of $h(z)$ so that the leading term of the determinants of $\mathbf{A}(t,x;\epsilon)$ and $\mathbf{B}(t,x;\epsilon)$ can be computed as the determinants of the leading terms of the matrix elements themselves. We claim that the following proposition is true (see Proposition~\ref{prop:W} for the general case with proof).

\begin{proposition}[Existence of steepest descent contours when $N=1$]\label{prop:contours N=1}
There exist contours $W_0,W_1$ related to contours of Definition~\ref{def:C} according to Remark~\ref{rem:change-contours} and such that one of the following alternatives hold.

If $J=0$, 
\begin{itemize}
\item on the contour $W_0$, $\operatorname{Re}(-\ii h(z))$ achieves its maximal value only at $y_0\in \mathbb{R}$;
\item on the contour $W_1$, $\operatorname{Re}(-\ii h(z))$ achieves its maximal value only at $z_1\in\mathbb{C}_+$ and only on one side of the branch cut $\Gamma_1$.
\end{itemize}

If $J=1$,
\begin{itemize}
\item on the contour $W_0$, $\operatorname{Re}(-\ii h(z))$ achieves its maximal value only at $y_0,y_1$ and $y_2$;
\item on the contour $W_1$, $\operatorname{Re}(-\ii h(z))$ achieves its maximal value only at $y_1,y_2$ and only on one side of the branch cut $\Gamma_1$.
\end{itemize}
\end{proposition}

\begin{figure}
\centering
\begin{tikzpicture}

\node (p1) at (3,3) [circle,fill=blue,minimum size=4pt,inner sep=0pt]{};

\node (y2) at (-4,0) [circle,fill=black,minimum size=4pt,inner sep=0pt]{};
\node (y1) at (-2,0) [circle,fill=black,minimum size=4pt,inner sep=0pt]{};
\node (y0) at (2,0) [circle,fill=black,minimum size=4pt,inner sep=0pt]{};

\draw[thick,postaction={mid arrow}] (y2) -- (-6,0);
\draw[thick,postaction={mid arrow}] (y2) -- (y1);
\draw[thick,postaction={mid arrow}] (y0) -- (y1);
\draw[thick,postaction={mid arrow}] (y0) -- (5,0);

\draw[thick,postaction={mid arrow}] (-4.1,-1) to[out=80,in=-90] (y2);
\draw[thick,postaction={midp arrow}] (y1) to[out=-90,in=95] (-1.9,-1);
\draw[thick,postaction={mid arrow}] (2.1,-1) to[out=100,in=-90] (y0);

\draw[thick,postaction={mid arrow}] (-4.1,1) to[out=-80,in=90] (y2);
\draw[thick] (-5,3.5) to[out=-70,in=100] (-4.1,1);
\draw[thick,postaction={midp arrow}] (y1) to[out=90,in=-95] (-1.9,1);
\draw[thick] (-1.9,1) to[out=85,in=-120] (-0.5,3.5);
\bonusspiral[thick](3,3)(90:90)(-1:0)[2];
\draw[thick,postaction={midpp arrow}] (3,2) to[out=-180,in=90] (y0);

\node (sd1) at (2.5,-0.5) [circle,fill=red,minimum size=3pt,inner sep=0pt]{};
\node (sd2) at (1.5,0.5) [circle,fill=red,minimum size=3pt,inner sep=0pt]{};

\node (sd3) at (-1.5,0.5) [circle,fill=red,minimum size=3pt,inner sep=0pt]{};
\node (sd4) at (-2.5,-0.5) [circle,fill=red,minimum size=3pt,inner sep=0pt]{};
\draw[thick,red] (sd3) -- (sd4);

\node (sd5) at (-3.5,-0.5) [circle,fill=red,minimum size=3pt,inner sep=0pt]{};
\node (sd6) at (-4.5,0.5) [circle,fill=red,minimum size=3pt,inner sep=0pt]{};
\draw[thick,red] (sd5) -- (sd6);

\bonusspiral[thick,red](3,3)(90:0)(-0.8:-0.4)[1];
\draw[thick,red] (2.6,3) -- (3,3);
\node (intermediate) at (2.6,3) [circle,fill=red,minimum size=3pt,inner sep=0pt]{};
\draw[thick,red] (3,2.2) to[out=180,in=55] (sd2);
\draw[thick,red] (sd2) .. controls (1.2,0.2) and (-1.2,0.2) .. (sd3);
\draw[thick,red] (sd4) to[out=135,in=45] (sd5);
\draw[thick,red] (sd6) to[out=-145,in=0] (-6,0.3);

\draw[thick,blue,postaction={mid arrow}] (-6,0.2) to[out=0,in=-135] (sd6);
\draw[thick,blue,dashed] (sd6) -- (sd5);
\draw[thick,blue,postaction={midp arrow}] (sd5) to[out=25,in=155] (sd4);
\draw[thick,blue,dashed] (sd4) -- (sd3);
\draw[thick,blue,postaction={midp arrow}] (sd3) .. controls (-1.2,0.1) and (1.2,0.1) .. (sd2);
\draw[thick,blue,postaction={midp arrow}] (sd2) to[out=45,in=-180] (3,2.1);
\bonusspiral[thick,blue](3,3)(90:15)(-0.9:-0.5)[1];
\draw[thick,blue] (2.51,2.86) -- (3,2.86);
\draw[thick,blue] (3,2.86) arc[start angle=-90,end angle=90,x radius=0.14,y radius=0.14];
\draw[thick,blue] (2.73,3.14) -- (3,3.14);
\bonusspiral[thick,blue](3,3)(90:-24)(-0.7:-0.3)[1];
\draw[thick,blue,postaction={midpp arrow}] (3,2.3) to[out=-180,in=25] (1.1,0.6);
\draw[thick,blue] (1.1,0.6) to[out=-155,in=-10] (sd3);
\draw[thick,blue,postaction={midp arrow}] (sd4) to[out=115,in=65] (sd5);
\draw[thick,blue,postaction={midp arrow}] (sd6) to[out=-155,in=0] (-6,0.4);

\draw[thick,blue,dashed] (sd2) -- (sd1);
\draw[thick,blue,postaction={mid arrow}] (sd1) to[out=35,in=-180] (5,-0.3);

\draw[dotted,thick] (p1) -- (2.3,2.5);
\node at (2.2,2.4) {$p_1$};
\node at (-4.3,-0.3) {$y_2$};
\node at (-1.6,-0.3) {$y_1$};
\node at (1.7,-0.3) {$y_0$};

\node[blue] at (3.5,-0.5) {$W_0$};
\node[blue] at (1,1) {$W_1$};

\node[red] at (2.5,0.5) {$y_{0,\beta}$};
\draw[red,dotted,thick] (sd2) -- (2.3,0.5);

\node[red] at (2.6,-0.75) {$y_{0,\alpha}$};



\node (p1) at (3,3) [circle,fill=blue,minimum size=4pt,inner sep=0pt]{};

\end{tikzpicture}
\caption{The union of the red curves is the branch cut $\Gamma_1$, emanating from the pole $p_1$ in a $J(t,x)=1$ scenario. The black curves are the critical trajectories $\mathrm{Re}(-\ii h(z))=0$ oriented according to increasing parametrization by $s$, and the blue curves $W_0$, $W_1$ are the integration contours.  The dashed blue line segments passing $y_2$, $y_1$, $y_0$ are to be understood as local steepest descent contours passing through their respective saddle points.}
\label{fig:branch-cut N=1}
\end{figure}
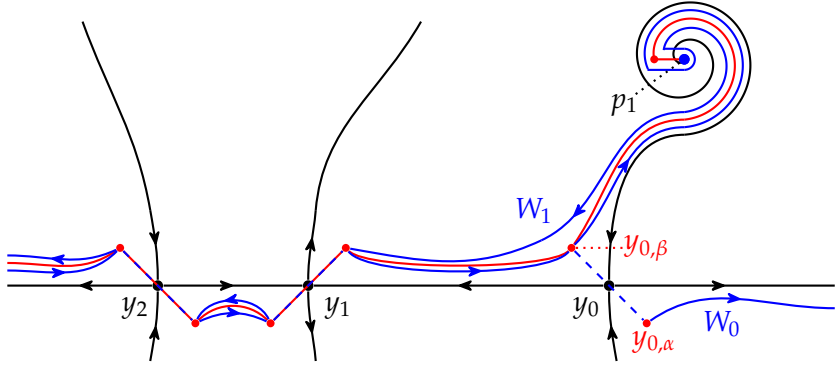

The idea of the proof of this Proposition (and its generalization in Section~\ref{section:StokesGraphs}) is to construct a systematic method (for each $(t,x)$) of selecting various extended level curves to create the branch cut $\Gamma_1$ emanating from the pole $p_1$ with desirable properties.  The steepest descent contour $W_1$ then travels along each side of the branch cut $\Gamma_1$ (with a small exception around $p_1$), see Remark~\ref{rem:move-branch-cuts}.  The remaining contour $W_0$ has a slightly different construction but is also specified in a similar systematic way.  See Figure~\ref{fig:branch-cut N=1} for an example of $\Gamma_1$ and $W_0$, $W_1$ in an $N=J(t,x)=1$ scenario.  We note that dashed blue paths on top of $\Gamma_1$ along the local steepest descent paths over $y_2$, $y_1$ are traversed once on each side of $\Gamma_1$ for $W_1$ and that $W_1$ (and $\Gamma_1$) tend to $\infty$ in the direction $\arg(z)=3\pi/4$.  The contour $W_0$ is the union of
\begin{itemize}
    \item the lower side of $W_1$ until reaching the point $y_{0,\alpha}$,
    \item the dashed blue local steepest descent contour over $y_0$ which connects $y_{0,\alpha}$ to $y_{0,\beta}$,
    \item the blue path from $y_{0,\beta}$ that tends to $\infty$ in the direction $\arg(z)=-\pi/4$. 
\end{itemize}

\subsection{Steepest descent analysis}

The second step, detailed in Section~\ref{sec:steepest-descent} below in the general case, relies on the steepest descent method  applied to the integrals in~\eqref{eq:AB-N1}.
When $N=1$, we choose $W_0$, $W_1$ as described in Proposition~\ref{prop:contours N=1}.

\begin{lemma}[Steepest descent method]\label{lem:steepest-descent}
Suppose $u_0$ and $(t,x)$ are in non-special configuration (see Definition~\ref{def:nonsingular-IC} below), and $\epsilon>0$ is sufficiently small.  Then,
\begin{equation}
    \int_W \ee^{-\ii h(z;t,x)/\epsilon}F(z)\,\dd z =\ee^{\ii\varphi(t,x)} F(\zeta)\sqrt{\frac{2\pi\epsilon}{|h''(\zeta;t,x)|}}\ee^{-\ii h(\zeta;t,x)/\epsilon}\left(1+\mathcal{O}(\epsilon)\right),
    \label{eq:steepest-descent}
\end{equation}
where $F(z)$ is analytic in a neighborhood of $W$, a contour arc along which $\mathrm{Im}(-\ii h(z;t,x))$ is constant and $\mathrm{Re}(-\ii h (z;t,x))$ achieves its maximum value at $z=\zeta\in W$. The point $z=\zeta$ is a simple critical point of $h(z;t,x)$ at which the tangent to the contour $W$ makes an angle of $\varphi(t,x)$ with the direction $\arg(z-\zeta)=0$. Moreover, the error term is uniform for $(t,x)$ in any connected component of the same set where $J(t,x)$ remains constant.
\end{lemma}

The uniformity of the error term follows from the fact that $u_0$ and $(t,x)$ being in non-special configuration guarantees that the critical points of $h(z;t,x)$ will be simple.  Hence Lemma~\ref{lem:steepest-descent} also applies for calculating the matrix elements of $\mathbf{A}(t,x;\epsilon)$ and $\mathbf{B}(t,x;\epsilon)$ if only $(t,x)\in\mathbb{R}^2\setminus\mathcal{D}$ as is required in Theorem~\ref{thm:u-app}.

\paragraph{Case $J=0$}
Applying the above Lemma to~\eqref{eq:Pi-N1} yields
\begin{equation}
 u(t,x;\epsilon)
 	=2\operatorname{Re}\left(\frac{\det(\mathbf{A}^0(t,x))}{\det(\mathbf{B}^0(t,x))}\right)+\mathcal{O}({\varepsilon}),
 	\label{eq:u-simple-N1}
\end{equation}
where
\begin{equation}
\mathbf{A}^0(t,x)
	=\begin{bmatrix}
u_0(y_0) & \displaystyle	\frac{1}{y_0-p_1} \\
u_0(y_2) & \displaystyle \frac{1}{y_2-p_1} \\
\end{bmatrix},
\quad
\mathbf{B}^0(t,x)
	=\begin{bmatrix}
1 & 	\displaystyle \frac{1}{y_0-p_1} \\
1 & 	\displaystyle \frac{1}{y_2-p_1} \\
\end{bmatrix}.
\end{equation}
According to~\eqref{eq:critical-points}, $u_0(y_j)=(x-y_j)/(2t)$, so evaluating the determinants in~\eqref{eq:u-simple-N1} gives
\begin{equation}
 u(t,x;\epsilon)=-\frac{1}{t}\operatorname{Re}(-x+y_0+y_2-p_1)+\mathcal{O}({\varepsilon}).
\end{equation}
In the special case $p_1=\ii$, then $y_0+y_1+y_2=x$ thanks to~\eqref{eq:critical-points-N1} and since $(x-y_0)/(2t)=u^B_0$, we conclude that
\begin{equation}
u(t,x;\epsilon)
	=u^B_0(t,x)+\mathcal{O}({\varepsilon}).
\end{equation}

\paragraph{Case $J=1$}
In this case, the steepest descent angle is $\phi=-\pi/4$ at $y_0$ and $y_2$, and $\phi=\pi/4$ at $y_1$.  Applying Lemma~\ref{lem:steepest-descent}, we remove the factor $\ee^{ -\ii{\pi}/{4}-\ii{h(y_0)}/{\varepsilon}} /\sqrt{|h''(y_0)|}$ from the first row of both determinants in~\eqref{eq:Pi-N1}, and $\ee^{ -\ii{\pi}/{4}-\ii{h(y_2)}/{\varepsilon}} /\sqrt{|h''(y_2)|}$ from the second row.  To simplify the notation, we define
\begin{equation}
\alpha
	:=-\ii\frac{\sqrt{|h''(y_1)|}}{\sqrt{|h''(y_2)|}}\ee^{\ii{(h(y_1)-h(y_2))}/{\varepsilon}}.
\end{equation}
Thus we obtain~\eqref{eq:u-simple-N1} with the understanding that now
\begin{equation}
\begin{split}
\mathbf{A}^0(t,x;\epsilon)
	&=\begin{bmatrix}
    u_0(y_2)+\alpha u_0(y_1) & \displaystyle{\frac{1}{y_2-p_1}+\frac{\alpha}{y_1-p_1}} \\ u_0(y_0) & \displaystyle{\frac{1}{y_0-p_1}}
\end{bmatrix}, \\
\mathbf{B}^0(t,x;\epsilon)
	&=\begin{bmatrix}
	1+\alpha & \displaystyle{\frac{1}{y_2-p_1}+\frac{\alpha}{y_1-p_1}} \\ 1 & \displaystyle{\frac{1}{y_0-p_1}}
\end{bmatrix}.
\end{split}
\end{equation}
Using that $u_0(y_j)=(x-y_j)/(2t)$, we have the simplicification
\begin{equation}
\mathbf{A}^0(t,x;\epsilon)
	=\frac{1}{2t}\begin{bmatrix}
x-y_1+\alpha(x-y_2)  & \displaystyle{\frac{1}{y_1-p_1}+\frac{\alpha}{y_2-p_1}} 
\\
x-y_0 & \displaystyle{\frac{1}{y_0-p_1}} \end{bmatrix}.
\end{equation}
We adopt the notation of~\cite[Lemma 8]{Gerard23} and denote
\begin{equation}
A(y)
	:=\begin{vmatrix}
y  & \displaystyle{\frac{1}{y-p_1}} 
\\
y_0 & \displaystyle{\frac{1}{y_0-p_1}} \end{vmatrix},
\quad
B(y)
	:=\begin{vmatrix}
1 & \displaystyle{\frac{1}{y-p_1}}
\\
1 & \displaystyle{\frac{1}{y_0-p_1}} \end{vmatrix}.
\end{equation}
By linearity, we can expand both determinants along the first row, so that
\begin{equation}
u(t,x;\epsilon)
	=\frac{1}{t}\operatorname{Re} \left(
	\frac{x \left( B(y_1)+ \alpha B(y_2)\right)- \left(A(y_1)+\alpha A(y_2)\right)}
	{B(y_1)+\alpha B(y_2)}\right)
	+\mathcal{O}(\epsilon).
\end{equation}
According to~\cite[Lemma 8]{Gerard23}, or via a direct calculation, we have that
\begin{equation}
A(y)=\left(y+y_0-p_1\right)B(y).
\end{equation}
Consequently,
\begin{equation}
 u(t,x;\varepsilon)
	=\frac{1}{t}\operatorname{Re}\left(x-y_0+p_1
	-\frac{y_1B(y_1)+\alpha y_2B(y_2)}{B(y_1)+\alpha B(y_2)}\right)
	+\mathcal{O}(\varepsilon).
\end{equation}
Taking the real part, since we assumed that $p_1=\ii$, we have
\begin{equation}
u(t,x;\epsilon)
    =u^B_0
    -
    \frac{(u^B_2-u^B_1)\left(|B(y_1)|^2-|\alpha|^2|B(y_2)|^2\right)}{|B(y_1)|^2+|\alpha|^2|B(y_2)|^2+\alpha B(y_2)\overline{B(y_1)}+\overline{\alpha} \overline{B(y_2)}B(y_1)}
    +\mathcal{O}(\varepsilon).
\end{equation}
Taking the derivative of equality~\eqref{eq:hprime} at $y_2$ and $y_1$, we know for $p_1=\ii$ that 
\begin{equation}
|\alpha|^2
	=\frac{|h''(y_1)|}{|h''(y_2)|}
	=\frac{y_0-y_1}{y_0-y_2}\cdot\frac{y_2^2+1}{y_1^2+1},
\end{equation}
and
\begin{equation}
B(y_2)=\frac{y_2-y_0}{(y_2-\ii)(y_0+\ii)}
,\quad
B(y_1)=\frac{y_1-y_0}{(y_1-\ii)(y_0+\ii)}.
\end{equation}
Consequently we simplify $u$ as
\begin{equation}
u(t,x;\epsilon)
    =u^B_0-\frac{(u^B_2-u^B_1)(y_2-y_1)}{ (2y_0-y_2-y_1)
	   -2\sqrt{(y_0-y_2)(y_0-y_1)}\cos(\theta/\varepsilon+\varphi)}
	   +\mathcal{O}(\varepsilon),
\end{equation}
where
\begin{equation}
\theta(t,x)
	=h(y_1)-h(y_2),
\quad
\varphi
	=\arg(y_1-\ii)-\arg(y_2-\ii)+\frac{\pi}{2}.
\end{equation}
Knowing that $(y_j-x)/(2t)=-u_0(y_j)$, and recalling the definition of $r$ in~\eqref{eq:r}, one can also write
\begin{equation}
u(t,x;\epsilon)
    =u^B_0(t,x)+\frac{(u^B_2(t,x)-u^B_1(t,x))(1-r(t,x)^2)}{ 1+r(t,x)^2
	   -2r\cos(\theta(t,x)/\varepsilon+\varphi(t,x))}
	   +\mathcal{O}(\varepsilon),
\end{equation}
which is in agreement with~\eqref{eq:uZD-J1}.  The quantity $\varphi(t,x)$ can then be further expressed in a form involving only the profile of $u_0$ on the real line by using  Lemma~\ref{lem:Phi} below.

\section{Extended level curves}\label{section:LevelSets}
The exponent in the contour integrals~\eqref{eq:Amatrix}, \eqref{eq:Bmatrix} is the function $E(z)$ defined by
\begin{equation}
    E(z):=-\ii h(z).
    \label{eq:exponent}
\end{equation}
The function $z\mapsto\mathrm{Re}(E(z))$ is well-defined and harmonic away from certain branch cuts emanating from $\{p_n,p_n^*\}_{n=1}^N$ consistent with Definition~\ref{def:B}.  In this domain the level set $\{z\in\mathbb{C}: \mathrm{Re}(E(z))=L\}$ for a given $L\in\mathbb{R}$ is a union of analytic arcs called \emph{level curves}.  A level curve is locally characterized by the zero-increment condition 
\begin{equation}
\mathrm{Re}(E'(z)\,\dd z)=\mathrm{Im}(h'(z)\,\dd z)=0.
\label{eq:zero-increment}
\end{equation} 
 Since $z\mapsto h'(z)$ defined by~\eqref{eq:hprime} with~\eqref{eq:rationalIC} is a meromorphic function on the whole $z$-plane, the zero-increment condition makes sense even on the branch cuts of $z\mapsto\mathrm{Re}(E(z))$ and allows the level curve to be analytically continued through any branch cut.  We call an analytic arc satisfying at every point the zero-increment condition, whether or not it intersects a branch cut, an \emph{extended level curve}.  If an extended level curve crosses a branch cut transversely, then $\mathrm{Re}(E(z))$ generally evaluates to two different constants on the arcs of the curve separated by the branch cut.  

The main point is that the extended level curves are well-defined even without specifying any branch cuts of $h(z)$ at all, so we can use them to actually determine branch cuts that will be suitable for subsequent steepest-descent analysis of the contour integrals.  The extended level curves have two different interpretations that will be useful to us:
\begin{itemize}
\item  as trajectories of a quadratic differential;
\item as integral curves of an analytic vector field.
\end{itemize}

\subsection{Extended level curves as trajectories of a quadratic differential}
In the field of geometric function theory, a \emph{quadratic differential} on the Riemann sphere with parameter $z\in\mathbb{C}\cup\{\infty\}$ is an expression of the form $Q(z)\,\dd z^2$ where $Q(z)$ is a rational function.    By a \emph{trajectory} of the quadratic differential we mean a maximal arc $s\mapsto z(s)$ along which $Q(z(s))z'(s)^2>0$, often abbreviated to $Q(z)\,\dd z^2>0$.  There are many results describing the trajectories locally near the poles and zeros of $Q(z)$ as well as globally on the sphere.  As general references, we recommend the books of Jenkins~\cite{Jenkins58} and Strebel~\cite{Strebel84}.  The zero-increment condition~\eqref{eq:zero-increment} fits into this framework if we define the rational function $Q(z)$ as a perfect square:  $Q(z):=h'(z)^2$; thus the extended level curves are precisely the trajectories of the quadratic differential $Q(z)\,\dd z^2=h'(z)^2\,\dd z^2$.

One can see by taking a square root that in general the trajectories of a quadratic differential are locally integral curves of a certain vector field; however since the sign of the square root is arbitrary there is no a priori orientation of these curves, and moreover if $\sqrt{Q(z)}$ has branch points it is impossible to assign a global orientation to the trajectories.  This can lead to some strange behavior not found in the theory of single-valued vector fields such as trajectories with topological closure having nonempty interior (so-called \emph{divergent} or \emph{recurrent} trajectories).  On the other hand, since $Q(z)$ is a rational function, trajectories of quadratic differentials also have some more regular behavior than typical vector fields guarantee.  For instance, a trajectory of a quadratic differential cannot approach a limit cycle in the plane in either direction.  

\subsection{Extended level curves as integral curves of an analytic vector field}
At each point $z=u+\ii v\in\mathbb{C}\setminus\{p_n,p_n^*\}_{n=1}^N$ the zero-increment condition~\eqref{eq:zero-increment} is equivalent to the renormalized form 
\begin{equation}
\mathrm{Im}( f(u,v)h'(u+\ii v)\,(\dd u+\ii \dd v))=0,
\label{eq:zero-renormalized}
\end{equation} where
$f(u,v)\ge 0$ is the function
\begin{equation}
    f(u,v):=\prod_{n=1}^N|u+\ii v-p_n|^2|u+\ii v-p_n^*|^2.
\end{equation}
Set
\begin{equation}
U(u,v):=f(u,v)\mathrm{Re}(h'(u+\ii v)),\quad V(u,v):=-f(u,v)\mathrm{Im}(h'(u+\ii v)),
\label{eq:vectorfieldcomponents}
\end{equation}
and we parametrize an extended level curve by $(u,v)=(u(s),v(s))$ for $s\in\mathbb{R}$.

\begin{lemma} The renormalized zero-increment condition~\eqref{eq:zero-renormalized} is equivalent to the first-order nonlinear autonomous system
\begin{equation}
\frac{\dd u}{\dd s}=U(u,v),\quad \frac{\dd v}{\dd s}=V(u,v).
\label{eq:1st-order-system}
\end{equation}
The extended level curves are also the integral curves (i.e., orbits) of the analytic vector field $(U,V)$ on $\mathbb{R}^2$ given by the first-order system~\eqref{eq:1st-order-system}. 
\end{lemma}

\begin{remark}
 This fact makes available a wide variety of tools from the theory of dynamical systems on the phase plane, for which general references include the book of Guckenheimer and Holmes~\cite{GuckenheimerH90} and the original work of Andronov, Vitt, and Khaiken~\cite{AndronovVK66}.
Note that while $\mathrm{Re}(E(z))$ is guaranteed to remain constant upon analytic continuation along a non-equilibrium orbit of~\eqref{eq:1st-order-system}, the Cauchy-Riemann equations imply that $\mathrm{Re}(E(z))$ will be increasing upon crossing the orbit transversely from right to left by the orientation induced by the parametrization with increasing $s$.
\end{remark}

\begin{proof}
The purpose of the renormalizing factor $f(u,v)$ is to cancel the poles from $h'(z)$ so that the vector field $(u,v)\mapsto (U(u,v),V(u,v))$ is analytic on the plane $\mathbb{R}^2$.  Indeed, near a pole $p_n$ we have $h'(z)=c_n(z-p_n)^{-1}+\mathcal{O}(1)$ where the error term represents a function analytic at $z=p_n$ with real and imaginary parts both real-analytic at $(u,v)=(\mathrm{Re}(p_n),\mathrm{Im}(p_n))$.  Then $f(u,v)/|u+\ii v-p_n|^2$ is real-analytic and nonvanishing at $(u,v)=(\mathrm{Re}(p_n),\mathrm{Im}(p_n))$ and hence for $z$ near $p_n$,
\begin{equation}
\begin{split}
U(u,v)&=\frac{f(u,v)}{|u+\ii v-p_n|^2}\cdot \left[\mathrm{Re}(c_n(u-\ii v-p_n^*))+\mathcal{O}(|u+\ii v-p_n|^2)\right]\\
V(u,v)&=-\frac{f(u,v)}{|u+\ii v-p_n|^2}\cdot\left[\mathrm{Im}(c_n(u-\ii v-p_n^*))+\mathcal{O}(|u+\ii v-p_n|^2)\right],
\end{split}
\label{eq:pole-expansion}
\end{equation}
where the error terms vanishing to second order are real-analytic in $(u,v)$.  A similar analysis is valid near $z=p_n^*$.  This shows that not only are the poles of $h'(z)$ points of analyticity of the vector field $(u,v)\mapsto (U(u,v),V(u,v))$, but also they are equilibrium points.  
\end{proof}

\begin{lemma}
The equilibrium points of~\eqref{eq:1st-order-system} are the poles $\{p_n,p_n^*\}_{n=1}^N$ and the zeros $\{y_k\}_{k=0}^{2N}$ of $h'(z)$.   If $\mathrm{Im}(c_n)\neq 0$ and $\mathrm{Re}(c_n)<0$, the equilibrium at $z=p_n$ and at $z=p_n^*$ is a \emph{spiral sink} while if $\mathrm{Im}(c_n)\neq 0$ and $\mathrm{Re}(c_n)>0$, it is instead a \emph{spiral source}. Each simple critical point $y_k$ is a hyperbolic equilibrium, a simple saddle point.
\label{lem:VF-critical-points}
\end{lemma}

\begin{proof}
According to~\eqref{eq:pole-expansion}, the Jacobian matrix of the vector field $(U,V)$ at the pole $z=p_n$ is proportional via a positive factor to 
\[
\begin{bmatrix}\mathrm{Re}(c_n) & \mathrm{Im}(c_n)\\-\mathrm{Im}(c_n) & \mathrm{Re}(c_n)\end{bmatrix}
\]
the eigenvalues of which are $\mathrm{Re}(c_n)\pm \ii\mathrm{Im}(c_n)$, i.e., the complex-conjugate pair $c_n$ and $c_n^*$.  Under the assumption that $\mathrm{Re}(c_n)\mathrm{Im}(c_n)\neq 0$, these eigenvalues are neither purely real nor purely imaginary.  If $\mathrm{Re}(c_n)<0$, the equilibrium at $z=p_n$ is therefore a \emph{spiral sink} while if $\mathrm{Re}(c_n)>0$, it is instead a \emph{spiral source}.  Exactly the same is true for the pole at $z=p_n^*$.  Similarly to~\eqref{eq:pole-expansion}, the expansion of $U(u,v)$ and $V(u,v)$ near a zero $y_k$ (assumed simple) of $h'(z)$ is
\begin{equation}
    \begin{split}
        U(u,v)&=f(u,v)\cdot\left[\mathrm{Re}(h''(y_k)(u+\ii v-y_k)) + O(|u+\ii v-y_k|^2)\right]\\
        V(u,v)&=-f(u,v)\cdot\left[\mathrm{Im}(h''(y_k)(u+\ii v-y_k))+O(|u+\ii v-y_k|^2)\right]
    \end{split}
    \label{eq:zero-expansion}
\end{equation}
and since $f(u,v)$ is nonzero at $u+\ii v=y_k$, the Jacobian matrix at this point is proportional via a positive factor to the real symmetric matrix
\[
\begin{bmatrix}
    \mathrm{Re}(h''(y_k)) & -\mathrm{Im}(h''(y_k))\\-\mathrm{Im}(h''(y_k)) & -\mathrm{Re}(h''(y_k))
\end{bmatrix}
\]
which has distinct eigenvalues $\pm|h''(y_k)|$ making $y_k$ a \emph{saddle point}, always a hyperbolic equilibrium.  Each such point therefore has exactly two orbits converging toward $y_k$ as $s\to+\infty$ tangent to the eigenspace for the negative eigenvalue (the \emph{stable manifold}) and exactly two orbits converging toward $y_k$ as $s\to -\infty$ tangent to the orthogonal eigenspace for the positive eigenvalue (the \emph{unstable manifold}).
\end{proof}

\subsection{Conditions guaranteeing typical behavior of extended level curves}
\label{sec:typical-behavior}
Assume that $t>0$.  We will be able to specify contours of integration $\{W_n\}_{j=0}^N$ suitable for our analysis if the extended level curves have the following properties.
\begin{definition}
    Let $u_0$ be a  rational initial condition as in \eqref{eq:rationalIC}.  For given $(t,x)$ with $t>0$, the extended level curves have typical behavior if
    \begin{itemize}
        \item the $2N+1$ zeros $\{y_k\}_{k=0}^{2N}$ of $h'(z)$ are all simple, and hence are simple saddle points of the autonomous system~\eqref{eq:1st-order-system};
        \item the autonomous system~\eqref{eq:1st-order-system} has no non-equilibrium periodic orbits;
        \item each of the four extended level curves connected to any zero (two in the limit $s\to+\infty$ and two in the limit $s\to-\infty$) is either an open real interval, or it tends in the opposite direction either to a pole $p_n$ or $p_n^*$ or to $z=\infty$.
    \end{itemize}
    \label{def:typical}
\end{definition}
Note that since $h'(z)=z/(2t)+\mathcal{O}(1)$ as $z\to\infty$, an orbit of~\eqref{eq:1st-order-system} tending to $\infty$ as $s\to+\infty$ does so in the left or right horizontal direction, while one originating from $\infty$ as $s\to-\infty$ does so from the  vertical direction.

We now want to describe simple conditions on $u_0$ and $(t,x)$ such that the extended level curves have typical behavior.  First note that since the zeros of $h'(z)$ are the roots of a polynomial that has degree $2N+1$ whenever $t>0$, the zeros of $h'(z)$ are the points $y$ at which~\eqref{eq:critical-points} holds. The simplicity of the zeros of $h'(z)$ is guaranteed by the nonvanishing of the corresponding polynomial discriminant. That discriminant is itself a polynomial in $(t,x)$ and the real and imaginary parts of $\{p_n\}_{n=1}^N$ and $\{c_n\}_{n=1}^N$, denoted $D((t,x);\{p_n\}_{n=1}^N;\{c_n\}_{n=1}^N)$.

\begin{definition}[Non-special configuration]
We say that the rational initial data $u_0$ and $(t,x)$ are in non-special configuration if the following holds:
\begin{itemize}
\item for every $1\leq n\leq N$, $\operatorname{Re}(c_n)\neq 0$ and $\operatorname{Im}(c_n)\neq 0$;
\item the discriminant is nonzero:
\begin{equation}
    D((t,x);\{p_n\}_{n=1}^N;\{c_n\}_{n=1}^N)\neq 0;
   \label{eq:discriminant}
\end{equation}
\item the inequalities
\begin{equation}
\sum_{n\in S}\mathrm{Re}(c_n)\neq 0,\quad \forall S\in 2^{\{1,\dots,N\}}\setminus\{\emptyset\}
\label{eq:no-homoclinics}
\end{equation}
and
\begin{equation}
    \sum_{n\in S_+}\mathrm{Re}(c_n) -\sum_{n\in S_-}\mathrm{Re}(c_n) \neq \frac{1}{2\pi}\mathrm{Im}\left(\int_{L'}h'(z)\,\dd z\right),\quad \forall S_\pm\in 2^{\{1,\dots,N\}},\; S_+\cap S_-=\emptyset
\label{eq:no-heteroclinics}
\end{equation}
hold, where $2^{\{1,\dots,N\}}$ denotes the set of all subsets of $\{1,\dots,N\}$.
\end{itemize}
These are generic conditions on the parameters  $\{p_n\}_{n=1}^N$, $\{c_n\}_{n=1}^N$ and on $(t,x)\in\mathbb{R}^2$ with $t>0$.
\label{def:nonsingular-IC}
\end{definition}

\begin{proposition}
    Let $u_0$ and $(t,x)$ be in non-special configuration (see Definition~\ref{def:nonsingular-IC}).   Then, the extended level curves have typical behavior in the sense of Definition~\ref{def:typical}.
    \label{prop:guarantee-typical}
\end{proposition}

\begin{proof}
The first condition in Definition~\ref{def:typical} is obviously guaranteed by \eqref{eq:discriminant}.  For the second condition in Definition~\ref{def:typical}, suppose that $C$ is a closed Jordan curve that is an orbit of~\eqref{eq:1st-order-system}.  Then, the zero-increment condition~\eqref{eq:zero-increment} $\mathrm{Im}(h'(z)\,\dd z)=0$ holds at every point of~$C$.  Integrating along the orbit therefore yields the identity
\begin{equation*}
\mathrm{Im}\left(\oint_C h'(z)\,\dd z\right)=0.
\end{equation*}
The Poincar\'e index of a closed orbit is always $1$, and therefore in the interior of $C$ there must lie at least one of the poles (spiral points --- Poincar\'e index $1$), and if there are $P>1$ of them there must also be $P-1$  enclosed critical points (saddle points --- Poincar\'e index $-1$). 
 Therefore, the integral on the left-hand side can also be evaluated by residues at the enclosed poles.  Moreover, since the real line is a finite union of trajectories and saddle points, the poles enclosed by $C$ must all lie either in the upper (residue $c_n$ of $h'(z)$) or lower (residue $c_n^*$ of $h'(z)$) half-plane.  In either case we deduce the identity
\begin{equation*}
\sum_{\text{$n$, $p_n$ is an enclosed pole of $C$}}\mathrm{Re}(c_n) = 0.
\end{equation*}
Therefore, without knowing a priori which of the poles are enclosed by $C$, we may guarantee that there can be no closed orbit if we assume that~\eqref{eq:no-homoclinics} holds.  This is a collection of $2^N-1$ inequalities.  For instance, for these conditions to hold it is sufficient but not necessary that the real parts of $c_n$ all have the same sign, i.e., the spiral points are either all sources or all sinks.  

For the third condition in Definition~\ref{def:typical}, let us assume the first and second conditions already hold.  If we introduce the local coordinate $\zeta=1/z$ to describe the neighborhood of $z=\infty$ on the Riemann sphere, the quadratic differential $h'(z)^2\,\dd z^2$ becomes $\zeta^{-4}h'(\zeta^{-1})^2\,\dd\zeta^2 = (1/(4t^2) + O(\zeta))\zeta^{-6}\,\dd\zeta^2$ so there is a pole of sixth order at $z=\infty$.  According to Jenkins~\cite[Theorem 3.3]{Jenkins58}, it follows that there exists $M>0$ such that any trajectory of $h'(z)^2\,\dd z^2$ containing a point $z$ with $|z|>M$ tends to $z=\infty$ in at least one direction.  This implies that any trajectory tending to a zero of $h'(z)$ in one direction and that is unbounded actually tends to infinity in the other direction (so cannot have any finite limit points).

Therefore, it suffices to characterize the limit sets of \emph{bounded} trajectories tending to a zero of $h'(z)$ in one direction.  For this, we can apply the theory of Andronov, Vitt, and Khaiken~\cite[Section VI.2]{AndronovVK66} (see also~\cite[Section 1.8]{GuckenheimerH90}) to the analytic vector field~\eqref{eq:1st-order-system} to deduce that the limit sets of each bounded trajectory forward and backward in ``time'' $s$ (the so-called $\omega$-limit set and $\alpha$-limit set respectively of the trajectory) consist of either 
\begin{itemize}
    \item a unique equilibrium point (i.e., a pole or zero of $h'(z)$);
    \item a closed orbit (closed Jordan curve); or
    \item a finite union of equilibrium points and homoclinic or heteroclinic orbits connecting them.
\end{itemize}
Under the assumption in force that the second condition from Definition~\ref{def:typical} holds, there are no closed orbits, discarding the second option. Let us show that the third option is also impossible.

We consider the specific trajectories emanating from the (simple) zeros of $h'(z)$, which are saddle points of the system~\eqref{eq:1st-order-system}.  For each saddle point $y$, we would like to classify the $\omega$-limit sets of the trajectories tending to $y$ as $s\to -\infty$ and the $\alpha$-limit sets of those tending to $y$ as $s\to+\infty$.  By the above analysis, each such limit set must be either:
\begin{itemize}
\item[(a)] empty (i.e., the trajectory escapes to $\infty$, otherwise there would be finite limit points);
\item[(b)] a pole $p_n$ or $p_n^*$; \item[(c)] a saddle point;
\item[(d)] a finite union of poles and/or saddles connected by trajectories.
\end{itemize}
According to Lemma~\ref{lem:VF-critical-points}, the conditions $\mathrm{Re}(c_n)\mathrm{Im}(c_n)\neq 0$ for all $n$ and \eqref{eq:discriminant} then guarantee that 
all of the poles and zeros of $h'(z)$ are hyperbolic equilibria (because the poles are all spiral sources or sinks while the zeros are saddles) and in this case referring to~\cite[Section 1.8]{GuckenheimerH90} the points connected by trajectories in case (d) must all be saddle points, i.e., zeros of $h'(z)$.  

To rule out cases (c) and (d) and hence show that the third condition in Definition~\ref{def:typical} is satisfied, it is enough to guarantee that there are no homoclinic orbits or heteroclinic connections in $\mathbb{C}\setminus\mathbb{R}$ between distinct saddle points.  First suppose that there exists a homoclinic orbit, i.e., a trajectory/orbit of~\eqref{eq:1st-order-system} whose $\alpha$ and $\omega$ limit sets both coincide with the same saddle point.  Adjoining the limiting saddle point to the trajectory produces a closed Jordan curve $C$ along which the zero increment condition holds at every point:  $\mathrm{Im}(h'(z)\,\dd z)=0$.  By the same argument as was used to analyze the second condition in Definition~\ref{def:typical}, the conditions~\eqref{eq:no-homoclinics} also rule out homoclinic orbits.

To preclude heteroclinic connections in $\mathbb{C}\setminus\mathbb{R}$ between a pair of distinct saddle points, we first dispense with the case that both saddle points are real, say $y<y'$, in which case there must be an even number (possibly zero) of saddle points on the real line in the open interval $(y,y')$.  Because each interval of the real line between consecutive real saddle points is itself a trajectory, if there is a heteroclinic connection in $\mathbb{C}\setminus\mathbb{R}$ joining $y$ and $y'$, we can again form a closed Jordan curve $C$ of trajectories and intervening saddle points from that arc and the real arc or arcs in between them.  Applying the same argument as above we see that the conditions~\eqref{eq:no-homoclinics} also preclude the existence of a complex heteroclinic connection between real saddles.  

We finally assume that there is a heteroclinic orbit $L$ connecting two saddle points, at most one of which is real.  Again, since the real line is a union of trajectories and saddle points, if both saddle points are complex they must both lie in the same (upper or lower) half-plane, and regardless of whether $L$ has a real endpoint, $L$ crosses no other real point.  Existence of this orbit implies as above that
\begin{equation*}
\mathrm{Im}\left(\int_L h'(z)\,\dd z\right)=0.
\end{equation*}
By the residue theorem, if $L'$ is any other simple curve in the same half-plane punctured at the poles with the same endpoints and orientation as $L$, then 
\begin{equation}
\mathrm{Im}\left(\int_{L'}h'(z)\,\dd z\right)= 2\pi\sum_{k=1}^N n_k\mathrm{Re}(c_k)
\label{eq:Lprime-integral}
\end{equation}
where $n_1,\dots,n_N$ are some integers, all of which are either zero or $\pm 1$.  Therefore, there can be no heteroclinic orbits with a complex endpoint if, for every pair of saddle points one of which is in the upper half-plane, a curve $L'$ in the punctured half-plane $\mathbb{C}_+$ is specified and one verifies that~\eqref{eq:no-heteroclinics} holds.
\end{proof}

\begin{remark}
The only condition appearing in the statement of Theorem~\ref{thm:u-app} is \eqref{eq:discriminant}, which for fixed rational $u_0$ defines the complement of the discriminant locus $\mathcal{D}$.  The proof of Theorem~\ref{thm:u-app} below relies on an algorithmic construction of suitable integration contours $W_0,\dots,W_N$ that is easiest to specify when the extended level curves have typical behavior in the sense of Definition~\ref{def:typical}.  However, we then provide an argument that constructs suitable contours if only $(t,x)\in\mathbb{R}^2\setminus \mathcal{D}$ with $t>0$.
\end{remark}

\section{Determination of suitable branch cuts and integration contours}
\label{section:StokesGraphs}

The goal of this section is to prove the following proposition.  Recall that $\mathcal{D}\subset\mathbb{R}^2$ is the discriminant locus consisting of points $(t,x)$ for which there exists a critical point of $h(z)$, real or complex, of higher multiplicity.

\begin{proposition}[Existence of steepest descent contours]
Let $u_0$ be rational and assume that $t>0$ and that $(t,x)\in\mathbb{R}^2\setminus\mathcal{D}$.
There exist contours $W_n$ according to Remark~\ref{rem:change-contours} (i.e., of the form $W_n=C_n-C_{n-1}$ where $C_n$ satisfies Definition~\ref{def:C} in the general case as $c_n\not\in\ii\epsilon\mathbb{N}$ or $W_n$ terminating at $p_n$ otherwise) and such that the following holds.
There is a $1-1$ correspondence between contours $W_n$, $n=1,\dots,N$, and the set consisting of the union of 
\begin{itemize}
\item the $N-J$ critical points $y$ of $h(\cdot)$ with $\mathrm{Im}(y)>0$;
\item the $J$ pairs of consecutive real critical points $(y_{2j},y_{2j-1})$ with $j=1,\dots,J$.
\end{itemize}
Moreover:
\begin{itemize}
\item
   if $W_n$ corresponds to a complex critical point $y\in\mathbb{C}_+$, then $\mathrm{Re}(E(z))$ achieves its maximum value only at $z=y$ and only on one side of the branch cut $\Gamma_n$;
\item 
    if $W_n$ corresponds to a pair of consecutive real critical points $(y_{2j},y_{2j-1})$, then $\mathrm{Re}(E(z))$ achieves its maximum value only at both $y_{2j-1}$ and $y_{2j}$ and possibly other pairs $(y_{2j'},y_{2j'-1})$ with $j<j'\le J$ and only on one side of the branch cut $\Gamma_n$.
\end{itemize}
The maximum value of $\mathrm{Re}(E(z))$ along $W_0$ is achieved precisely at each of the real critical points $y_0,\dots,y_{2J}$.
In all cases, the maxima are achieved at simple critical points.
\label{prop:W}
\end{proposition}

The reader willing to take Proposition~\ref{prop:W} as given can jump directly to Section~\ref{sec:steepest-descent}.  The rest of this section is devoted to the proof of Proposition~\ref{prop:W}.  Until the very end of the proof, we will actually assume the stronger condition that $u_0$ and $(t,x)$ are in nonspecial configuration as in Definition~\ref{def:nonsingular-IC}.  Then we will remove this condition assuming only that $(t,x)\in\mathbb{R}^2\setminus\mathcal{D}$ holds.

\subsection{Stokes graphs and their properties}
In the theory of quadratic differentials, trajectories that are asymptotic at one end or the other to a zero of $Q(z)=h'(z)^2$ are called \emph{critical}, and the union of these trajectories is denoted $\mathcal{K}$.  The conditions of Proposition~\ref{prop:guarantee-typical} guarantee typical behavior of the extended level curves implying in particular that there are no divergent/recurrent critical trajectories, and this in turn implies that the closure of $\mathcal{K}$ has an empty interior.  In this situation, the Basic Structure Theorem~\cite[Theorem 3.5]{Jenkins58} asserts that  the complement of $\overline{\mathcal{K}}$ in the Riemann sphere $\overline{\mathbb{C}}$ is a finite disjoint union of ``end'', ``strip'', ``ring'', and ``circle'' domains.  Moreover, since these conditions rule out the existence of any closed trajectories, there can be no ring or circle domains, since by definition these are foliated by closed trajectories.

\begin{definition}[End domain]
An \emph{end domain} $\mathscr{E}$ for $h'(z)^2\,\dd z^2$~\cite[Definition 3.6]{Jenkins58} is a maximal union of trajectories on the Riemann sphere $\overline{\mathbb{C}}$, each of which tends to $\infty$ in both directions, that contains no pole $p_j,p_j^*$ or $\infty$ or zero $y_j$ of $h'(z)$, and that is mapped conformally by a univalent branch of $h(z)$ onto an open half-plane $\mathrm{Im}(h(z))>\text{const}$ or $\mathrm{Im}(h(z))<\text{const}$. 
\end{definition}
 An end domain for $h'(z)^2\,\dd z^2$ necessarily has the point $z=\infty$, but none of the finite poles, on its boundary.  It also has one or more zeros of $h'(z)$ on its boundary, however under the conditions of Proposition~\ref{prop:guarantee-typical} the boundary contains either one or two zeros of $h'(z)$, and if it contains two, then they are both real.  According to~\cite[Theorem 3.5(iii)]{Jenkins58} and the fact that $z=\infty$ is a pole of $h'(z)^2\,\dd z^2$ of order six while all finite poles are double, there are exactly four end domains on $\overline{\mathbb{C}}$.  Moreover due to Schwarz symmetry of $h'(z)$ and the fact that the real line is a union of trajectories and saddle points, two of the end domains lie in the upper half-plane and the remaining two are obtained from these by Schwarz reflection.  

\begin{definition}[Strip domain]
A \emph{strip domain} $\mathscr{S}$ for $h'(z)^2\,\dd z^2$~\cite[Definition 3.8]{Jenkins58} is a maximal union of trajectories on the Riemann sphere $\overline{\mathbb{C}}$, each of which tends to different finite poles (one source and one sink) in opposite directions or to $\infty$ in both directions, that contains no pole $p_j,p_j^*$ or $\infty$ or zero $y_j$ of $h'(z)$, and that is mapped conformally by a univalent branch of $h(z)$ onto a horizontal strip $-\infty<a<\mathrm{Im}(h(z))<b<\infty$.
\end{definition}
  The boundary of the strip domain on $\overline{\mathbb{C}}$ contains the common limiting endpoints of all included trajectories but no other finite or infinite poles of $h'(z)^2\,\dd z^2$. There also has to be at least one zero of $h'(z)$ on its boundary, to provide the necessary obstruction to either $a=-\infty$ or $b=+\infty$.   (It can be that a single critical point of $h(z)$ is the obstruction on both sides of the strip, as shown in the left-hand two panels of Figure~\ref{fig:S-faces-strips-enclosing}.)  According to~\cite[Theorem 3.5(iv)]{Jenkins58} and the absence of circle domains according to the conditions of Proposition~\ref{prop:guarantee-typical}, each finite pole of $h'(z)$ has a neighborhood covered by the union of $\overline{\mathcal{K}}$ and finitely-many strip domains.  

\begin{remark}
End and strip domains have an interpretation in terms of ideal (inviscid, irrotational, and incompressible) plane-parallel steady fluid flow, for which there exists a \emph{stream function} that is harmonic and whose level curves are fluid particle paths.  That stream function is locally defined up to a constant on each end and strip domain as a branch $\psi(u,v)$ of $\mathrm{Im}(h(u+\ii v))$.  The saddle points $y$ then have the interpretation of \emph{stagnation points} at which the fluid velocity vanishes.  Likewise, the finite poles $p_n$ with $\mathrm{Re}(c_n)>0$ are point sources of fluid while if instead $\mathrm{Re}(c_n)<0$ they are fluid sinks.  In this language, an end domain $\mathscr{E}$ can be viewed as an \emph{eddy} of recirculating fluid while a strip domain $\mathscr{S}$ is a \emph{sluice} (or \emph{stream}) transporting fluid from a single point source to a single point sink.  
\end{remark}

The set $\overline{\mathcal{K}}$ is Schwarz symmetric and $\mathbb{R}\subset\overline{\mathcal{K}}$, so it is sufficient to restrict attention to the part of $\overline{\mathcal{K}}$ contained in the closed upper half-plane of the Riemann sphere.  For purposes of visualization, we conformally map the latter onto the closed unit disk $\overline{\mathbb{D}}$ via a fractional linear mapping:
\[
w=-\ii\frac{1+\ii z/\rho}{1-\ii z/\rho},\quad \rho>0.
\]
In our plots, we choose $\rho=10$.  We refer to $\overline{\mathcal{K}}$ or its image restricted to the closed unit disk in the $w$-plane and annotated with arrows on each arc indicating direction of increasing parameter $s$ as the \emph{Stokes graph} of $h'(z)$.
Topologically, the Stokes graph is equivalent to an \emph{abstract Stokes graph} $\mathcal{S}$, which has the following definition.

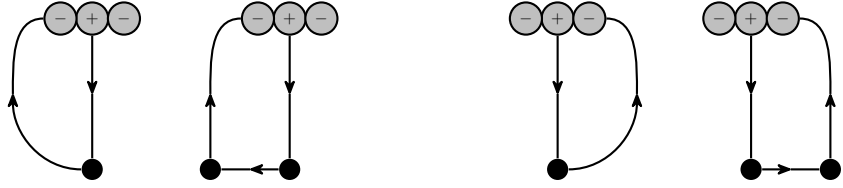
\begin{figure}
\hfill\begin{tikzpicture}
    \node (inftyC) at (0,2cm) [thick,black,circle,fill=lightgray,draw=black,inner sep=2pt] {\tiny $+$};
    \node (inftyL) at (-12pt,2cm) [thick,black,circle,fill=lightgray,draw=black,inner sep=2pt] {\tiny $-$};
    \node (inftyR) at (12pt,2cm) [thick,black,circle,fill=lightgray,draw=black,inner sep=2pt] {\tiny $-$};
    \node (saddle) at (0,0) [circle,fill=black,minimum size=8pt,inner sep=0pt]{};
    \draw [->,>={Stealth[round]},thick] (inftyC) -- (0,1cm);
    \draw [thick] (0,1cm) -- (saddle);
    \draw [->,>={Stealth[round]},thick] (saddle) to[out=180,in=270] (-30pt,1cm);
    \draw [thick] (-30pt,1cm) to[out=90, in=180] (inftyL);
\end{tikzpicture}\hspace{0.05\linewidth}%
\begin{tikzpicture}
    \node (inftyC) at (0,2cm) [thick,black,circle,fill=lightgray,draw=black,inner sep=2pt] {\tiny $+$};
    \node (inftyL) at (-12pt,2cm) [thick,black,circle,fill=lightgray,draw=black,inner sep=2pt] {\tiny $-$};
    \node (inftyR) at (12pt,2cm) [thick,black,circle,fill=lightgray,draw=black,inner sep=2pt] {\tiny $-$};
    \node (saddle1) at (0,0) [circle,fill=black,minimum size=8pt,inner sep=0pt]{};
    \node (saddle2) at (-30pt,0) [circle,fill=black,minimum size=8pt,inner sep=0pt]{};
    \draw [->,>={Stealth[round]},thick] (inftyC) -- (0,1cm);
    \draw [thick] (0,1cm) -- (saddle1);
    \draw [->,>={Stealth[round]},thick] (saddle1) -- (-15pt,0);
    \draw [thick] (-15pt,0) -- (saddle2);
    \draw [->,>={Stealth[round]},thick] (saddle2) to[out=90,in=270] (-30pt,1cm);
    \draw [thick] (-30pt,1cm) to[out=90, in=180] (inftyL);
\end{tikzpicture}\hfill%
\begin{tikzpicture}
    \node (inftyC) at (0,2cm) [thick,black,circle,fill=lightgray,draw=black,inner sep=2pt] {\tiny $+$};
    \node (inftyL) at (-12pt,2cm) [thick,black,circle,fill=lightgray,draw=black,inner sep=2pt] {\tiny $-$};
    \node (inftyR) at (12pt,2cm) [thick,black,circle,fill=lightgray,draw=black,inner sep=2pt] {\tiny $-$};
    \node (saddle) at (0,0) [circle,fill=black,minimum size=8pt,inner sep=0pt]{};
    \draw [->,>={Stealth[round]},thick] (inftyC) -- (0,1cm);
    \draw [thick] (0,1cm) -- (saddle);
    \draw [->,>={Stealth[round]},thick] (saddle) to[out=0,in=270] (30pt,1cm);
    \draw [thick] (30pt,1cm) to[out=90, in=0] (inftyR);
\end{tikzpicture}\hspace{0.05\linewidth}%
\begin{tikzpicture}
    \node (inftyC) at (0,2cm) [thick,black,circle,fill=lightgray,draw=black,inner sep=2pt] {\tiny $+$};
    \node (inftyL) at (-12pt,2cm) [thick,black,circle,fill=lightgray,draw=black,inner sep=2pt] {\tiny $-$};
    \node (inftyR) at (12pt,2cm) [thick,black,circle,fill=lightgray,draw=black,inner sep=2pt] {\tiny $-$};
    \node (saddle1) at (0,0) [circle,fill=black,minimum size=8pt,inner sep=0pt]{};
    \node (saddle2) at (30pt,0) [circle,fill=black,minimum size=8pt,inner sep=0pt]{};
    \draw [->,>={Stealth[round]},thick] (inftyC) -- (0,1cm);
    \draw [thick] (0,1cm) -- (saddle1);
    \draw [->,>={Stealth[round]},thick] (saddle1) -- (15pt,0);
    \draw [thick] (15pt,0) -- (saddle2);
    \draw [->,>={Stealth[round]},thick] (saddle2) to[out=90,in=270] (30pt,1cm);
    \draw [thick] (30pt,1cm) to[out=90, in=0] (inftyR);
\end{tikzpicture}\hfill
\caption{Possible faces of an abstract Stokes graph $\mathcal{S}$ representing end domains.  $\mathcal{S}$ contains precisely one of the left two faces and precisely one of the right two faces.}
    \label{fig:S-faces-ends}
\end{figure}

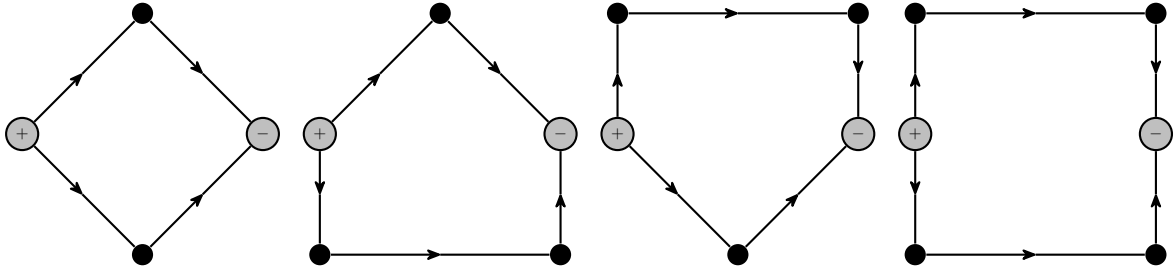
\begin{figure}
\begin{tikzpicture}
    \node (source) at (-1.6cm,0) [thick,black,circle,fill=lightgray,draw=black,inner sep=2pt]{\tiny $+$};
    \node (sink) at (1.6cm,0) [thick,black,circle,fill=lightgray,draw=black,inner sep=2pt]{\tiny $-$};
    \node (saddle1) at (0,1.6cm) [circle,fill=black,minimum size=8pt,inner sep=0pt]{};
    \node (saddle2) at (0,-1.6cm) [circle,fill=black,minimum size=8pt,inner sep=0pt]{};
    \draw [->,>={Stealth[round]},thick] (source) -- (-0.8cm,0.8cm);
    \draw [thick] (-0.8cm,0.8cm) -- (saddle1);
    \draw [->,>={Stealth[round]},thick] (source) -- (-0.8cm,-0.8cm);
    \draw [thick] (-0.8cm,-0.8cm) -- (saddle2);
    \draw[->,>={Stealth[round]},thick] (saddle1) -- (0.8cm,0.8cm);
    \draw [thick] (0.8cm,0.8cm) -- (sink);
    \draw[->,>={Stealth[round]},thick] (saddle2) -- (0.8cm,-0.8cm);
    \draw [thick] (0.8cm,-0.8cm) -- (sink);
\end{tikzpicture}\hfill%
\begin{tikzpicture}
    \node (source) at (-1.6cm,0) [thick,black,circle,fill=lightgray,draw=black,inner sep=2pt]{\tiny $+$};
    \node (sink) at (1.6cm,0) [thick,black,circle,fill=lightgray,draw=black,inner sep=2pt]{\tiny $-$};
    \node (saddle1) at (0,1.6cm) [circle,fill=black,minimum size=8pt,inner sep=0pt]{};
   \node (saddle2) at (-1.6cm,-1.6cm) [circle,fill=black,minimum size=8pt,inner sep=0pt]{};
   \node (saddle3) at (1.6cm,-1.6cm) [circle,fill=black,minimum size=8pt,inner sep=0pt]{};
   \draw [->,>={Stealth[round]},thick] (source) -- (-0.8cm,0.8cm);
   \draw [thick] (-0.8cm,0.8cm) -- (saddle1);
   \draw[->,>={Stealth[round]},thick] (saddle1) -- (0.8cm,0.8cm);
   \draw [thick] (0.8cm,0.8cm) -- (sink);
   \draw[->,>={Stealth[round]},thick] (source) -- (-1.6cm,-0.8cm);
   \draw [thick] (-1.6cm,-0.8cm) -- (saddle2);
   \draw [->,>={Stealth[round]},thick] (saddle2) -- (0,-1.6cm);
   \draw [thick] (0,-1.6cm) -- (saddle3);
   \draw[->,>={Stealth[round]},thick] (saddle3) -- (1.6cm,-0.8cm);
   \draw [thick] (1.6cm,-0.8cm) -- (sink);
\end{tikzpicture}\hfill%
\begin{tikzpicture}
    \node (source) at (-1.6cm,0) [thick,black,circle,fill=lightgray,draw=black,inner sep=2pt]{\tiny $+$};
    \node (sink) at (1.6cm,0) [thick,black,circle,fill=lightgray,draw=black,inner sep=2pt]{\tiny $-$};
    \node (saddle1) at (0,-1.6cm) [circle,fill=black,minimum size=8pt,inner sep=0pt]{};
   \node (saddle2) at (-1.6cm,1.6cm) [circle,fill=black,minimum size=8pt,inner sep=0pt]{};
   \node (saddle3) at (1.6cm,1.6cm) [circle,fill=black,minimum size=8pt,inner sep=0pt]{};
   \draw [->,>={Stealth[round]},thick] (source) -- (-0.8cm,-0.8cm);
   \draw [thick] (-0.8cm,-0.8cm) -- (saddle1);
   \draw[->,>={Stealth[round]},thick] (saddle1) -- (0.8cm,-0.8cm);
   \draw [thick] (0.8cm,-0.8cm) -- (sink);
   \draw[->,>={Stealth[round]},thick] (source) -- (-1.6cm,0.8cm);
   \draw [thick] (-1.6cm,0.8cm) -- (saddle2);
   \draw [->,>={Stealth[round]},thick] (saddle2) -- (0,1.6cm);
   \draw [thick] (0,1.6cm) -- (saddle3);
   \draw[->,>={Stealth[round]},thick] (saddle3) -- (1.6cm,0.8cm);
   \draw [thick] (1.6cm,0.8cm) -- (sink);    
\end{tikzpicture}\hfill%
\begin{tikzpicture}
    \node (source) at (-1.6cm,0) [thick,black,circle,fill=lightgray,draw=black,inner sep=2pt]{\tiny $+$};
    \node (sink) at (1.6cm,0) [thick,black,circle,fill=lightgray,draw=black,inner sep=2pt]{\tiny $-$};
    \node (saddleNW) at (-1.6cm,1.6cm) [circle,fill=black,minimum size=8pt,inner sep=0pt]{};
    \node (saddleNE) at (1.6cm,1.6cm) [circle,fill=black,minimum size = 8pt,inner sep=0pt]{};
    \node (saddleSW) at (-1.6cm,-1.6cm) [circle,fill=black,minimum size=8pt,inner sep=0pt]{};
    \node (saddleSE) at (1.6cm,-1.6cm) [circle,fill=black,minimum size=8pt,inner sep=0pt]{};
    \draw [->,>={Stealth[round]},thick] (source) -- (-1.6cm,0.8cm);
    \draw [thick] (-1.6cm,0.8cm) -- (saddleNW);
    \draw [->,>={Stealth[round]},thick] (saddleNW) -- (0cm,1.6cm);
    \draw [thick] (0cm,1.6cm) -- (saddleNE);
    \draw [->,>={Stealth[round]},thick] (saddleNE) -- (1.6cm,0.8cm);
    \draw [thick] (1.6cm,0.8cm) -- (sink);
    \draw [->,>={Stealth[round]},thick] (source) -- (-1.6cm,-0.8cm);
    \draw [thick] (-1.6cm,-0.8cm) -- (saddleSW);
    \draw [->,>={Stealth[round]},thick] (saddleSW) -- (0cm,-1.6cm);
    \draw [thick] (0cm,-1.6cm) -- (saddleSE);
    \draw [->,>={Stealth[round]},thick] (saddleSE) -- (1.6cm,-0.8cm);
    \draw [thick] (1.6cm,-0.8cm) -- (sink);
\end{tikzpicture}
    \caption{Possible faces of an abstract Stokes graph $\mathcal{S}$ representing strip domains with each source and sink having two boundary edges. The symbols \saddlenode\ represent a saddle point, and \sourcenode\ or \sinknode\ represent poles or a part of the compound vertex \inftynode.}
    \label{fig:S-faces-strips-simple}
\end{figure}

\begin{figure}
\begin{tikzpicture}
    \node (saddleS) at (0cm,-1.6cm) [circle,fill=black,minimum size=8pt,inner sep=0pt]{};
    \node (sink) at (0cm,0) [thick,black,circle,fill=lightgray,draw=black,inner sep=2pt]{\tiny $-$};
    \node (source) at (0,1.6cm) [thick,black,circle,fill=lightgray,draw=black,inner sep=2pt]{\tiny $+$};
    \draw [->,>={Stealth[round]},thick] (saddleS) -- (0cm,-0.8cm);
    \draw [thick] (0cm,-0.8cm) -- (sink);
    \draw [->,>={Stealth[round]},thick] (source) to[out=180+45,in=90] (-0.8cm,0cm);
    \draw [thick] (-0.8cm,0cm) to[out=-90,in=180] (saddleS);
    \draw [->,>={Stealth[round]},thick] (source) to[out=0-45,in=90] (0.8cm,0cm);
    \draw [thick] (0.8cm,0cm) to[out=-90,in=0] (saddleS);    
\end{tikzpicture}\hfill%
\begin{tikzpicture}
    \node (saddleS) at (0cm,-1.6cm) [circle,fill=black,minimum size=8pt,inner sep=0pt]{};
    \node (source) at (0cm,0) [thick,black,circle,fill=lightgray,draw=black,inner sep=2pt]{\tiny $+$};
    \node (sink) at (0,1.6cm) [thick,black,circle,fill=lightgray,draw=black,inner sep=2pt]{\tiny $-$};
    \draw [->,>={Stealth[round]},thick] (source) -- (0cm,-0.8cm);
    \draw [thick] (0cm,-0.8cm) -- (saddleS);
    \draw [->,>={Stealth[round]},thick] (saddleS) to[out=180,in=270] (-0.8cm,0cm);
    \draw [thick] (-0.8cm,0cm) to[out=90,in=180+45] (sink);
    \draw [->,>={Stealth[round]},thick] (saddleS) to[out=0,in=270] (0.8cm,0cm);
    \draw [thick] (0.8cm,0cm) to[out=90,in=0-45] (sink);
\end{tikzpicture}\hfill%
\begin{tikzpicture}
    \node (saddleS) at (0cm,-1.6cm) [circle,fill=black,minimum size=8pt,inner sep=0pt]{};
    \node (saddleSW) at (-0.8cm,-1.6cm) [circle,fill=black,minimum size=8pt,inner sep=0pt]{};
    \node (saddleSE) at (0.8cm,-1.6cm) [circle,fill=black,minimum size=8pt,inner sep=0pt]{};
    \node (sink) at (0cm,0) [thick,black,circle,fill=lightgray,draw=black,inner sep=2pt]{\tiny $-$};
    \node (source) at (0,1.6cm) [thick,black,circle,fill=lightgray,draw=black,inner sep=2pt]{\tiny $+$};
    \draw [->,>={Stealth[round]},thick] (saddleS) -- (0cm,-0.8cm);
    \draw [thick] (0cm,-0.8cm) -- (sink);
    \draw [->,>={Stealth[round]},thick] (saddleSW) -- (-0.3cm,-1.6cm);
    \draw [thick] (-0.3cm,-1.6cm) -- (saddleS);
    \draw [->,>={Stealth[round]},thick] (saddleSE) -- (0.3cm,-1.6cm);
    \draw [thick] (0.3cm,-1.6cm) -- (saddleS);
    \draw [->,>={Stealth[round]},thick] (source) to[out=180+45,in=90] (-0.8cm,0cm);
    \draw [thick] (-0.8cm,0cm) to[out=-90,in=90] (saddleSW);
    \draw [->,>={Stealth[round]},thick] (source) to[out=0-45,in=90] (0.8cm,0cm);
    \draw [thick] (0.8cm,0cm) to[out=-90,in=90] (saddleSE);
    \draw [thick] (-0.8cm,0cm) to[out=-90,in=90] (saddleSW);   
\end{tikzpicture}\hfill%
\begin{tikzpicture}
    \node (saddleS) at (0cm,-1.6cm) [circle,fill=black,minimum size=8pt,inner sep=0pt]{};
    \node (saddleSW) at (-0.8cm,-1.6cm) [circle,fill=black,minimum size=8pt,inner sep=0pt]{};
    \node (saddleSE) at (0.8cm,-1.6cm) [circle,fill=black,minimum size=8pt,inner sep=0pt]{};
    \node (source) at (0cm,0) [thick,black,circle,fill=lightgray,draw=black,inner sep=2pt]{\tiny $+$};
    \node (sink) at (0,1.6cm) [thick,black,circle,fill=lightgray,draw=black,inner sep=2pt]{\tiny $-$};
    \draw [->,>={Stealth[round]},thick] (source) -- (0cm,-0.8cm);
    \draw [thick] (0cm,-0.8cm) -- (saddleS);
    \draw [->,>={Stealth[round]},thick] (saddleS) -- (-0.5cm,-1.6cm);
    \draw [thick] (-0.5cm,-1.6cm) -- (saddleSW);
    \draw [->,>={Stealth[round]},thick] (saddleSW) -- (-0.8cm,0cm);
    \draw [thick] (-0.8cm,0cm) to[out=90,in=180+45] (sink);
    \draw [->,>={Stealth[round]},thick] (saddleS) --(0.5cm,-1.6cm);
    \draw [thick] (0.5cm,-1.6cm) -- (saddleSE);
    \
   \draw [->,>={Stealth[round]},thick] (saddleSE) -- (0.8cm,0cm);
    \draw [thick] (0.8cm,0cm) to[out=90,in=0-45] (sink);    
\end{tikzpicture}\hfill%
\begin{tikzpicture}
\node (inftyC) at (0,1.6cm) [thick,black,circle,fill=lightgray,draw=black,inner sep=2pt] {\tiny $+$};
\node (inftyL) at (-12pt,1.6cm) [thick,black,circle,fill=lightgray,draw=black,inner sep=2pt] {\tiny $-$};
\node (inftyR) at (12pt,1.6cm) [thick,black,circle,fill=lightgray,draw=black,inner sep=2pt] {\tiny $-$};
    \node (saddleS) at (0cm,-1.6cm) [circle,fill=black,minimum size=8pt,inner sep=0pt]{};
    \node (saddleSE) at (0.8cm,-1.6cm) [circle,fill=black,minimum size=8pt,inner sep=0pt]{};
    \node (source) at (0cm,0) [thick,black,circle,fill=lightgray,draw=black,inner sep=2pt]{\tiny $+$};
    \draw [->,>={Stealth[round]},thick] (source) -- (0cm,-0.8cm);
    \draw [thick] (0cm,-0.8cm) -- (saddleS);
    \draw [->,>={Stealth[round]},thick] (saddleS) to[out=180,in=-90] (-1cm,0cm);
    \draw[thick] (-1cm,0cm) to[out=90,in=180] (inftyL);
    \draw [->,>={Stealth[round]},thick] (saddleS) --(0.5cm,-1.6cm);
    \draw [thick] (0.5cm,-1.6cm) -- (saddleSE);
    \draw [->,>={Stealth[round]},thick] (saddleSE) to[out=90,in=-45] (0cm,0.5cm);
    \draw[thick] (0cm,0.5cm) to[out=135,in=180] (inftyL);
\end{tikzpicture}\hfill%
\begin{tikzpicture}
\node (inftyC) at (0,1.6cm) [thick,black,circle,fill=lightgray,draw=black,inner sep=2pt] {\tiny $+$};
\node (inftyL) at (-12pt,1.6cm) [thick,black,circle,fill=lightgray,draw=black,inner sep=2pt] {\tiny $-$};
\node (inftyR) at (12pt,1.6cm) [thick,black,circle,fill=lightgray,draw=black,inner sep=2pt] {\tiny $-$};
    \node (saddleS) at (0cm,-1.6cm) [circle,fill=black,minimum size=8pt,inner sep=0pt]{};
    \node (saddleSW) at (-0.8cm,-1.6cm) [circle,fill=black,minimum size=8pt,inner sep=0pt]{};
    \node (source) at (0cm,0) [thick,black,circle,fill=lightgray,draw=black,inner sep=2pt]{\tiny $+$};
    \draw [->,>={Stealth[round]},thick] (source) -- (0cm,-0.8cm);
    \draw [thick] (0cm,-0.8cm) -- (saddleS);
    \draw [->,>={Stealth[round]},thick] (saddleS) to[out=0,in=-90] (1cm,0cm);
    \draw[thick] (1cm,0cm) to[out=90,in=0] (inftyR);
    \draw [->,>={Stealth[round]},thick] (saddleS) --(-0.5cm,-1.6cm);
    \draw [thick] (-0.5cm,-1.6cm) -- (saddleSW);
    \draw [->,>={Stealth[round]},thick] (saddleSW) to[out=90,in=-135] (0cm,0.5cm);
    \draw[thick] (0cm,0.5cm) to[out=45,in=0] (inftyR);
\end{tikzpicture}
    \caption{Possible faces of an abstract Stokes graph $\mathcal{S}$ representing strip domains having a source or sink connected by only one edge. The doubly-connected vertex
    \sourcenode\
    or 
    \sinknode\ 
    on the exterior boundary of each of the faces shown in the left four diagrams can represent a part of the compound vertex \inftynode.}
    \label{fig:S-faces-strips-enclosing}
\end{figure}
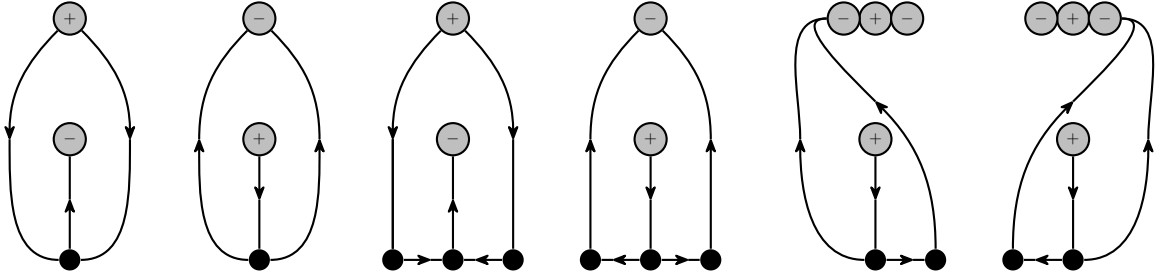

\begin{definition}[Abstract Stokes graph]
An abstract Stokes graph $\mathcal{S}$ with indices $(N,J)$, $N\in\mathbb{Z}_{>0}$ and $J\in\{0,\dots,N\}$ is a connected planar directed graph with the following properties:

\noindent $\bullet$ {\bf Vertices:}  $\mathcal{S}$ has the following vertices:
    \begin{itemize}[-]
        \item $N$ vertices representing the images in $\mathbb{D}$ of the finite poles $p_1,\dots,p_N$.  Each of these vertices either has in-degree or out-degree zero, depending on whether $p_j$ is a source (denoted \tikz[baseline=-3pt]{\node (A) at (0,0) [thick,black, circle,fill=lightgray,draw=black,inner sep=2pt] {\tiny $+$}}) or a sink (denoted \tikz[baseline=-3pt]{\node (A) at (0,0) [thick,black, circle,fill=lightgray,draw=black,inner sep=2pt] {\tiny $-$}}), respectively;
        \item One ``compound'' vertex denoted \tikz[baseline=-3pt]{\node at (0,0) [thick,black,circle,fill=lightgray,draw=black,inner sep=2pt] {\tiny $\infty$}} or \tikz[baseline=-3pt]{\node at (0,0) [thick,black,circle,fill=lightgray,draw=black,inner sep=2pt] {\tiny $+$};\node at (-12pt,0) [thick,black,circle,fill=lightgray,draw=black,inner sep=2pt] {\tiny $-$};\node at (12pt,0) [thick,black,circle,fill=lightgray,draw=black,inner sep=2pt] {\tiny $-$}} representing the point at infinity in the upper half $z$-plane and consisting of two sink-type vertices with an intermediate source-type vertex, each of which has degree at least 1;
        \item $N-J$ vertices denoted \tikz[baseline=-3pt]{\node at (0,0) [circle,fill=black,minimum size=8pt,inner sep=0pt]{}} representing the images in $\mathbb{D}$ of saddle points $y$ belonging to the open upper half $z$-plane.  Each such vertex has in-degree 2 and out-degree 2, with orientation alternating about the vertex.
        \item $2J+1$ vertices also denoted \tikz[baseline=-3pt]{\node at (0,0) [circle,fill=black,minimum size=8pt,inner sep=0pt]{}} representing the images on $\partial\mathbb{D}$ of real saddle points $y\in\mathbb{R}$.  Each such vertex has degree $3$ since there is no need to keep track of the trajectories in the lower half-plane. $J$ of these vertices have out-degree $1$ and $J+1$ of them have in-degree $1$.
    \end{itemize}

\noindent$\bullet$ {\bf Boundary cycle:}  The boundary of the unbounded face of $\mathcal{S}$ is a cycle representing the unit circle in the $w$-plane, i.e., the images of the real trajectories and their limit points on the Riemann sphere.  The boundary cycle therefore consists of the $2J+1$ degree-3 saddle-type vertices and the compound vertex only with alternating edges: \tikz[baseline=-3pt]{\node at (0,0) [thick,black,circle,fill=lightgray,draw=black,inner sep=2pt] {\tiny $\infty$}} $\longleftarrow$ \tikz[baseline=-3pt]{\node at (0,0) [circle,fill=black,minimum size=8pt,inner sep=0pt]{}} $\longrightarrow$ \tikz[baseline=-3pt]{\node at (0,0) [circle,fill=black,minimum size=8pt,inner sep=0pt]{}} $\longleftarrow \cdots \longleftarrow $ \tikz[baseline=-3pt]{\node at (0,0) [circle,fill=black,minimum size=8pt,inner sep=0pt]{}} $ \longrightarrow$ \tikz[baseline=-3pt]{\node at (0,0) [thick,black,circle,fill=lightgray,draw=black,inner sep=2pt] {\tiny $\infty$}}.  The saddle(s) adjacent to \inftynode\ on the boundary cycle have in-degree $1$ and consecutive saddles alternate around the cycle between in-degree $1$ and out-degree $1$.  All outward directed edges from \tikz[baseline=-3pt]{\node at (0,0) [thick,black,circle,fill=lightgray,draw=black,inner sep=2pt] {\tiny $\infty$}} 
    (there is at least one) are consecutive about the vertex and 
    separate two collections of inward directed edges as suggested by the alternate notation \tikz[baseline=-3pt]{\node at (0,0) [thick,black,circle,fill=lightgray,draw=black,inner sep=2pt] {\tiny $+$};\node at (-12pt,0) [thick,black,circle,fill=lightgray,draw=black,inner sep=2pt] {\tiny $-$};\node at (12pt,0) [thick,black,circle,fill=lightgray,draw=black,inner sep=2pt] {\tiny $-$}}.  Also, with the exception of these $2J$ alternating edges of the boundary cycle, there are no edges connecting two (possibly identical) saddle-type vertices.

\noindent$\bullet$ {\bf Faces:}  the faces of the bounded component of $\mathcal{S}$ are among those shown in Figure~\ref{fig:S-faces-ends} (exactly one of the left-hand pair and exactly one of the right-hand pair) and in Figures~\ref{fig:S-faces-strips-simple} and \ref{fig:S-faces-strips-enclosing}.  In particular, every edge is connected to at least one saddle-type vertex \saddlenode.  Also, treating \inftynode\ as three distinct vertices \tikz[baseline=-3pt]{\node at (0,0) [thick,black,circle,fill=lightgray,draw=black,inner sep=2pt] {\tiny $+$};\node at (-12pt,0) [thick,black,circle,fill=lightgray,draw=black,inner sep=2pt] {\tiny $-$};\node at (12pt,0) [thick,black,circle,fill=lightgray,draw=black,inner sep=2pt] {\tiny $-$}} at most one edge connects any two vertices, and there are no edges connecting a vertex to itself.

\label{def:abstractStokesgraph}
\end{definition}

\subsection{From a Stokes graph to a Stokes tree}
Given the Stokes graph of $h'(z)$, we wish to identify the images in the closed unit disk of certain trajectories that we can subsequently use to define suitable branch cuts of $h(z)$ and paths of integration $\{W_n\}_{n=0}^N$ (see Theorem~\ref{thm:inversion-formula} and Remark~\ref{rem:change-contours}).  For this purpose, we start by identifying the Stokes graph of $h'(z)$ with its topologically-equivalent abstract Stokes graph $\mathcal{S}$ with indices $(N,J)$. We will show that we can extract from this graph a Stokes tree.

\begin{definition}[Stokes tree]
We say that a subgraph $\mathcal{T}$ of $\mathcal{S}$ obtained by selecting certain edges and vertices  is a {\emph Stokes tree} if it is a polytree satisfying the following conditions.  Letting $y_0$, $y_1$, \dots, $y_{2J}$ denote the saddle-type vertices on the boundary cycle of $\mathcal{S}$ enumerated in clockwise order starting at \inftynode, we impose the following rules on $\mathcal{T}$:
\begin{enumerate}
    \item $\mathcal{T}$ contains all vertices of $\mathcal{S}$ with the possible exception of $y_0$, and \inftynode\ is considered as a single vertex of $\mathcal{T}$, designated as the root.
    \item $\mathcal{T}$ contains no edges of $\mathcal{S}$ that are joined to the right-most component \sinknode\ of the compound vertex \inftynode\ (i.e., the component first from the central component \sourcenode\ in the clockwise direction). 
    \item If $J\ge 1$, then $\mathcal{T}$ contains the edge of the boundary cycle joining $y_{2j}$ and $y_{2j-1}$ for $j=1,\dots,J$.
    \item Each saddle-type vertex \saddlenode\ not on the boundary cycle has degree two in $\mathcal{T}$.  If one edge is incoming and the other is outgoing, these edges have the same face of $\mathcal{S}$ on their left.
\end{enumerate}
\end{definition}

The goal of this subsection is to prove the following result.
\begin{proposition}[Extraction of a Stokes tree]
It is possible to extract a Stokes tree from any abstract Stokes graph $\mathcal{S}$ with indices $(N,J)$.
\end{proposition}

\subsubsection{Definition of the Stokes tree when \texorpdfstring{$N=1$}{N=1}}\label{subsub:N1}

\begin{figure}
\begin{center}
    \includegraphics[width=0.3\linewidth]{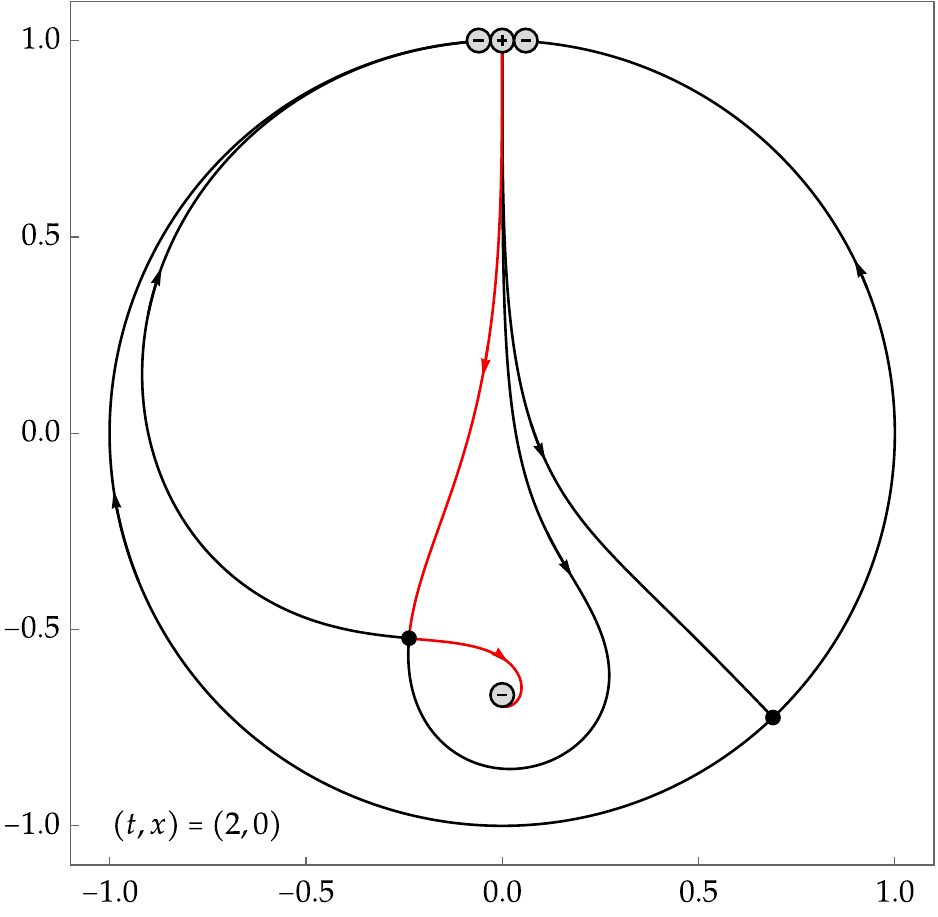}\hspace{0.03\linewidth}%
    \includegraphics[width=0.3\linewidth]{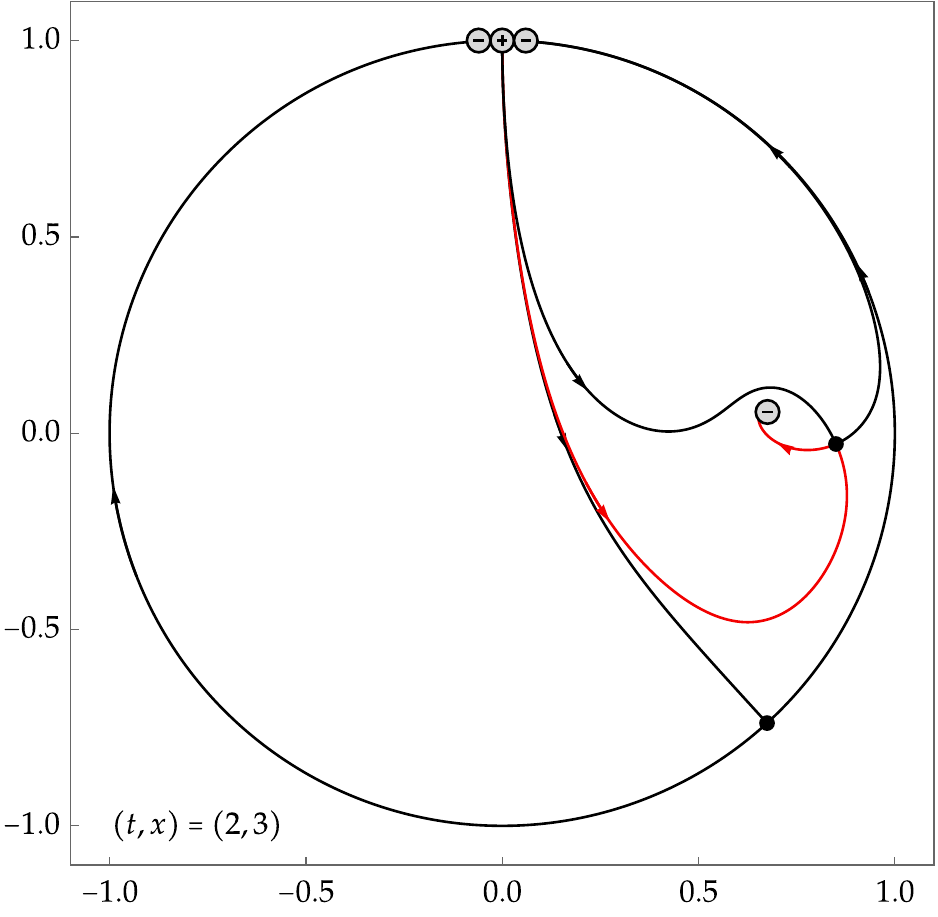}\\
    \includegraphics[width=0.3\linewidth]{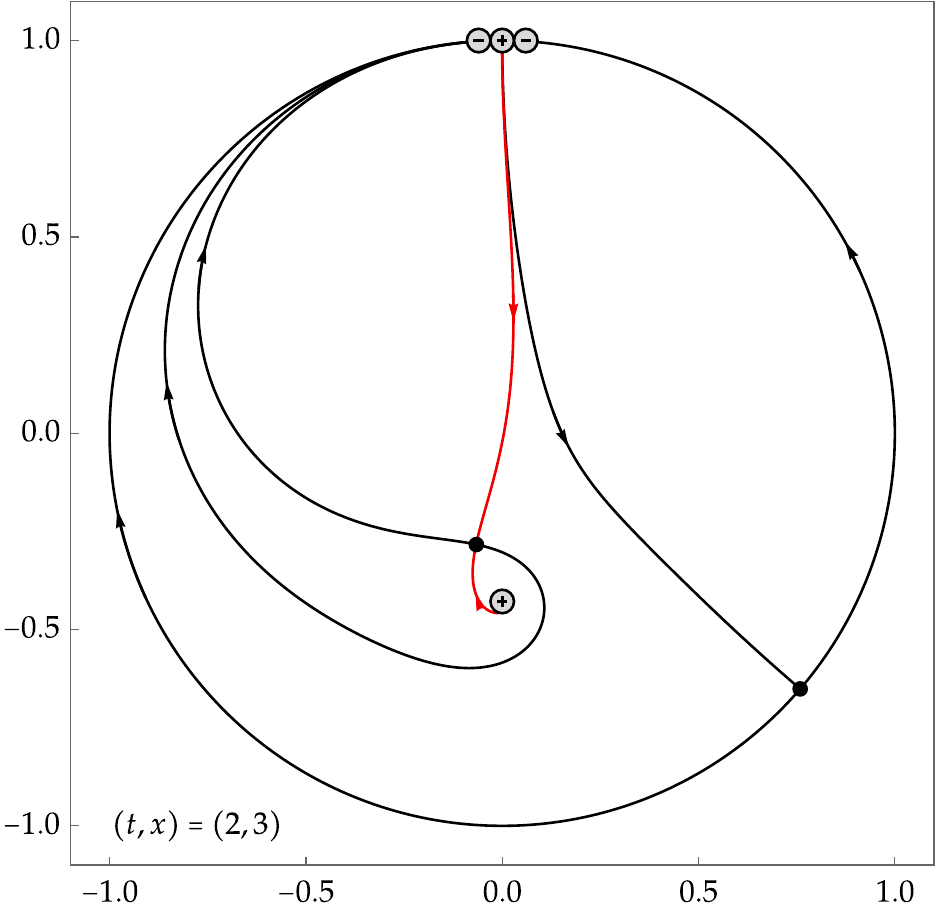}\hspace{0.03\linewidth}%
    \includegraphics[width=0.3\linewidth]{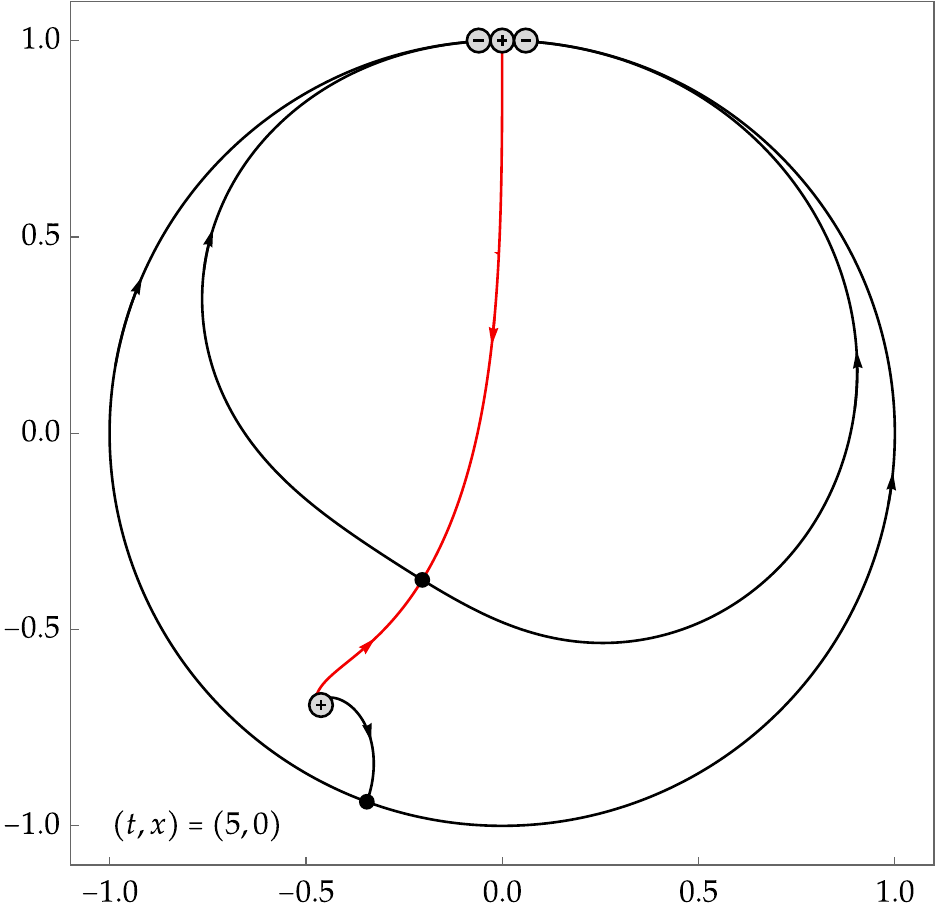}\hspace{0.03\linewidth}%
    \includegraphics[width=0.3\linewidth]{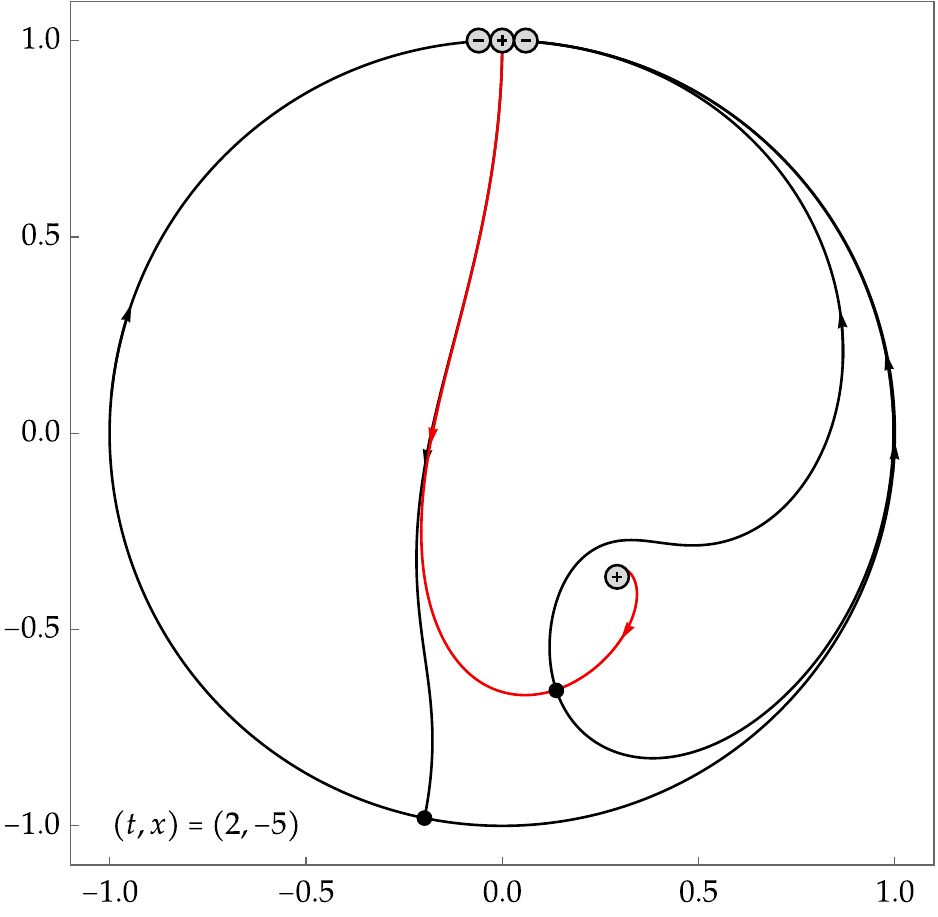}
\end{center}
\caption{Representative Stokes graphs with red indicating the edges belonging to $\mathcal{T}$ for $N=1$ and $J=0$ at the indicated values of $(t,x)$.  Top row ($p_1$ is a spiral sink):  $p_1=2\ii$, $c_1=-1+3\ii$ (left); $p_1=10+4\ii$, $c_1=-1+3\ii$.  Bottom row ($p_1$ is a spiral source):  $p_1=4\ii$, $c_1=1+3\ii$ (left); $p_1=-3+\ii$, $c_1=1+\ii$ (center); $p_1=3+4\ii$, $c_1=1+3\ii$ (right).  }
\label{fig:N1J0-Stokes}
\end{figure}

\begin{figure}
\begin{center}
    \includegraphics[width=0.3\linewidth]{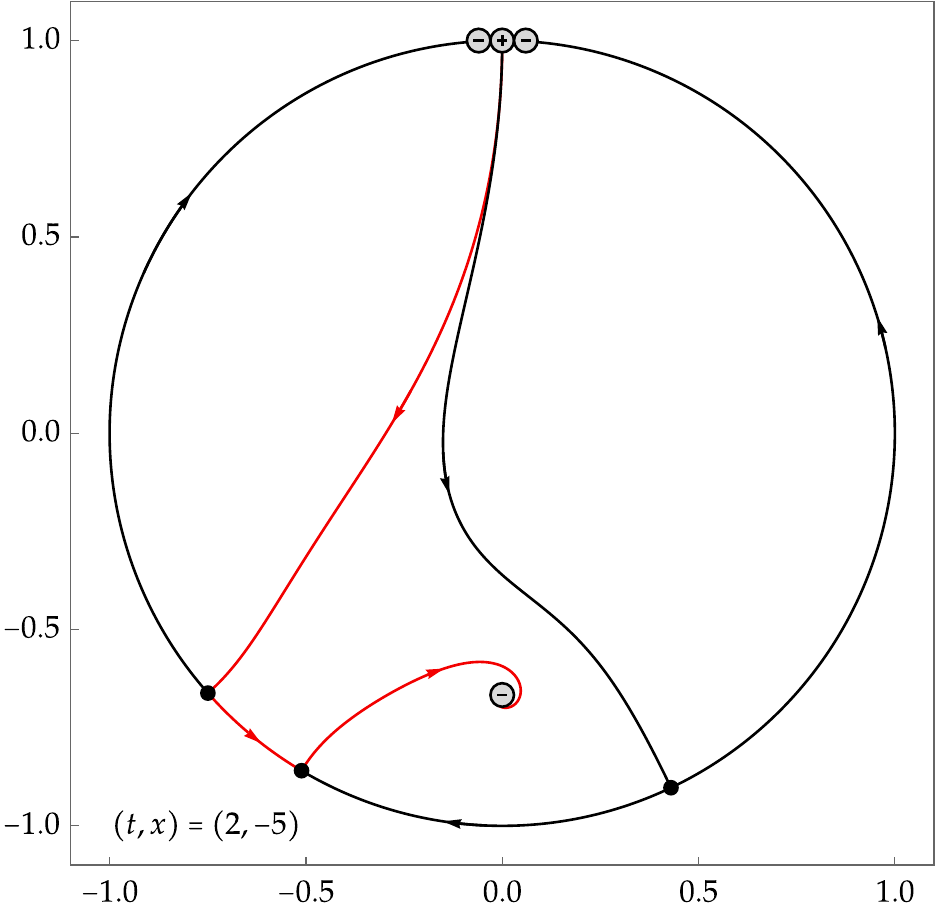}\hspace{0.03\linewidth}%
    \includegraphics[width=0.3\linewidth]{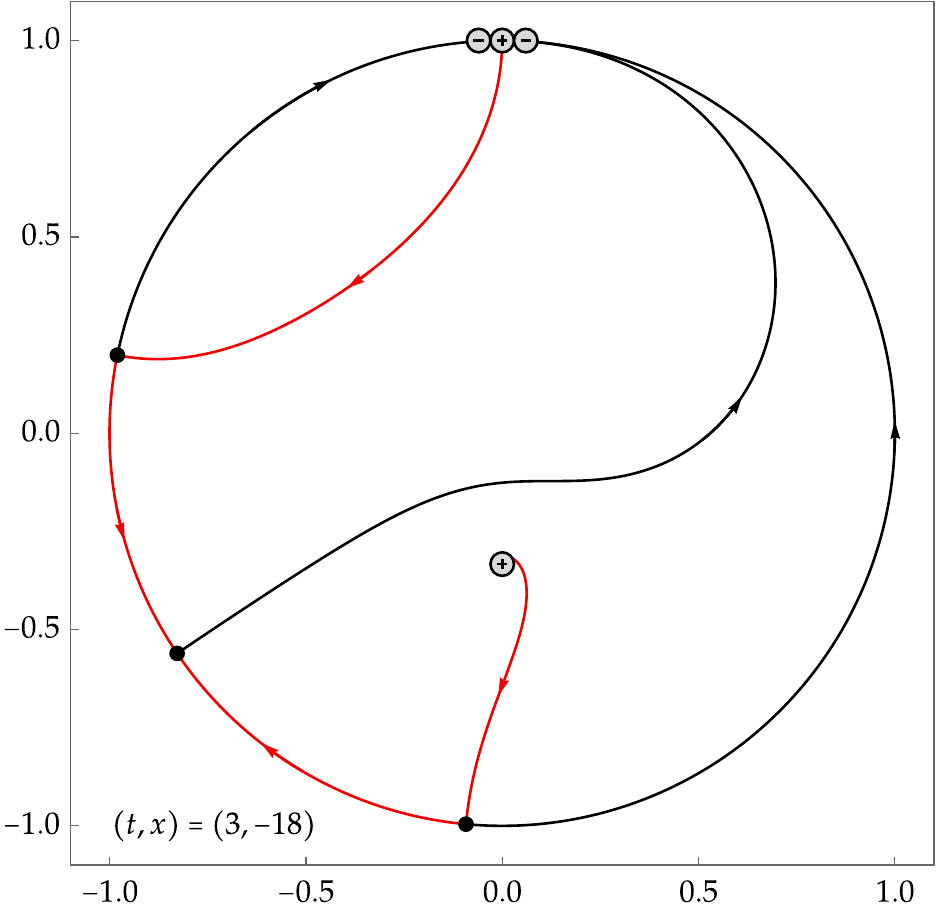}\hspace{0.03\linewidth}%
    \includegraphics[width=0.3\linewidth]{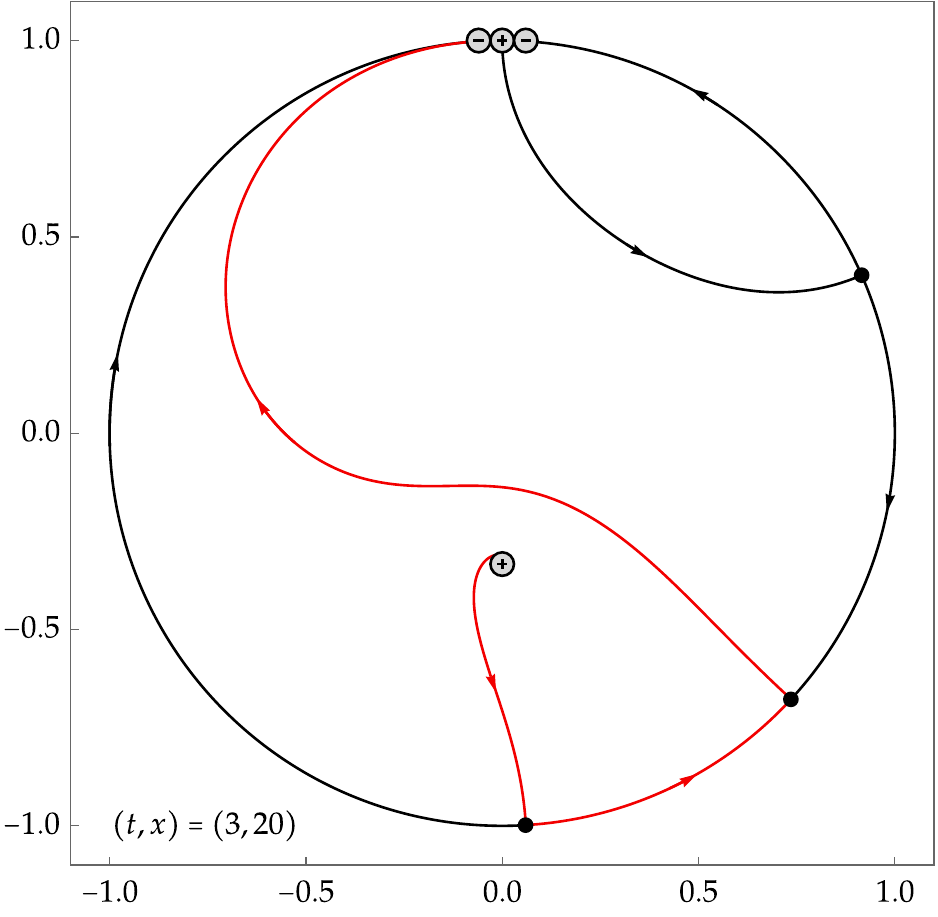}
\end{center}
\caption{Representative Stokes graphs with red indicating the edges belonging to $\mathcal{T}$ for $N=1$ and $J=1$ at the indicated values of $(t,x)$.  Left ($p_1$ is a spiral sink):  $p_1=2\ii$, $c_1=-1+3\ii$.  Center ($p_1$ is a spiral source):  $p_1=5\ii$, $c_1=4+7\ii$.  Right ($p_1$ is a spiral source):  $p_1=5\ii$, $c_1=4-8\ii$.}
\label{fig:N1J1-Stokes}
\end{figure}

To define $\mathcal{T}$, we first consider the case that $N=1$.  It turns out that there are eight distinct abstract Stokes graphs for $N=1$, five of type $(N,J)=(1,0)$ (see Figure~\ref{fig:N1J0-Stokes}) and three of type $(N,J)=(1,1)$ (see Figure~\ref{fig:N1J1-Stokes}).  

\paragraph{Case $J=0$ and $p_1$ is a sink}
The single saddle-type vertex on the boundary cycle of $\mathcal{S}$ is necessarily connected to all three component vertices of the compound vertex \inftynode, because its inward-directed edge is not part of the boundary cycle and cannot terminate at \sinknode.  That edge therefore divides the interior of the boundary cycle into two disjoint components.  Since \sinknode\ is certainly connected to a saddle point \saddlenode\ as $\mathcal{S}$ is connected, \saddlenode\ must be the degree-four saddle vertex $y_1$. Both \sinknode\ and \saddlenode\ lie within the same component, and the other component is an end-type face $\mathscr{E}$.  The three edges joined to \saddlenode\ other than that terminating at \sinknode\ necessarily terminate at \inftynode.  The four edges joined to \saddlenode\ subdivide the component containing \saddlenode\ into a second end-type face and two strip-type faces.  The two possible abstract Stokes graphs (depending on whether \sinknode\ and \saddlenode\ lie on the left or the right of the edge from \inftynode\ to the saddle on the boundary cycle) are represented by actual Stokes graphs for certain parameters as shown in the first row of Figure~\ref{fig:N1J0-Stokes}.  In these two cases, the edges selected to form $\mathcal{T}$ are (i) the edge from \saddlenode\ to \sinknode\ and (ii) the unique edge from \inftynode\ to \saddlenode\ that makes a left turn at \saddlenode\ to join edge (i).  These edges are colored red in the first row of Figure~\ref{fig:N1J0-Stokes}.

\paragraph{Case $J=0$ and $p_1$ is a source}  Then it is still possible that the saddle-type vertex $y_0$ on the boundary cycle is connected to all three components of \inftynode\ (see the left and right plots of second row in Figure~\ref{fig:N1J0-Stokes}), with the inward-directed edge separating the interior of the cycle into two disjoint components, exactly one of which must contain both the degree-four saddle \saddlenode\ and \sourcenode\ joined by an edge with three additional edges connected to \saddlenode\ terminating at \inftynode.

A third possibility when $p_1$ is a source \sourcenode\ is that the saddle-type vertex $y_0$ on the boundary cycle is connected in the interior of the cycle to the pole  \sourcenode\ representing $p_1$ (see the middle plot of the second row in Figure~\ref{fig:N1J0-Stokes}).  Since there must be at least one outward-directed edge from \inftynode, it can only terminate at the degree-four saddle vertex \saddlenode\ representing $y_1$. The remaining incoming edge to $y_1$ must come from $p_1$, thus the edges from $p_1$ to $y_0$, $p_1$ to $y_1$, and \inftynode~to $y_1$ divide the interior of $\mathcal{S}$ into two disjoint components, each of which must contain exactly one of the outward-directed edges from $y_1$ that necessarily terminate at \inftynode.  The abstract Stokes graph for this case is represented by the actual Stokes graph in the central panel of the second row in Figure~\ref{fig:N1J0-Stokes}.  To obtain the polytree $\mathcal{T}$ from $\mathcal{S}$ in all three cases that $(N,J)=(1,0)$ and $p_1$ is a source, we select the edge joining \sourcenode\ and the degree-four saddle vertex \saddlenode, as well as the other edge incoming to \saddlenode\ from \inftynode\ as shown in the three plots on the second row of Figure~\ref{fig:N1J0-Stokes}.

\paragraph{Case $J=1$} 
Next, suppose that $(N,J)=(1,1)$, in which case the only saddle-type vertices are on the boundary cycle.  We denote them $y_0$, $y_1$, and $y_2$ in clockwise order starting from \inftynode, with adjacent edges in the interior of the boundary cycle being directed toward $y_0$ and $y_2$ and away from $y_1$.  If $p_1$ is a sink vertex \sinknode, then it is necessarily joined to $y_1$, while $y_0$ and $y_2$ must be joined by edges in the interior of the boundary cycle to \inftynode.  This determines a unique abstract Stokes graph, that is represented by the actual Stokes graph shown in the left-hand panel of Figure~\ref{fig:N1J1-Stokes}.  The directed edges selected to construct $\mathcal{T}$ are:  \inftynode\ to $y_2$; $y_2$ to $y_1$; and $y_1$ to \sinknode.  If instead $p_1$ is a source vertex \sourcenode, then it is necessarily connected to either $y_0$ or $y_2$, with the remaining two saddle-type vertices on the boundary cycle joined to different components of the compound vertex \inftynode.  Thus we obtain two additional abstract Stokes graphs, concrete representatives of which are illustrated in the two right-hand panels of Figure~\ref{fig:N1J1-Stokes}.  In these two cases, we select edges for $\mathcal{T}$ by ensuring that $\mathcal{T}$ includes always the directed edge from $y_2$ to $y_1$ and that the edges of $\mathcal{T}$ join \sourcenode\ to either the left-hand sink-type vertex or the central source-type vertex at \inftynode.  The selected edges are shown in red in the corresponding panels of Figure~\ref{fig:N1J1-Stokes}.

\subsubsection{Inductive definition of the Stokes tree when \texorpdfstring{$N>1$}{N>1}}
Let $N>1$ be a given integer, let $J$ be an integer with $0\le J\le N$, and let $\mathcal{S}$ denote an arbitrary abstract Stokes graph with indices $(N,J)$.  Assume that a polytree $\mathcal{T}'$ has been extracted from every abstract Stokes graph $\mathcal{S}'$ with indices $(N-1,J')$ with $0\le J'\le N-1$.  We will construct the polytree $\mathcal{T}$ corresponding to $\mathcal{S}$ by induction on $N$.  

\paragraph{Case $J<N$} Then $\mathcal{S}$ contains a degree-four saddle-type vertex \saddlenode\ that is not in the boundary cycle.  There are either three or four distinct faces of $\mathcal{S}$ adjacent to \saddlenode.  

\begin{figure}
    \hfill\begin{tikzpicture}
     \node (source1) at (-0.8cm,0cm) [thick,black,circle,fill=lightgray,draw=black,inner sep=2pt]{\tiny $+$};
     \node (saddle) at (0cm,0cm) [circle,fill=black,minimum size=8pt,inner sep=0pt]{};
     \node (source2) at (0.8cm,0cm) [thick,black,circle,fill=lightgray,draw=black,inner sep=2pt]{\tiny $+$};
     \node (sink) at (1.6cm,0cm) [thick,black,circle,fill=lightgray,draw=black,inner sep=2pt]{\tiny $-$};
     \draw [->,>={Stealth[round]},thick,red] (source1) -- (-0.25cm,0cm);
     \draw [thick,red] (-0.25cm,0cm) -- (saddle);
     \draw [->,>={Stealth[round]},thick,red] (source2) -- (0.25cm,0cm);
     \draw [thick,red] (0.25cm,0cm) -- (saddle);
     \draw [->,>={Stealth[round]},thick,blue] (saddle) to[out=90,in=180] (0.8cm,0.8cm);
     \draw [thick,blue] (0.8cm,0.8cm) to[out=0,in=90] (sink);
     \draw [->,>={Stealth[round]},thick,blue] (saddle) to[out=-90,in=180] (0.8cm,-0.8cm);
     \draw [thick,blue] (0.8cm,-0.8cm) to[out=0,in=-90] (sink);
     \draw[thick] (source1) to[out=90,in=180] (-0.2cm,1cm);
     \draw[dashed,thick,->,>={Stealth[round]}] (-0.2cm,1cm) -- (0.5,1cm);
     \draw[dashed,thick] (0.5cm,1cm) -- (1cm,1cm);
     \draw[thick] (1cm,1cm) to[out=0,in=90] (sink);
     \draw[thick] (source1) to[out=-90,in=180] (-0.2cm,-1cm);
     \draw[dashed,thick,->,>={Stealth[round]}] (-0.2cm,-1cm) -- (0.5,-1cm);
     \draw[dashed,thick] (0.5cm,-1cm) -- (1cm,-1cm);
     \draw[thick] (1cm,-1cm) to[out=0,in=-90] (sink);
    \end{tikzpicture}\hfill%
    \begin{tikzpicture}
     \node (sink1) at (-0.8cm,0cm) [thick,black,circle,fill=lightgray,draw=black,inner sep=2pt]{\tiny $-$};
     \node (saddle) at (0cm,0cm) [circle,fill=black,minimum size=8pt,inner sep=0pt]{};
     \node (sink2) at (0.8cm,0cm) [thick,black,circle,fill=lightgray,draw=black,inner sep=2pt]{\tiny $-$};
     \node (source) at (1.6cm,0cm) [thick,black,circle,fill=lightgray,draw=black,inner sep=2pt]{\tiny $+$};
     \draw [->,>={Stealth[round]},thick,blue] (saddle) -- (-0.45cm,0cm);
     \draw [thick,blue] (-0.45cm,0cm) -- (sink1);
     \draw [->,>={Stealth[round]},thick,red] (saddle) -- (0.45cm,0cm);
     \draw [thick,red] (0.45cm,0cm) -- (sink2);
     \draw [->,>={Stealth[round]},thick,red] (source) to[out=90,in=0] (0.8cm,0.8cm);
     \draw [thick,red] (0.8cm,0.8cm) to[out=180,in=90] (saddle);
     \draw [->,>={Stealth[round]},thick,blue] (source) to[out=-90,in=0] (0.8cm,-0.8cm);
     \draw [thick,blue] (0.8cm,-0.8cm) to[out=180,in=-90] (saddle);
     \draw[thick] (source) to[out=90,in=0] (1cm,1cm);
     \draw[dashed,thick,->,>={Stealth[round]}] (1cm,1cm) -- (0.3,1cm);
     \draw[dashed,thick] (0.3cm,1cm) -- (-0.2cm,1cm);
     \draw[thick] (-0.2cm,1cm) to[out=180,in=90] (sink1);
     \draw[thick] (source) to[out=-90,in=0] (1cm,-1cm);
     \draw[dashed,thick,->,>={Stealth[round]}] (1cm,-1cm) -- (0.3,-1cm);
     \draw[dashed,thick] (0.3cm,-1cm) -- (-0.2cm,-1cm);
     \draw[thick] (-0.2cm,-1cm) to[out=180,in=-90] (sink1);
    \end{tikzpicture}\hfill
        \caption{Three strip-type faces adjacent to \saddlenode. Note that on the right, the vertex \sinknode\ opposite the saddle from the isolated pole might be the right-hand component of the compound vertex \inftynode, so the edge between it and the saddle should then not be selected as part of $\mathcal{T}$.}
        \label{fig:3strips}
\end{figure}
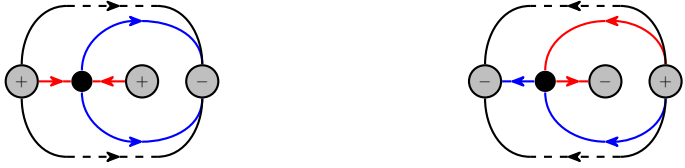

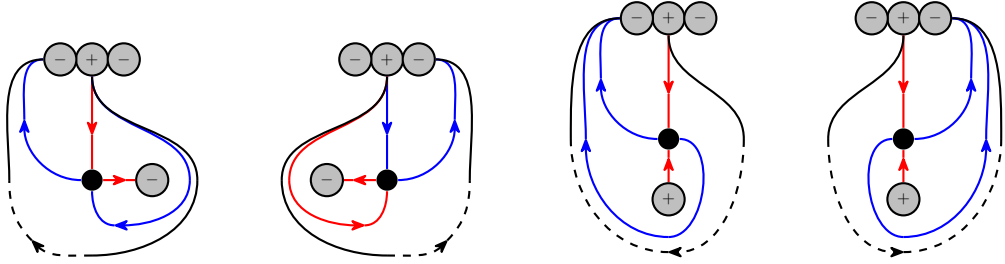
\begin{figure}
\hfill
\begin{tikzpicture}
\node (saddle) at (0cm,0cm) [circle,fill=black,minimum size=8pt,inner sep=0pt]{};
\node (sink) at (0.8cm,0cm) [thick,black,circle,fill=lightgray,draw=black,inner sep=2pt]{\tiny $-$};
\node (inftyC) at (0,1.6cm) [thick,black,circle,fill=lightgray,draw=black,inner sep=2pt] {\tiny $+$};
\node (inftyL) at (-12pt,1.6cm) [thick,black,circle,fill=lightgray,draw=black,inner sep=2pt] {\tiny $-$};
\node (inftyR) at (12pt,1.6cm) [thick,black,circle,fill=lightgray,draw=black,inner sep=2pt] {\tiny $-$};
\draw[->,>={Stealth[round]},thick,red] (inftyC) -- (0cm,0.6cm);
\draw[thick,red] (0cm,0.6cm) -- (saddle);
\draw[->,>={Stealth[round]},thick,red] (saddle) -- (0.45cm,0cm);
\draw[thick,red] (0.45cm,0cm) -- (sink);
\draw[->,>={Stealth[round]},thick,blue] (saddle) to[out=180,in=-90] (-0.9cm,0.8cm);
\draw[thick,blue] (-0.9cm,0.8cm) to[out=90,in=180] (inftyL);
\draw[thick,blue] (inftyC) to[out=-90,in=90] (1.3cm,0cm);
\draw[->,>={Stealth[round]},thick,blue] (1.3cm,0cm) to[out=-90,in=0] (0.3cm,-0.6cm);
\draw[thick,blue] (0.3cm,-0.6cm) to[out=180,in=-90] (saddle);
\draw[thick] (inftyC) to[out=-90,in=90] (1.4cm,0cm);
\draw[thick] (1.4cm,0cm) to[out=-90,in=0] (0cm,-1cm);
\draw[thick,dashed,->,>={Stealth[round]}] (0cm,-1cm) to[out=180,in=-45] (-0.8cm,-0.8cm);
\draw[thick,dashed] (-0.8cm,-0.8cm) to[out=135,in=-90] (-1.1cm,0cm);
\draw[thick] (-1.1cm,0cm) to[out=90,in=180] (inftyL);
\end{tikzpicture}\hfill%
\begin{tikzpicture}
\node (saddle) at (0cm,0cm) [circle,fill=black,minimum size=8pt,inner sep=0pt]{};
\node (sink) at (-0.8cm,0cm) [thick,black,circle,fill=lightgray,draw=black,inner sep=2pt]{\tiny $-$};
\node (inftyC) at (0,1.6cm) [thick,black,circle,fill=lightgray,draw=black,inner sep=2pt] {\tiny $+$};
\node (inftyL) at (-12pt,1.6cm) [thick,black,circle,fill=lightgray,draw=black,inner sep=2pt] {\tiny $-$};
\node (inftyR) at (12pt,1.6cm) [thick,black,circle,fill=lightgray,draw=black,inner sep=2pt] {\tiny $-$};
\draw[->,>={Stealth[round]},thick,blue] (inftyC) -- (0cm,0.6cm);
\draw[thick,blue] (0cm,0.6cm) -- (saddle);
\draw[->,>={Stealth[round]},thick,red] (saddle) -- (-0.45cm,0cm);
\draw[thick,red] (-0.45cm,0cm) -- (sink);
\draw[->,>={Stealth[round]},thick,blue] (saddle) to[out=0,in=-90] (0.9cm,0.8cm);
\draw[thick,blue] (0.9cm,0.8cm) to[out=90,in=0] (inftyR);
\draw[thick,red] (inftyC) to[out=-90,in=90] (-1.3cm,0cm);
\draw[->,>={Stealth[round]},thick,red] (-1.3cm,0cm) to[out=-90,in=180] (-0.3cm,-0.6cm);
\draw[thick,red] (-0.3cm,-0.6cm) to[out=0,in=-90] (saddle);
\draw[thick] (inftyC) to[out=-90,in=90] (-1.4cm,0cm);
\draw[thick] (-1.4cm,0cm) to[out=-90,in=180] (0cm,-1cm);
\draw[thick,dashed,->,>={Stealth[round]}] (0cm,-1cm) to[out=0,in=-135] (0.8cm,-0.8cm);
\draw[thick,dashed] (0.8cm,-0.8cm) to[out=45,in=-90] (1.1cm,0cm);
\draw[thick] (1.1cm,0cm) to[out=90,in=0] (inftyR);
\end{tikzpicture}\hfill%
\begin{tikzpicture}
\node (saddle) at (0cm,0cm) [circle,fill=black,minimum size=8pt,inner sep=0pt]{};
\node (inftyC) at (0,1.6cm) [thick,black,circle,fill=lightgray,draw=black,inner sep=2pt] {\tiny $+$};
\node (inftyL) at (-12pt,1.6cm) [thick,black,circle,fill=lightgray,draw=black,inner sep=2pt] {\tiny $-$};
\node (inftyR) at (12pt,1.6cm) [thick,black,circle,fill=lightgray,draw=black,inner sep=2pt] {\tiny $-$};    
\node (source) at (0cm,-0.8cm) [thick,black,circle,fill=lightgray,draw=black,inner sep=2pt]{\tiny $+$};
\draw[->,>={Stealth[round]},thick,red] (inftyC) -- (0cm,0.6cm);
\draw[thick,red] (0cm,0.6cm) -- (saddle);
\draw[->,>={Stealth[round]},thick,blue] (saddle) to[out=180,in=-90] (-0.9cm,0.8cm);
\draw[thick,blue] (-0.9cm,0.8cm) to[out=90,in=180] (inftyL);
\draw[->,>={Stealth[round]},thick,red] (source) -- (0cm,-0.25cm);
\draw[thick,red] (0cm,-0.25cm) -- (saddle);
\draw[thick,blue] (saddle) to[out=0,in=0] (0cm,-1.3cm);
\draw[->,>={Stealth[round]},thick,blue] (0cm,-1.3cm) to[out=180,in=-90] (-1.1cm,0cm);
\draw[thick,blue] (-1.1cm,0cm) to[out=90,in=180] (inftyL);
\draw[thick] (inftyC) to[out=-90,in=90] (1cm,0cm);
\draw[thick,dashed,->,>={Stealth[round]}] (1cm,0cm) to[out=-90,in=0] (0cm,-1.5cm);
\draw[thick,dashed] (0cm,-1.5cm) to[out=180,in=-90] (-1.3cm,0cm);
\draw[thick] (-1.3cm,0cm) to[out=90,in=180] (inftyL);
\end{tikzpicture}\hfill
\begin{tikzpicture}
\node (saddle) at (0cm,0cm) [circle,fill=black,minimum size=8pt,inner sep=0pt]{};
\node (inftyC) at (0,1.6cm) [thick,black,circle,fill=lightgray,draw=black,inner sep=2pt] {\tiny $+$};
\node (inftyL) at (-12pt,1.6cm) [thick,black,circle,fill=lightgray,draw=black,inner sep=2pt] {\tiny $-$};
\node (inftyR) at (12pt,1.6cm) [thick,black,circle,fill=lightgray,draw=black,inner sep=2pt] {\tiny $-$};    
\node (source) at (0cm,-0.8cm) [thick,black,circle,fill=lightgray,draw=black,inner sep=2pt]{\tiny $+$};
\draw[->,>={Stealth[round]},thick,red] (inftyC) -- (0cm,0.6cm);
\draw[thick,red] (0cm,0.6cm) -- (saddle);
\draw[->,>={Stealth[round]},thick,blue] (saddle) to[out=0,in=-90] (0.9cm,0.8cm);
\draw[thick,blue] (0.9cm,0.8cm) to[out=90,in=0] (inftyR);
\draw[->,>={Stealth[round]},thick,red] (source) -- (0cm,-0.25cm);
\draw[thick,red] (0cm,-0.25cm) -- (saddle);
\draw[thick,blue] (saddle) to[out=180,in=180] (0cm,-1.3cm);
\draw[->,>={Stealth[round]},thick,blue] (0cm,-1.3cm) to[out=0,in=-90] (1.1cm,0cm);
\draw[thick,blue] (1.1cm,0cm) to[out=90,in=0] (inftyR);
\draw[thick] (inftyC) to[out=-90,in=90] (-1cm,0cm);
\draw[thick,dashed,->,>={Stealth[round]}] (-1cm,0cm) to[out=-90,in=180] (0cm,-1.5cm);
\draw[thick,dashed] (0cm,-1.5cm) to[out=0,in=-90] (1.3cm,0cm);
\draw[thick] (1.3cm,0cm) to[out=90,in=0] (inftyR);
\end{tikzpicture}\hfill
    \caption{Two strip-type faces and one end-type face adjacent to \saddlenode.}
    \label{fig:2strips1end}
\end{figure}

First suppose that there are three distinct faces of $\mathcal{S}$ adjacent to \saddlenode.  Then exactly one of the faces is a strip-type face matching one of the two diagrams on the left of Figure~\ref{fig:S-faces-strips-enclosing}.  One or both of the remaining faces are strip-type faces from among those in Figure~\ref{fig:S-faces-strips-simple}, and if it is just one, then the other is one of the end-type faces from Figure~\ref{fig:S-faces-ends} having just one saddle-type vertex~\saddlenode. The possible configurations of the adjacent three domains are shown in Figures~\ref{fig:3strips} and~\ref{fig:2strips1end}.  In these diagrams, the dashed arcs indicate chains of edges with the indicated direction mediated by one or two saddle-type vertices forming part of the boundary of one of the strip-type faces shown in Figure~\ref{fig:S-faces-strips-simple}.  For each of the possible configurations from Figure~\ref{fig:3strips} we observe that upon erasing the red and blue edges as well as \saddlenode\ and the isolated source/sink-type vertex attached to it by a red edge, and keeping only the sequences of edges indicated by the dashed lines and the source and sink vertex forming their endpoints we are left with a single strip-type face from among the cases from Figure~\ref{fig:S-faces-strips-simple}.  We therefore obtain another abstract Stokes graph $\mathcal{S}'$ with indices $(N',J')=(N-1,J)$ by retaining all edges and vertices from $\mathcal{S}$ other than those erased from the relevant diagram in Figure~\ref{fig:3strips}.  Thus in $\mathcal{S}'$, three strip-type faces (two from the first three diagrams in Figure~\ref{fig:S-faces-strips-simple} and one from the first two diagrams in Figure~\ref{fig:S-faces-strips-enclosing}) have been replaced by a single strip-type face that could be any of the diagrams from Figure~\ref{fig:S-faces-strips-simple}.  Similarly, for each of the possible configurations from Figure~\ref{fig:2strips1end} we obtain $\mathcal{S}'$ with indices $(N',J')=(N-1,J)$ by erasing all nodes and edges bounded by the sequence of edges indicated by the dashed lines path from \inftynode\ to itself, hence replacing two strip-type faces and one end-type face with a single end-type face.  By the inductive hypothesis, we have the polytree $\mathcal{T}'$ for the abstract Stokes graph $\mathcal{S}'$, so to obtain $\mathcal{T}$ we simply include the edges shown in red in the diagrams from Figures~\ref{fig:3strips} and \ref{fig:2strips1end}.

\begin{figure}
\hfill\begin{tikzpicture}
\node (saddle) at (0cm,0cm) [circle,fill=black,minimum size=8pt,inner sep=0pt]{};
\node (sourceS) at (0cm,-1.6cm) [thick,black,circle,fill=lightgray,draw=black,inner sep=2pt]{\tiny $+$};
\node (sourceN) at (0cm,1.6cm) [thick,black,circle,fill=lightgray,draw=black,inner sep=2pt]{\tiny $+$};
\node (sinkW) at (-1.6cm,0cm) [thick,black,circle,fill=lightgray,draw=black,inner sep=2pt]{\tiny $-$};
\node (sinkE) at (1.6cm,0cm) [thick,black,circle,fill=lightgray,draw=black,inner sep=2pt]{\tiny $-$};
\draw[thick,red,->,>={Stealth[round]}] (sourceN) -- (0cm,0.6cm);
\draw[thick,red] (0cm,0.6cm) -- (saddle);
\draw[thick,red,->,>={Stealth[round]}] (sourceS) -- (0cm,-0.6cm);
\draw[red,thick] (0cm,-0.6cm) -- (saddle);
\draw[thick,blue,->,>={Stealth[round]}] (saddle) -- (-0.8cm,0cm);
\draw[thick,blue] (-0.8cm,0cm) -- (sinkW);
\draw[thick,blue,->,>={Stealth[round]}] (saddle) -- (0.8cm,0cm);
\draw[thick,blue] (0.8cm,0cm) -- (sinkE);
\draw[thick] (sourceN) -- (-0.6cm,1.6cm);
\draw[thick,dashed,->,>={Stealth[round]}] (-0.6cm,1.6cm) to[out=180,in=45] (-1.4cm,1.4cm);
\draw[thick,dashed] (-1.4cm,1.4cm) to[out=-135,in=-90] (-1.6cm,0.6cm);
\draw[thick] (-1.6cm,0.6cm) -- (sinkW);
\draw[thick] (sourceN) -- (0.6cm,1.6cm);
\draw[thick,dashed,->,>={Stealth[round]}] (0.6cm,1.6cm) to[out=0,in=135] (1.4cm,1.4cm);
\draw[thick,dashed] (1.4cm,1.4cm) to[out=-45,in=-90] (1.6cm,0.6cm);
\draw[thick] (1.6cm,0.6cm) -- (sinkE);
\draw[thick] (sourceS) -- (-0.6cm,-1.6cm);
\draw[thick,dashed,->,>={Stealth[round]}] (-0.6cm,-1.6cm) to[out=180,in=-45] (-1.4cm,-1.4cm);
\draw[thick,dashed] (-1.4cm,-1.4cm) to[out=135,in=-90] (-1.6cm,-0.6cm);
\draw[thick] (-1.6cm,-0.6cm) -- (sinkW);
\draw[thick] (sourceS) -- (0.6cm,-1.6cm);
\draw[thick,dashed,->,>={Stealth[round]}] (0.6cm,-1.6cm) to[out=0,in=-135] (1.4cm,-1.4cm);
\draw[thick,dashed] (1.4cm,-1.4cm) to[out=45,in=-90] (1.6cm,-0.6cm);
\draw[thick] (1.6cm,-0.6cm) -- (sinkE);
\end{tikzpicture}\hfill%
\begin{tikzpicture}
\node (inftyC) at (0,1.6cm) [thick,black,circle,fill=lightgray,draw=black,inner sep=2pt] {\tiny $+$};
\node (inftyL) at (-12pt,1.6cm) [thick,black,circle,fill=lightgray,draw=black,inner sep=2pt] {\tiny $-$};
\node (inftyR) at (12pt,1.6cm) [thick,black,circle,fill=lightgray,draw=black,inner sep=2pt] {\tiny $-$};    
\node (saddle) at (0cm,0cm) [circle,fill=black,minimum size=8pt,inner sep=0pt]{};
\node (source) at (0cm,-1.6cm) [thick,black,circle,fill=lightgray,draw=black,inner sep=2pt]{\tiny $+$};
\node (sink) at (1.6cm,0cm) [thick,black,circle,fill=lightgray,draw=black,inner sep=2pt]{\tiny $-$};
\draw[->,>={Stealth[round]},thick,blue] (saddle) to[out=180,in=-90] (-0.9cm,0.8cm);
\draw[thick,blue] (-0.9cm,0.8cm) to[out=90,in=180] (inftyL);
\draw[thick,red,->,>={Stealth[round]}] (inftyC) -- (0cm,0.6cm);
\draw[thick,red] (0cm,0.6cm) -- (saddle);
\draw[thick,red,->,>={Stealth[round]}] (source) -- (0cm,-0.6cm);
\draw[red,thick] (0cm,-0.6cm) -- (saddle);
\draw[thick,blue,->,>={Stealth[round]}] (saddle) -- (0.8cm,0cm);
\draw[thick,blue] (0.8cm,0cm) -- (sink);
\draw[thick] (source) -- (0.6cm,-1.6cm);
\draw[thick,dashed,->,>={Stealth[round]}] (0.6cm,-1.6cm) to[out=0,in=-135] (1.4cm,-1.4cm);
\draw[thick,dashed] (1.4cm,-1.4cm) to[out=45,in=-90] (1.6cm,-0.6cm);
\draw[thick] (1.6cm,-0.6cm) -- (sink);
\draw[thick] (source) -- (-1.1cm,-1.6cm);
\draw[thick,dashed,->,>={Stealth[round]}] (-1.1cm,-1.6cm) to[out=180,in=-90] (-1.6cm,0cm);
\draw[thick,dashed] (-1.6cm,0cm) to[out=90,in=180] (-1.1cm,1.6cm);
\draw[thick] (-1.1cm,1.6cm) -- (inftyL);
\draw[thick] (inftyC) to[out=-90,in=180] (0.3cm,0.8cm);
\draw[thick,dashed,->,>={Stealth[round]}] (0.3cm,0.8cm) -- (0.9cm,0.8cm);
\draw[thick,dashed] (0.9cm,0.8cm) -- (1.3cm,0.8cm);
\draw[thick] (1.3cm,0.8cm) to[out=0,in=90] (sink);
\end{tikzpicture}\hfill%
\begin{tikzpicture}
\node (inftyC) at (0,1.6cm) [thick,black,circle,fill=lightgray,draw=black,inner sep=2pt] {\tiny $+$};
\node (inftyL) at (-12pt,1.6cm) [thick,black,circle,fill=lightgray,draw=black,inner sep=2pt] {\tiny $-$};
\node (inftyR) at (12pt,1.6cm) [thick,black,circle,fill=lightgray,draw=black,inner sep=2pt] {\tiny $-$};    
\node (saddle) at (0cm,0cm) [circle,fill=black,minimum size=8pt,inner sep=0pt]{};
\node (source) at (0cm,-1.6cm) [thick,black,circle,fill=lightgray,draw=black,inner sep=2pt]{\tiny $+$};
\node (sink) at (-1.6cm,0cm) [thick,black,circle,fill=lightgray,draw=black,inner sep=2pt]{\tiny $-$};
\draw[->,>={Stealth[round]},thick,blue] (saddle) to[out=0,in=-90] (0.9cm,0.8cm);
\draw[thick,blue] (0.9cm,0.8cm) to[out=90,in=0] (inftyR);
\draw[thick,red,->,>={Stealth[round]}] (inftyC) -- (0cm,0.6cm);
\draw[thick,red] (0cm,0.6cm) -- (saddle);
\draw[thick,red,->,>={Stealth[round]}] (source) -- (0cm,-0.6cm);
\draw[red,thick] (0cm,-0.6cm) -- (saddle);
\draw[thick,blue,->,>={Stealth[round]}] (saddle) -- (-0.8cm,0cm);
\draw[thick,blue] (-0.8cm,0cm) -- (sink);
\draw[thick] (source) -- (-0.6cm,-1.6cm);
\draw[thick,dashed,->,>={Stealth[round]}] (-0.6cm,-1.6cm) to[out=0,in=-45] (-1.4cm,-1.4cm);
\draw[thick,dashed] (-1.4cm,-1.4cm) to[out=135,in=-90] (-1.6cm,-0.6cm);
\draw[thick] (-1.6cm,-0.6cm) -- (sink);
\draw[thick] (source) -- (1.1cm,-1.6cm);
\draw[thick,dashed,->,>={Stealth[round]}] (1.1cm,-1.6cm) to[out=0,in=-90] (1.6cm,0cm);
\draw[thick,dashed] (1.6cm,0cm) to[out=90,in=0] (1.1cm,1.6cm);
\draw[thick] (1.1cm,1.6cm) -- (inftyR);
\draw[thick] (inftyC) to[out=-90,in=0] (-0.3cm,0.8cm);
\draw[thick,dashed,->,>={Stealth[round]}] (-0.3cm,0.8cm) -- (-0.9cm,0.8cm);
\draw[thick,dashed] (-0.9cm,0.8cm) -- (-1.3cm,0.8cm);
\draw[thick] (-1.3cm,0.8cm) to[out=180,in=90] (sink);
\end{tikzpicture}\hfill%
\begin{tikzpicture}
\node (inftyC) at (0,1.6cm) [thick,black,circle,fill=lightgray,draw=black,inner sep=2pt] {\tiny $+$};
\node (inftyL) at (-12pt,1.6cm) [thick,black,circle,fill=lightgray,draw=black,inner sep=2pt] {\tiny $-$};
\node (inftyR) at (12pt,1.6cm) [thick,black,circle,fill=lightgray,draw=black,inner sep=2pt] {\tiny $-$};    
\node (saddle) at (0cm,0cm) [circle,fill=black,minimum size=8pt,inner sep=0pt]{};
\node (source) at (0cm,-1.6cm) [thick,black,circle,fill=lightgray,draw=black,inner sep=2pt]{\tiny $+$};
\draw[->,>={Stealth[round]},thick,blue] (saddle) to[out=0,in=-90] (0.9cm,0.8cm);
\draw[thick,blue] (0.9cm,0.8cm) to[out=90,in=0] (inftyR);
\draw[->,>={Stealth[round]},thick,blue] (saddle) to[out=180,in=-90] (-0.9cm,0.8cm);
\draw[thick,blue] (-0.9cm,0.8cm) to[out=90,in=180] (inftyL);
\draw[thick,red,->,>={Stealth[round]}] (inftyC) -- (0cm,0.6cm);
\draw[thick,red] (0cm,0.6cm) -- (saddle);
\draw[thick,red,->,>={Stealth[round]}] (source) -- (0cm,-0.6cm);
\draw[red,thick] (0cm,-0.6cm) -- (saddle);
\draw[thick] (source) -- (1.1cm,-1.6cm);
\draw[thick,dashed,->,>={Stealth[round]}] (1.1cm,-1.6cm) to[out=0,in=-90] (1.6cm,0cm);
\draw[thick,dashed] (1.6cm,0cm) to[out=90,in=0] (1.1cm,1.6cm);
\draw[thick] (1.1cm,1.6cm) -- (inftyR);
\draw[thick] (source) -- (-1.1cm,-1.6cm);
\draw[thick,dashed,->,>={Stealth[round]}] (-1.1cm,-1.6cm) to[out=180,in=-90] (-1.6cm,0cm);
\draw[thick,dashed] (-1.6cm,0cm) to[out=90,in=180] (-1.1cm,1.6cm);
\draw[thick] (-1.1cm,1.6cm) -- (inftyL);
\end{tikzpicture}\hfill
\caption{Four faces adjacent to \saddlenode.}
\label{fig:4faces}
\end{figure}
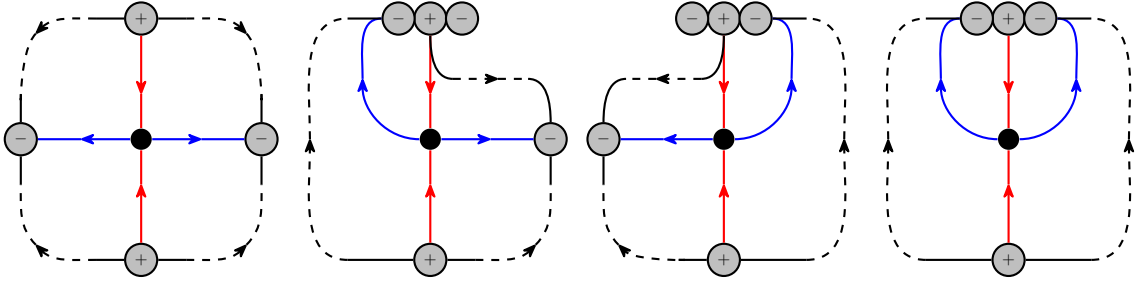

If instead there are four distinct faces of $\mathcal{S}$ adjacent to \saddlenode, then among these faces there can be zero, one, or two end-type faces (specifically, the left-hand diagram from either pair in Figure~\ref{fig:S-faces-ends}), and the remaining faces must be strip-type faces from among the first three of those shown in Figure~\ref{fig:S-faces-strips-simple}.
There are four distinct configurations possible, as shown in Figure~\ref{fig:4faces}.  
For each of these configurations, we pass from the abstract Stokes graph $\mathcal{S}$ to a smaller abstract Stokes graph $\mathcal{S}'$ by omitting the blue outward-directed edges from \saddlenode\ and then contracting along the red edges directed inward to \saddlenode\ so that the two source vertices labeled \sourcenode\ become identified and the vertex \saddlenode\ is also deleted.  Thus in each case, four faces enclosed by dashed arcs indicating chains of edges with one or two intervening saddle-type vertices are reduced to two faces bounded by the same edges, and the number and orientation of end-type faces among the initial four is preserved in the process.  Thus whereas $\mathcal{S}$ had indices $(N,J)$, $\mathcal{S}'$ is also an abstract Stokes graph and with one fewer saddle-type vertex \saddlenode\ not on the boundary cycle, and $\mathcal{S}'$ has indices $(N',J')=(N-1,J)$.  Invoking the inductive hypothesis, we obtain a polytree $\mathcal{T}'$ from $\mathcal{S}'$, and the identified vertex \sourcenode\ has a unique path in the tree to any other vertex of $\mathcal{T}'$. To obtain a polytree $\mathcal{T}$ subgraph of $\mathcal{S}$ from $\mathcal{T}'$, we once again separate the original two vertices that were identified into two distinct vertices of type \sourcenode\ only one of which belongs to $\mathcal{T}'$; then to obtain $\mathcal{T}$, we adjoin to $\mathcal{T}'$ the other vertex of type \sourcenode, the vertex \saddlenode, and the red edges of $\mathcal{S}$ joining them as shown in the diagrams of Figure~\ref{fig:4faces}.  

\paragraph{Case $J=N>1$} In this case, all saddle-type vertices are on the boundary cycle of $\mathcal{S}$, labeled in clockwise order from \inftynode\ as $y_0$, $y_1$, \dots, $y_{2J}$.  We first explain how the pair of vertices $y_{2J}$, $y_{2J-1}$ can be removed along with one of the $N$ source/sink vertices and some edges to result in an abstract Stokes graph $\mathcal{S}'$ with indices $(N',J')$ with $J'=N'=N-1$.  The inward-directed edge to $y_{2J}$ originates either 
\begin{enumerate}
    \item from the central source-type vertex at \inftynode, in which case the outward-directed edge from $y_{2J-1}$ terminates either
    \begin{enumerate}
        \item at the right-hand sink-type vertex at \inftynode, or \label{case1a}
        \item at a finite sink-type vertex \sinknode\ that is connected to another saddle-type vertex $y_j$ with $j<2J-2$ odd, or \label{case1b}
        \item at a finite sink-type vertex \sinknode\ that is not connected to any other vertices;
        \label{case1c}
    \end{enumerate}
    or
    \item  \label{case2} from a finite source-type vertex \sourcenode, in which case the outward-directed edge from $y_{2J-1}$ terminates either
    \begin{enumerate}
        \item at the left-hand sink-type vertex at \inftynode, or 
        \label{case2a}
        \item \label{case2b} at a finite sink-type vertex \sinknode\ that is  connected to another saddle-type vertex $y_j$ with $j<2J-2$ odd, or
        \item at a finite sink-type vertex \sinknode\ that is not connected to any other vertices.
        \label{case2c}
    \end{enumerate}
\end{enumerate}  
Note that in cases \ref{case2b} and \ref{case2c}, \sourcenode\ is  connected to at least one other saddle-type vertex $y_{j}$ with $j<2J-1$ even.

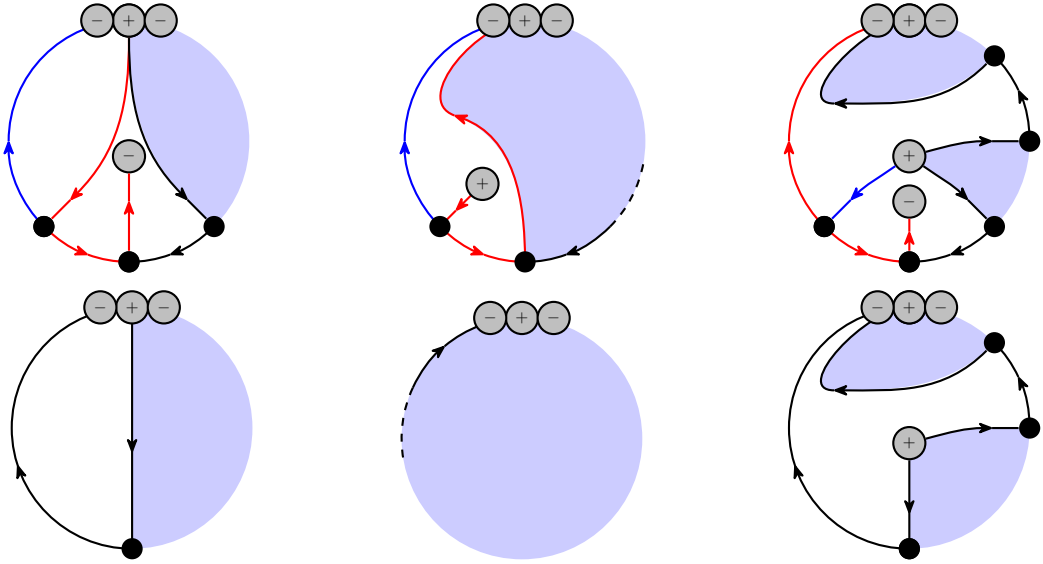
\begin{figure}
\begin{center}
\hfill\begin{tikzpicture}
\fill[blue!20!white] (315:1.6cm) arc[start angle=315,end angle=450,radius=1.6cm] to[out=-90,in=315-180] (315:1.1cm) -- cycle;
\node (inftyC) at (0,1.6cm) [thick,black,circle,fill=lightgray,draw=black,inner sep=2pt] {\tiny $+$};
\node (saddleC) at (0cm,-1.6cm) [circle,fill=black,minimum size=8pt,inner sep=0pt]{};
\node (saddleL) at (225:1.6cm) [circle,fill=black,minimum size=8pt,inner sep=0pt]{};
\node (saddleR) at (315:1.6cm) [circle,fill=black,minimum size=8pt,inner sep=0pt]{};
\node (sink) at (0cm,-0.2cm) [thick,black,circle,fill=lightgray,draw=black,inner sep=2pt]{\tiny $-$};
\draw[thick,blue,->,>={Stealth[round]}] (saddleL) arc[start angle=225,end angle=180,radius=1.6cm];
\draw[thick,blue] (180:1.6cm) arc[start angle=180,end angle=90,radius=1.6cm];
\draw[thick,red,->,>={Stealth[round]}] (saddleL) arc[start angle=225, end angle=250,radius=1.6cm];
\draw[thick,red] (250:1.6cm) arc[start angle=250,end angle=270,radius=1.6cm];
\draw[thick,->,>={Stealth[round]}] (saddleR) arc[start angle=315,end angle=290,radius=1.6cm];
\draw[thick] (290:1.6cm) arc[start angle=290,end angle=270,radius=1.6cm];
\draw[thick,red,->,>={Stealth[round]}] (saddleC) -- (0cm,-0.8cm);
\draw[thick,red] (0cm,-0.8cm) -- (sink);
\draw[thick,red,->,>={Stealth[round]}] (inftyC) to[out=-90,in=225-180] (225:1.1cm);
\draw[thick,red] (225:1.1cm) -- (saddleL);
\draw[thick,->,>={Stealth[round]}] (inftyC) to[out=-90,in=315-180] (315:1.1cm);
\draw[thick] (315:1.1cm) -- (saddleR);
\node (inftyC) at (0,1.6cm) [thick,black,circle,fill=lightgray,draw=black,inner sep=2pt] {\tiny $+$};
\node (inftyL) at (-12pt,1.6cm) [thick,black,circle,fill=lightgray,draw=black,inner sep=2pt] {\tiny $-$};
\node (inftyR) at (12pt,1.6cm) [thick,black,circle,fill=lightgray,draw=black,inner sep=2pt] {\tiny $-$};    
\node (saddleC) at (0cm,-1.6cm) [circle,fill=black,minimum size=8pt,inner sep=0pt]{};
\node (saddleL) at (225:1.6cm) [circle,fill=black,minimum size=8pt,inner sep=0pt]{};
\end{tikzpicture}\hfill%
\begin{tikzpicture}
\fill[blue!20!white] (-90:1.6cm) to[out=90,in=160-180] (160:1cm) to[out=160,in=180] (90:1.6cm) arc[start angle=90, end angle=-90,radius=1.6cm] -- cycle;
\draw[thick,red,->,>={Stealth[round]}] (-90:1.6cm) to[out=90,in=160-180] (160:1cm);
\draw[thick,red] (160:1cm) to[out=160,in=180] (90:1.6cm);
\draw[thick,red,->,>={Stealth[round]}] (225:1.6cm) arc[start angle=225, end angle=250,radius=1.6cm];
\draw[thick,red] (250:1.6cm) arc[start angle=250,end angle=270,radius=1.6cm];
\draw[thick,dashed] (315:1.6cm) arc[start angle=315, end angle=350,radius=1.6cm];
\draw[thick,->,>={Stealth[round]}] (315:1.6cm) arc[start angle=315,end angle=290,radius=1.6cm];
\draw[thick] (290:1.6cm) arc[start angle=290,end angle=270,radius=1.6cm];
\draw[thick,red,->,>={Stealth[round]}] (225:0.8cm)--(225:1.3cm);
\draw[thick,red] (225:1.3cm) -- (225:1.6cm);
\draw[thick,blue,->,>={Stealth[round]}] (225:1.6cm) arc[start angle=225,end angle=180,radius=1.6cm];
\draw[thick,blue] (180:1.6cm) arc[start angle=180,end angle=90,radius=1.6cm];
\node (inftyC) at (0,1.6cm) [thick,black,circle,fill=lightgray,draw=black,inner sep=2pt] {\tiny $+$};
\node (inftyL) at (-12pt,1.6cm) [thick,black,circle,fill=lightgray,draw=black,inner sep=2pt] {\tiny $-$};
\node (inftyR) at (12pt,1.6cm) [thick,black,circle,fill=lightgray,draw=black,inner sep=2pt] {\tiny $-$};    
\node (saddleC) at (0cm,-1.6cm) [circle,fill=black,minimum size=8pt,inner sep=0pt]{};
\node (saddleL) at (225:1.6cm) [circle,fill=black,minimum size=8pt,inner sep=0pt]{};
\node (source) at (225:0.8cm) [thick,black,circle,fill=lightgray,draw=black,inner sep=2pt] {\tiny $+$};
\end{tikzpicture}\hfill%
\begin{tikzpicture}
\fill[blue!20!white] (90:1.6cm) arc[start angle=90,end angle=45,radius=1.6cm] to[out=-135,in=0] (-1cm,0.5cm) 
to[out=180,in=180] cycle;
\fill[blue!20!white] (-90:0.2cm) to[out=15,in=180] (1.1cm,0cm) -- (0:1.6cm) arc[start angle=0,end angle=-45,radius=1.6cm] -- (315:1.1cm) to[out=135,in=180-35] cycle;
\node (inftyC) at (90:1.6cm) [thick,black,circle,fill=lightgray,draw=black,inner sep=2pt] {\tiny $+$};
\node (saddleC) at (0cm,-1.6cm) [circle,fill=black,minimum size=8pt,inner sep=0pt]{};
\node (saddleL) at (225:1.6cm) [circle,fill=black,minimum size=8pt,inner sep=0pt]{};
\node (saddleR) at (315:1.6cm) [circle,fill=black,minimum size=8pt,inner sep=0pt]{};
\node (saddleRR) at (0:1.6cm) [circle,fill=black,minimum size=8pt,inner sep=0pt]{};
\node (saddleRRR) at (45:1.6cm) [circle,fill=black,minimum size=8pt,inner sep=0pt]{};
\node (sink) at (0cm,-0.8cm) [thick,black,circle,fill=lightgray,draw=black,inner sep=2pt]{\tiny $-$};
\node (source) at (-90:0.2cm) [thick,black,circle,fill=lightgray,draw=black,inner sep=2pt] {\tiny $+$};
\draw[thick,->,>={Stealth[round]}] (source) to[out=15,in=180] (1.1cm,0cm);
\draw[thick] (1.1cm,0cm) -- (saddleRR);
\draw[thick,->,>={Stealth[round]}] (saddleRR) arc[start angle=0,end angle=25,radius=1.6cm];
\draw[thick] (25:1.6cm) arc [start angle=25,end angle=45,radius=1.6cm];
\draw[thick,->,>={Stealth[round]}] (saddleRRR) to[out=-135,in=0] (-1cm,0.5cm);
\draw[thick] (-1cm,0.5cm) to[out=180,in=180] (90:1.6cm);
\draw[thick,red,->,>={Stealth[round]}] (saddleL) arc[start angle=225,end angle=180,radius=1.6cm];
\draw[thick,red] (180:1.6cm) arc[start angle=180,end angle=90,radius=1.6cm];
\draw[thick,red,->,>={Stealth[round]}] (saddleL) arc[start angle=225, end angle=250,radius=1.6cm];
\draw[thick,red] (250:1.6cm) arc[start angle=250,end angle=270,radius=1.6cm];
\draw[thick,->,>={Stealth[round]}] (saddleR) arc[start angle=315,end angle=290,radius=1.6cm];
\draw[thick] (290:1.6cm) arc[start angle=290,end angle=270,radius=1.6cm];
\draw[thick,red,->,>={Stealth[round]}] (saddleC) -- (0cm,-1.2cm);
\draw[thick,red] (0cm,-1.2cm) -- (sink);
\draw[thick,blue,->,>={Stealth[round]}] (source) to[out=-145,in=225-180] (225:1.1cm);
\draw[thick,blue] (225:1.1cm) -- (saddleL);
\draw[thick,->,>={Stealth[round]}] (source) to[out=-35,in=315-180] (315:1.1cm);
\draw[thick] (315:1.1cm) -- (saddleR);
\node (inftyC) at (0,1.6cm) [thick,black,circle,fill=lightgray,draw=black,inner sep=2pt] {\tiny $+$};
\node (inftyL) at (-12pt,1.6cm) [thick,black,circle,fill=lightgray,draw=black,inner sep=2pt] {\tiny $-$};
\node (inftyR) at (12pt,1.6cm) [thick,black,circle,fill=lightgray,draw=black,inner sep=2pt] {\tiny $-$};    
\node (saddleC) at (0cm,-1.6cm) [circle,fill=black,minimum size=8pt,inner sep=0pt]{};
\node (saddleL) at (225:1.6cm) [circle,fill=black,minimum size=8pt,inner sep=0pt]{};
\end{tikzpicture}\hfill \\

\vspace{0.2cm}

\hfill\begin{tikzpicture}

\fill[blue!20!white] (90:1.6cm) arc[start angle=90, end angle=-90, radius=1.6cm];
\node (saddleC) at (0cm,-1.6cm) [circle,fill=black,minimum size=8pt,inner sep=0pt]{};
\draw[thick,postaction={midm arrow}] (-90:1.6cm) arc[start angle=-90, end angle=-270, radius=1.6cm];
\draw[thick,postaction={midp arrow}] (0,1.6) -- (0,-1.6);
\node (inftyC) at (0,1.6cm) [thick,black,circle,fill=lightgray,draw=black,inner sep=2pt] {\tiny $+$};
\node (inftyL) at (-12pt,1.6cm) [thick,black,circle,fill=lightgray,draw=black,inner sep=2pt] {\tiny $-$};
\node (inftyR) at (12pt,1.6cm) [thick,black,circle,fill=lightgray,draw=black,inner sep=2pt] {\tiny $-$};  

\end{tikzpicture}\hfill%
\begin{tikzpicture}

\fill[blue!20!white] (90:1.6cm) arc[start angle=90, end angle=-270, radius=1.6cm];

\draw[thick,dashed] (515:1.6cm) arc[start angle=515, end angle=550,radius=1.6cm];
\draw[thick,->,>={Stealth[round]}] (515:1.6cm) arc[start angle=515,end angle=490,radius=1.6cm];
\draw[thick] (490:1.6cm) arc[start angle=490,end angle=470,radius=1.6cm];

\node (inftyC) at (0,1.6cm) [thick,black,circle,fill=lightgray,draw=black,inner sep=2pt] {\tiny $+$};
\node (inftyL) at (-12pt,1.6cm) [thick,black,circle,fill=lightgray,draw=black,inner sep=2pt] {\tiny $-$};
\node (inftyR) at (12pt,1.6cm) [thick,black,circle,fill=lightgray,draw=black,inner sep=2pt] {\tiny $-$};  

\end{tikzpicture}\hfill%
\begin{tikzpicture}
\fill[blue!20!white] (90:1.6cm) arc[start angle=90,end angle=45,radius=1.6cm] to[out=-135,in=0] (-1cm,0.5cm) 
to[out=180,in=180] cycle;
\fill[blue!20!white] (-90:0.2cm) to[out=15,in=180] (1.1cm,0cm) -- (0:1.6cm) arc[start angle=0,end angle=-90,radius=1.6cm];
\node (inftyC) at (90:1.6cm) [thick,black,circle,fill=lightgray,draw=black,inner sep=2pt] {\tiny $+$};
\node (saddleC) at (0cm,-1.6cm) [circle,fill=black,minimum size=8pt,inner sep=0pt]{};
\node (saddleRR) at (0:1.6cm) [circle,fill=black,minimum size=8pt,inner sep=0pt]{};
\node (saddleRRR) at (45:1.6cm) [circle,fill=black,minimum size=8pt,inner sep=0pt]{};

\node (source) at (-90:0.2cm) [thick,black,circle,fill=lightgray,draw=black,inner sep=2pt] {\tiny $+$};
\draw[thick,->,>={Stealth[round]}] (source) to[out=15,in=180] (1.1cm,0cm);
\draw[thick] (1.1cm,0cm) -- (saddleRR);
\draw[thick,->,>={Stealth[round]}] (saddleRR) arc[start angle=0,end angle=25,radius=1.6cm];
\draw[thick] (25:1.6cm) arc [start angle=25,end angle=45,radius=1.6cm];
\draw[thick,->,>={Stealth[round]}] (saddleRRR) to[out=-135,in=0] (-1cm,0.5cm);
\draw[thick] (-1cm,0.5cm) to[out=180,in=180] (90:1.6cm);

\draw[thick,postaction={midm arrow}] (0,-1.6) arc[start angle=-90, end angle=-270, radius=1.6];

\draw[thick,postaction={midp arrow}] (source) to[out=-90,in=90] (0,-1.6);

\node (inftyC) at (0,1.6cm) [thick,black,circle,fill=lightgray,draw=black,inner sep=2pt] {\tiny $+$};
\node (inftyL) at (-12pt,1.6cm) [thick,black,circle,fill=lightgray,draw=black,inner sep=2pt] {\tiny $-$};
\node (inftyR) at (12pt,1.6cm) [thick,black,circle,fill=lightgray,draw=black,inner sep=2pt] {\tiny $-$};    
\node (saddleC) at (0cm,-1.6cm) [circle,fill=black,minimum size=8pt,inner sep=0pt]{};
\end{tikzpicture}\hfill
\end{center}
\caption{Configurations for $\mathcal{S}$ of type \ref{case1c} (left),  \ref{case2a} (center), and \ref{case2c} (right) with $N=J>1$ appear in the top row.  The rest of the abstract Stokes graph occupies the shaded region(s).  (In case \ref{case2c} the black edges shown emanating from \sourcenode\ may coincide.)  The corresponding reduced abstract Stokes graphs $\mathcal{S}'$ appear in the bottom row.}
\label{fig:ReduceJa}
\end{figure}

The cases \ref{case1c}, \ref{case2a}, and \ref{case2c} are the easiest to explain. These configurations are shown in Figure~\ref{fig:ReduceJa}.
In all three of these cases, there is a finite source or sink vertex joined only to $y_{2J}$ or $y_{2J-1}$ respectively and hence embedded in a strip-type face, specifically the third diagram from the left in Figure~\ref{fig:S-faces-strips-enclosing} for cases \ref{case1c} and \ref{case2c} and the second diagram from the right in Figure~\ref{fig:S-faces-strips-enclosing} for case \ref{case2a}.  In all three cases we obtain an abstract Stokes graph $\mathcal{S}'$ with indices $(N',J')$ for $J'=N'=N-1$ by essentially omitting the indicated face.  More precisely, we delete all edges colored red or blue in Figure~\ref{fig:ReduceJa}, the saddle-type vertices $y_{2J}$ and $y_{2J-1}$, and in cases \ref{case1c} and \ref{case2c} we omit the vertex \sinknode\ while in case \ref{case2a} we omit instead the vertex \sourcenode.  Then in all three cases we extend the edge that was incoming from the right to $y_{2J-1}$ to complete the boundary cycle by terminating at the left-hand sink-type vertex at \inftynode.  In case \ref{case1c} the end-type face on the left has thus been enlarged, while in case \ref{case2c} the strip-type face above the source vertex \sourcenode\ has been enlarged instead. All vertices and edges in the shaded parts of the abstract Stokes graphs are retained in $\mathcal{S}'$ as shown in the second row of Figure~\ref{fig:ReduceJa}.  Note that in case \ref{case2a}, the bounded face of $\mathcal{S}'$ that is adjacent to the edge shown entering the left sink-type vertex at \inftynode\ corresponds to a face of the same type (end or strip) in $\mathcal{S}$ that had two consecutive saddle-type vertices along its boundary (one of which was $y_{2J-1}$) instead of one; all other faces within the shaded part of $\mathcal{S}$ are exactly the same in $\mathcal{S}'$.
\begin{figure}
\hfill\begin{tikzpicture}


\fill[blue!20!white] (90:1.6cm) arc[start angle=90,end angle=-10,radius=1.6cm] to[out=170,in=0] (90:1.6cm);

\fill[blue!20!white] (0.2,-0.5) to[out=-60,in=130] (-50:1.6cm) -- (-50:1.6cm) arc[start angle=-50,end angle=-90,radius=1.6cm] -- (-90:1.6cm) to[out=90,in=-100] cycle;

\draw[thick,postaction={midm arrow}] (180:1.6cm) arc [start angle=180, end angle=90, radius=1.6cm];
\draw[thick,red,postaction={midp arrow}] (0,1.6) to[out=-90,in=0] (180:1.6cm);

\draw[thick,red,postaction={mid arrow}] (180:1.6cm) arc [start angle=180, end angle=220, radius=1.6cm];

\draw[thick,blue,postaction={midm arrow}] (220:1.6cm) to[out=40,in=0] (90:1.6);

\draw[thick,red,postaction={mid arrow}] (-90:1.6cm) arc [start angle=-90, end angle=-140, radius=1.6cm];

\draw[thick,postaction={midp arrow}] (-50:1.6cm) arc [start angle=-50, end angle=-10, radius=1.6cm];

\draw[thick,postaction={mid arrow}] (-10:1.6cm) to[out=170,in=0] (90:1.6);

\draw[thick,postaction={midp arrow}] (0.2,-0.5) to[out=-100,in=90] (-90:1.6cm);
\draw[thick,postaction={midp arrow}] (0.2,-0.5) to[out=-60,in=130] (-50:1.6cm);

\node (saddleL) at (180:1.6cm) [circle,fill=black,minimum size=8pt,inner sep=0pt]{};
\node (saddleC) at (220:1.6cm) [circle,fill=black,minimum size=8pt,inner sep=0pt]{};
\node (saddleR) at (-90:1.6cm) [circle,fill=black,minimum size=8pt,inner sep=0pt]{};
\node (saddleRR) at (-50:1.6cm) [circle,fill=black,minimum size=8pt,inner sep=0pt]{};
\node (saddleRRR) at (-10:1.6cm) [circle,fill=black,minimum size=8pt,inner sep=0pt]{};
\node (source) at (0.2,-0.5) [thick,black,circle,fill=lightgray,draw=black,inner sep=2pt] {\tiny $+$};

\node (inftyC) at (0,1.6cm) [thick,black,circle,fill=lightgray,draw=black,inner sep=2pt] {\tiny $+$};
\node (inftyL) at (-12pt,1.6cm) [thick,black,circle,fill=lightgray,draw=black,inner sep=2pt] {\tiny $-$};
\node (inftyR) at (12pt,1.6cm) [thick,black,circle,fill=lightgray,draw=black,inner sep=2pt] {\tiny $-$};

\end{tikzpicture}\hfill%
\begin{tikzpicture}

\fill[blue!20!white] (90:1.6) arc[start angle=90,end angle=-10,radius=1.6cm] to[out=170,in=-90] (90:1.6);

\fill[blue!20!white] (-50:1.6cm) arc[start angle=-50,end angle=-140,radius=1.6cm] to[out=40,in=-135] (-0.3,-0.4) to[out=-20,in=130] (-50:1.6cm);

\draw[thick,blue,postaction={midm arrow}] (180:1.6cm) arc [start angle=180, end angle=90, radius=1.6cm];
\draw[thick,red,postaction={midp arrow}] (0,1.6) to[out=-90,in=0] (180:1.6cm);

\draw[thick,red,postaction={mid arrow}] (180:1.6cm) arc [start angle=180, end angle=220, radius=1.6cm];

\draw[thick,red,postaction={midp arrow}] (220:1.6cm) to[out=40,in=-135] (-0.3,-0.4);

\draw[thick,postaction={midp arrow}] (-10:1.6cm) arc [start angle=-10, end angle=-50, radius=1.6cm];

\draw[thick,postaction={mid arrow}] (90:1.6cm) to[out=-90,in=170] (-10:1.6);

\draw[thick,postaction={midp arrow}] (-50:1.6cm) to[out=130,in=-20] (-0.3,-0.4);

\draw[thick,postaction={mid arrow}] (-110:1.6cm) arc [start angle=-110,end angle=-140,radius=1.6cm];
\draw[thick,dashed] (-90:1.6cm) arc [start angle=-90,end angle=-110,radius=1.6cm];

\node (saddleL) at (180:1.6cm) [circle,fill=black,minimum size=8pt,inner sep=0pt]{};
\node (saddleC) at (220:1.6cm) [circle,fill=black,minimum size=8pt,inner sep=0pt]{};
\node (saddleR) at (-50:1.6cm) [circle,fill=black,minimum size=8pt,inner sep=0pt]{};
\node (saddleRR) at (-10:1.6cm) [circle,fill=black,minimum size=8pt,inner sep=0pt]{};
\node (source) at (-0.3,-0.4) [thick,black,circle,fill=lightgray,draw=black,inner sep=2pt] {\tiny $-$};

\node (inftyC) at (0,1.6cm) [thick,black,circle,fill=lightgray,draw=black,inner sep=2pt] {\tiny $+$};
\node (inftyL) at (-12pt,1.6cm) [thick,black,circle,fill=lightgray,draw=black,inner sep=2pt] {\tiny $-$};
\node (inftyR) at (12pt,1.6cm) [thick,black,circle,fill=lightgray,draw=black,inner sep=2pt] {\tiny $-$};

\end{tikzpicture}\hfill%
\begin{tikzpicture}

\fill[blue!20!white] (90:1.6cm) arc[start angle=90,end angle=45,radius=1.6cm] to[out=-135,in=0] (-1cm,0.5cm) 
to[out=180,in=180] cycle;

\fill[blue!20!white] (0:1.6cm) arc[start angle=0,end angle=-50,radius=1.6cm] to[out=130,in=-30] (-0.5,0.1) 
to[out=10,in=180] (0:1.6cm);

\fill[blue!20!white] (-80:1.6cm) arc[start angle=-80,end angle=-140,radius=1.6cm] to[out=40,in=-135] (-0.3,-0.6) 
to[out=-75,in=100] (-80:1.6cm);

\draw[thick,red,postaction={midm arrow}] (-180:1.6cm) arc[start angle=-180,end angle=-270,radius=1.6cm];
\draw[thick,red,postaction={midp arrow}] (-180:1.6cm) arc[start angle=-180,end angle=-140,radius=1.6cm];

\draw[thick,postaction={midp arrow}] (-50:1.6cm) arc[start angle=-50,end angle=-80,radius=1.6cm];

\draw[thick,postaction={midp arrow}] (0:1.6cm) arc[start angle=0,end angle=45,radius=1.6cm];

\draw[thick,->,>={Stealth[round]}] (45:1.6cm) to[out=-135,in=0] (-1cm,0.5cm);
\draw[thick] (-1cm,0.5cm) to[out=180,in=180] (90:1.6cm);

\draw[thick,red,postaction={midp arrow}] (-140:1.6cm) to[out=40,in=-135] (-0.3,-0.6);
\draw[thick,postaction={midp arrow}] (-80:1.6cm) to[out=100,in=-75] (-0.3,-0.6);

\draw[thick,blue,postaction={midp arrow}] (-0.5,0.1) to[out=-170,in=0] (-180:1.6cm);
\draw[thick,postaction={midp arrow}] (-0.5,0.1) to[out=-30,in=130] (-50:1.6cm);
\draw[thick,postaction={midp arrow}] (-0.5,0.1) to[out=10,in=180] (0:1.6cm);

\draw[thick,postaction={mid arrow}] (-120:1.6cm) arc [start angle=-120,end angle=-140,radius=1.6cm];
\draw[thick,dashed] (-105:1.6cm) arc [start angle=-105,end angle=-120,radius=1.6cm];

\node (saddleL) at (-180:1.6cm) [circle,fill=black,minimum size=8pt,inner sep=0pt]{};
\node (saddleC) at (-140:1.6cm) [circle,fill=black,minimum size=8pt,inner sep=0pt]{};
\node (saddleR) at (-80:1.6cm) [circle,fill=black,minimum size=8pt,inner sep=0pt]{};
\node (saddleRR) at (-50:1.6cm) [circle,fill=black,minimum size=8pt,inner sep=0pt]{};
\node (saddleRRR) at (0:1.6cm) [circle,fill=black,minimum size=8pt,inner sep=0pt]{};
\node (saddleRRRR) at (45:1.6cm) [circle,fill=black,minimum size=8pt,inner sep=0pt]{};
\node (source) at (-0.5,0.1) [thick,black,circle,fill=lightgray,draw=black,inner sep=2pt] {\tiny $+$};
\node (sink) at (-0.3,-0.6) [thick,black,circle,fill=lightgray,draw=black,inner sep=2pt] {\tiny $-$};

\node (inftyC) at (0,1.6cm) [thick,black,circle,fill=lightgray,draw=black,inner sep=2pt] {\tiny $+$};
\node (inftyL) at (-12pt,1.6cm) [thick,black,circle,fill=lightgray,draw=black,inner sep=2pt] {\tiny $-$};
\node (inftyR) at (12pt,1.6cm) [thick,black,circle,fill=lightgray,draw=black,inner sep=2pt] {\tiny $-$};

\end{tikzpicture}\hfill \\

\hfill\begin{tikzpicture}

\fill[blue!20!white] (90:1.6cm) arc[start angle=90,end angle=-10,radius=1.6cm] to[out=170,in=0] (90:1.6cm);

\fill[blue!20!white] (-50:1.6cm) arc[start angle=-50,end angle=-130,radius=1.6cm] to[out=50,in=-90] (90:1.6cm) to[out=-90,in=130] (-50:1.6cm);

\draw[thick,postaction={midp arrow}] (-50:1.6cm) arc [start angle=-50, end angle=-10, radius=1.6cm];

\draw[thick,postaction={midm arrow}] (-130:1.6cm) arc [start angle=-130, end angle=-270, radius=1.6cm];

\draw[thick,postaction={midp arrow}] (90:1.6cm) to[out=-90,in=50] (-130:1.6cm);
\draw[thick,postaction={midp arrow}] (90:1.6cm) to[out=-90,in=130] (-50:1.6cm);

\draw[thick,postaction={mid arrow}] (-10:1.6cm) to[out=170,in=0] (90:1.6);

\node (saddleR) at (-130:1.6cm) [circle,fill=black,minimum size=8pt,inner sep=0pt]{};
\node (saddleRR) at (-50:1.6cm) [circle,fill=black,minimum size=8pt,inner sep=0pt]{};
\node (saddleRRR) at (-10:1.6cm) [circle,fill=black,minimum size=8pt,inner sep=0pt]{};

\node (inftyC) at (0,1.6cm) [thick,black,circle,fill=lightgray,draw=black,inner sep=2pt] {\tiny $+$};
\node (inftyL) at (-12pt,1.6cm) [thick,black,circle,fill=lightgray,draw=black,inner sep=2pt] {\tiny $-$};
\node (inftyR) at (12pt,1.6cm) [thick,black,circle,fill=lightgray,draw=black,inner sep=2pt] {\tiny $-$};

\end{tikzpicture}\hfill%
\begin{tikzpicture}


\fill[blue!20!white] (90:1.6) arc[start angle=90,end angle=-10,radius=1.6cm] to[out=170,in=-90] (90:1.6);

\fill[blue!20!white] (-50:1.6cm) arc[start angle=-50,end angle=-270,radius=1.6cm] to[out=-180,in=90] (160:1.2) to[out=-90,in=130] (-50:1.6cm);

\draw[thick,postaction={midp arrow}] (-10:1.6cm) arc [start angle=-10, end angle=-50, radius=1.6cm];

\draw[thick,postaction={mid arrow}] (90:1.6cm) to[out=-90,in=170] (-10:1.6);

\draw[thick,postaction={midp arrow}] (-50:1.6cm) to[out=130,in=-90] (160:1.2);
\draw[thick] (160:1.2) to[out=90,in=-180] (90:1.6);

\draw[thick,postaction={midm arrow}] (-220:1.6cm) arc [start angle=-220,end angle=-270,radius=1.6cm];
\draw[thick,dashed] (-190:1.6cm) arc [start angle=-190,end angle=-220,radius=1.6cm];

\node (saddleR) at (-50:1.6cm) [circle,fill=black,minimum size=8pt,inner sep=0pt]{};
\node (saddleRR) at (-10:1.6cm) [circle,fill=black,minimum size=8pt,inner sep=0pt]{};

\node (inftyC) at (0,1.6cm) [thick,black,circle,fill=lightgray,draw=black,inner sep=2pt] {\tiny $+$};
\node (inftyL) at (-12pt,1.6cm) [thick,black,circle,fill=lightgray,draw=black,inner sep=2pt] {\tiny $-$};
\node (inftyR) at (12pt,1.6cm) [thick,black,circle,fill=lightgray,draw=black,inner sep=2pt] {\tiny $-$};

\end{tikzpicture}\hfill%
\begin{tikzpicture}

\fill[blue!20!white] (90:1.6cm) arc[start angle=90,end angle=45,radius=1.6cm] to[out=-135,in=0] (-1cm,0.5cm) 
to[out=180,in=180] cycle;

\fill[blue!20!white] (0:1.6cm) arc[start angle=0,end angle=-50,radius=1.6cm] to[out=130,in=-30] (-0.5,0.1) 
to[out=10,in=180] (0:1.6cm);

\fill[blue!20!white] (-80:1.6cm) arc[start angle=-80,end angle=-270,radius=1.6cm] to[out=-180,in=90] (160:1.4) 
to[out=-90,in=100] (-80:1.6cm);

\draw[thick,postaction={midp arrow}] (-50:1.6cm) arc[start angle=-50,end angle=-80,radius=1.6cm];

\draw[thick,postaction={midp arrow}] (0:1.6cm) arc[start angle=0,end angle=45,radius=1.6cm];

\draw[thick,->,>={Stealth[round]}] (45:1.6cm) to[out=-135,in=0] (-1cm,0.5cm);
\draw[thick] (-1cm,0.5cm) to[out=180,in=180] (90:1.6cm);

\draw[thick,postaction={midm arrow}] (-80:1.6cm) to[out=100,in=-90] (160:1.4) to[out=90,in=-180] (90:1.6);

\draw[thick,postaction={midp arrow}] (-0.5,0.1) to[out=-30,in=130] (-50:1.6cm);
\draw[thick,postaction={midp arrow}] (-0.5,0.1) to[out=10,in=180] (0:1.6cm);

\draw[thick,postaction={midm arrow}] (-180:1.6cm) arc [start angle=-180,end angle=-270,radius=1.6cm];
\draw[thick,dashed] (-150:1.6cm) arc [start angle=-150,end angle=-180,radius=1.6cm];

\node (saddleR) at (-80:1.6cm) [circle,fill=black,minimum size=8pt,inner sep=0pt]{};
\node (saddleRR) at (-50:1.6cm) [circle,fill=black,minimum size=8pt,inner sep=0pt]{};
\node (saddleRRR) at (0:1.6cm) [circle,fill=black,minimum size=8pt,inner sep=0pt]{};
\node (saddleRRRR) at (45:1.6cm) [circle,fill=black,minimum size=8pt,inner sep=0pt]{};
\node (source) at (-0.5,0.1) [thick,black,circle,fill=lightgray,draw=black,inner sep=2pt] {\tiny $+$};

\node (inftyC) at (0,1.6cm) [thick,black,circle,fill=lightgray,draw=black,inner sep=2pt] {\tiny $+$};
\node (inftyL) at (-12pt,1.6cm) [thick,black,circle,fill=lightgray,draw=black,inner sep=2pt] {\tiny $-$};
\node (inftyR) at (12pt,1.6cm) [thick,black,circle,fill=lightgray,draw=black,inner sep=2pt] {\tiny $-$};
\end{tikzpicture}\hfill

\caption{Configurations of type \ref{case1a} (left),  \ref{case1b} (center), and \ref{case2b} (right) with $N=J>1$ appear in the top row.  The rest of the abstract Stokes graph occupies the shaded region(s). (In cases \ref{case1a} and \ref{case2b} the black edges shown emanating from \sourcenode\ may coincide.)  The corresponding reduced abstract Stokes graphs $\mathcal{S}'$ appear in the bottom row.}
\label{fig:ReduceJa2}
\end{figure}
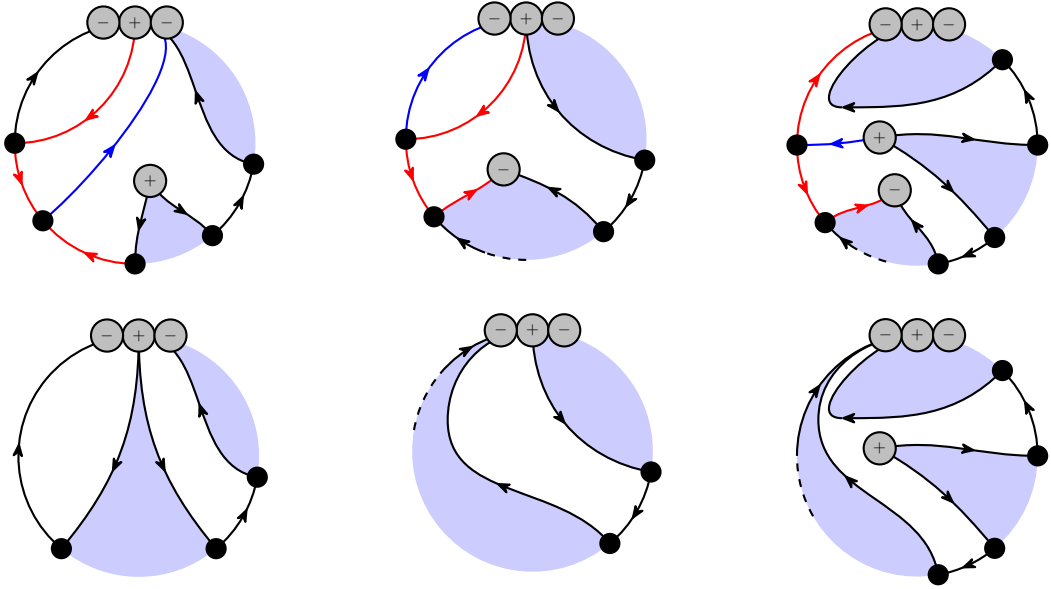

To obtain a reduced abstract Stokes graph $\mathcal{S}'$ with indices $(N-1,N-1)$ in cases \ref{case1a}, \ref{case1b}, and \ref{case2b} (see Figure \ref{fig:ReduceJa2}, first row), we first remove the red and blue colored edges and the saddle-type vertices $y_{2J}$, $y_{2J-1}$.  In case \ref{case1a}, we extend the edge from $y_{2J}$ (now deleted) to \inftynode\ backwards along the boundary cycle to $y_{2J-2}$ and contract the finite source node \sourcenode\ vertically to the source node at \inftynode\ while retaining all remaining edges.  In cases \ref{case1b} and \ref{case2b}, we contract the sink node \sinknode\ along the boundary cycle to the left sink node at \inftynode\ while retaining all remaining edges.  As before, all vertices and edges in the shaded parts of the abstract Stokes graphs are retained in $\mathcal{S}'$.  The resulting reduced abstract Stokes graphs are displayed in the second row of Figure \ref{fig:ReduceJa2}.  Note that in cases \ref{case1b} and \ref{case2b}, the bounded face of $\mathcal{S}'$ adjacent to the edge on the boundary cycle  entering the left sink-type vertex at \inftynode\ is a strip-type face whose corresponding strip-type face in $\mathcal{S}$ had two consecutive saddle-type vertices (one of which was $y_{2J-1}$) instead of one, while all other faces within the shaded part of $\mathcal{S}$ are exactly the same in $\mathcal{S}'$ (as in case \ref{case2a} shown in Figure~\ref{fig:ReduceJa}).

In all six cases of $N=J>1$, we invoke the inductive hypothesis to extract a polytree $\mathcal{T}'$ from the reduced abstract Stokes graph $\mathcal{S}'$, and then identify the edges of $\mathcal{T}'$ with specific edges of $\mathcal{S}$, to which we add also the red-colored edges in the top row diagrams of Figures~\ref{fig:ReduceJa}--\ref{fig:ReduceJa2} to obtain the polytree $\mathcal{T}$.

The edges of the polytree $\mathcal{T}$ correspond to actual critical trajectories of the quadratic differential $h'(z)^2\,\dd z^2$ in the $z$-plane, each oriented in the direction of increasing parameter $s$.  The union of the selected oriented critical trajectories will be called the \emph{Stokes tree} of $h'(z)$.  Thus, starting from each pole $p_n\in\mathbb{C}_+$ we may proceed along a well-defined sequence of trajectories toward the root, which means that the final trajectory in the sequence tends to $z=\infty$ either in the leftward horizontal direction or in the upward vertical direction.

\subsection{From a Stokes tree to suitable branch cuts of \texorpdfstring{$h(z)$}{h(z)} and paths of integration}
Let $p_n\in\mathbb{C}_+$, $n=1,\dots,N$, be one of the finite poles of $h'(z)$.  We now define a branch cut $\Gamma_n$ of $h(z)$ emanating from $z=p_n$ and a contour $W_n$ surrounding the branch cut as follows. Let $\delta>0$ be a small positive parameter. 
To each of the critical points $y$ of $h(z)$ we associate two nearby points $y_\alpha$ and $y_\beta$ defined by the property that for every locally-defined single-valued branch of $E(z)=-\ii h(z)$ we have $\mathrm{Im}(E(y_{\alpha,\beta}))=\mathrm{Im}(E(y))$ and $\mathrm{Re}(E(y_{\alpha,\beta}))=\mathrm{Re}(E(y))-\delta$.  The points $y_\alpha$ and $y_\beta$ are joined by a local steepest descent path $S_y$ passing from valley to valley over the saddle point $y$ and consisting of \emph{orthogonal} trajectory arcs.   For each Stokes tree trajectory along the path from $p_n$ back to the root  we define the noncritical trajectory arc $L$ (i.e., extended level curve) that lies on the right side of the Stokes tree trajectory by its orientation and terminates in each direction at either
\begin{itemize}
\item[(i)] one of the points $y_{\alpha,\beta}$ if $y$ is an endpoint of the Stokes tree trajectory;
\item[(ii)] the unique point where $|z-p|=\delta$ if the Stokes tree trajectory has a pole $p$ as an endpoint;
\item[(iii)] the unique point where $|z|=\delta^{-1}\gg 1$ if the Stokes tree trajectory tends to $z=\infty$.
\end{itemize}
Just as the trajectories of the Stokes tree can be ordered from $p_n$ back to the root, so can the corresponding noncritical trajectories $L$.

\begin{figure}
\centering
\begin{tikzpicture}

\node (p1) at (5,-2) [circle,fill=blue,minimum size=4pt,inner sep=0pt]{};
\node (p2) at (-1,0) [circle,fill=blue,minimum size=4pt,inner sep=0pt]{};

\node (s1) at (2,-0.5) [circle,fill=black,minimum size=4pt,inner sep=0pt]{};
\node (s2) at (-4,3) [circle,fill=black,minimum size=4pt,inner sep=0pt]{};

\node (sd1) at (2.3,-0.5) [circle,fill=red,minimum size=3pt,inner sep=0pt]{};
\node (sd2) at (1.7,-0.5) [circle,fill=red,minimum size=3pt,inner sep=0pt]{};
\draw[thick,red] (sd1) -- (sd2);

\node (sd3) at (-3.8,2.8) [circle,fill=red,minimum size=3pt,inner sep=0pt]{};
\node (sd4) at (-4.2,3.2) [circle,fill=red,minimum size=3pt,inner sep=0pt]{};
\draw[thick,red] (sd3) -- (sd4);

\draw[thick,red,postaction={midp arrow}] (sd4) to[out=180,in=-45] (-5.75,4);

\spiral[black,thick](5,-2)(90:3:0.5);
\draw[thick,postaction={midp arrow}] (5,-1.5) to[out=175,in=-45] (s1);

\bonusspiral[red,thick](5,-2)(0:90)(0.38:0.65)[1];
\draw[thick,red] (5.38,-2) -- (5,-2);
\node (intermediate1) at (5.38,-2) [circle,fill=red,minimum size=3pt,inner sep=0pt]{};
\draw[thick,red,postaction={midp arrow}] (5,-1.35) to[out=170,in=-45] (sd1);
 
\spiral[black,thick](-1,0)(-90:3:0.8);
\draw[thick,postaction={midp arrow}] (-1,-0.8) to[out=-10,in=135] (s1);

\spiral[black,thick](-1,0)(135:3:0.8);
\draw[thick,postaction={midp arrow}] (-1.57,0.57) to[out=-135,in=-90] (s2);

\draw[red,thick] (sd2) to[out=135,in=25] (0.4,-0.7);
\draw[red,thick,postaction={midp arrow}] (0.4,-0.7) to[out=-155,in=-4] (-1,-0.85);
\bonusspiral[red,thick](-1,0)(0:270)(0.65:0.85)[0];
\draw[red,thick] (-0.35,0) arc [start angle=0,end angle=-180,x radius=0.65cm,y radius=0.65cm];
\node (intermediate2) at (-0.35,0) [circle,fill=red,minimum size=3pt,inner sep=0pt]{};
\node (intermediate3) at (-1.65,0) [circle,fill=red,minimum size=3pt,inner sep=0pt]{};
\bonusspiral[red,thick](-1,0)(180:450)(0.65:0.86)[0];
\draw[red,thick] (-1,0.86) to[out=180,in=0] (-2.1,0.5);
\draw[red,thick,postaction={mid arrow}] (-2.1,0.5) .. controls (-2.9,0.6) and (-3.8,2) .. (sd3);


\node (t2) at (6,2) [circle,fill=white,minimum size=4pt,inner sep=0pt]{};
\node (t3) at (-3,-3) [circle,fill=white,minimum size=4pt,inner sep=0pt]{};

\node (t4) at (-3,4) [circle,fill=white,minimum size=4pt,inner sep=0pt]{};
\node (t5) at (-6,4) [circle,fill=white,minimum size=4pt,inner sep=0pt]{};
\node (t6) at (0,4) [circle,fill=white,minimum size=4pt,inner sep=0pt]{};

\draw[thick,black,postaction={midp arrow}] (t4) to[out=-135,in=90] (s2);
\draw[thick,black,postaction={midp arrow}] (s2) to[out=180,in=-45] (t5);
\draw[thick,black,postaction={midp arrow}] (s2) to[out=0,in=-90] (t6);

\draw[thick,postaction={midp arrow}] (s1) to[out=45,in=-135] (t2);

\draw[thick,postaction={midp arrow}] (s1) to[out=-135,in=45] (t3);

\node (p1) at (5,-2) [circle,fill=blue,minimum size=4pt,inner sep=0pt]{};
\node (p2) at (-1,0) [circle,fill=blue,minimum size=4pt,inner sep=0pt]{};

\node[red] at (5,-1.) {$L$};
\node[red] at (-0.8,-1.1) {$L_1$};
\node[red] at (-2.7,1.2) {$L_2$};
\node[red] at (-5,3.7) {$L$};

\node[red] at (2.7,-0.4) {$y_\alpha$};
\node[red] at (1.4,-0.6) {$y_\beta$};

\node[red] at (-3.5,2.7) {$y_\alpha$};
\node[red] at (-4.3,3.4) {$y_\beta$};

\node[blue] at (4.3,-2.4) {$p_n$};
\draw[blue,dotted,thick] (4.4,-2.3) -- (5,-2);

\node[blue] at (-2,-0.4) {$p$};
\draw[blue,dotted,thick] (-1.85,-0.3) -- (-1,0);

\end{tikzpicture}
\caption{The union of the red curves is the branch cut $\Gamma_n$, emanating from the pole $p_n$, as described by Definition~\ref{def:branchcuts}. The red points separate the different concatenated arcs of $\Gamma_n$.  The black points denote complex critical points (in the upper half plane) and the black curves are the critical trajectories.}
\label{fig:branch-cut}
\end{figure}
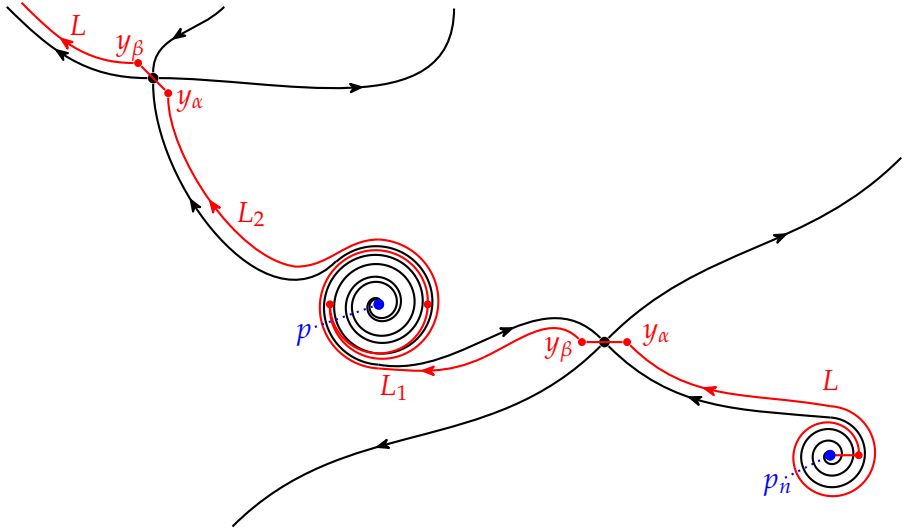

\begin{definition}[Choice of branch cuts]\label{def:branchcuts}
  We define the branch cut $\Gamma_n$ as the concatenation of:
\begin{itemize}
    \item all trajectory arcs $L$ corresponding to critical trajectories along the path in the Stokes tree from $p_n$ back to the root;
    \item the line segment from $p_n$ to the endpoint $z$ with $|z-p_n|=\delta$ of the first arc $L$ along the route back to the root;
    \item all local steepest descent paths $S_y$ whose endpoints $y_\alpha$ and $y_\beta$ are also endpoints of consecutive arcs $L$;
   \item the orthogonal trajectory $S_\infty$ from the large endpoint $z$ with $|z|=\delta^{-1}$ of the last arc $L$ toward $z=\infty$ in the direction of steepest decrease of $\mathrm{Re}(E(z))$ (this arc is asymptotic to the ray $\arg(z)=3\pi/4$);
   \item for each finite pole $p$ different from $p_n$ that is an intermediate endpoint of a Stokes tree trajectory on the path back to the root, a circular arc $A$ with $|z-p|=\delta$ from the terminal endpoint $z=z_1$ of the preceding arc $L_1$ to the initial endpoint $z=z_2$ of the subsequent arc $L_2$,  oriented in the clockwise (resp., counterclockwise) direction if $p$ is a source (resp., sink).
\end{itemize}
\end{definition}
See Figure~\ref{fig:branch-cut} for an example of the cut $\Gamma_n$ and Figure~\ref{fig:branch-cut N=1} for an example of the steepest descent contours $W_0$, $W_1$ and the branch cut $\Gamma_1$ in the $N=1$ setting corresponding to the middle pane of Figure~\ref{fig:N1J1-Stokes}.  With the branch cuts defined, we can now construct the steepest descent contours.

\begin{definition}[Choice of steepest descent contours]
For $n=1,\dots,N$, the contour $W_n$ is formed from arcs lying immediately against and on either side of the branch cut $\Gamma_n$ except for the straight-line segment of length $\delta$ from $p_n$; we then complete $W_n$ with a positively-oriented circular arc denoted $\Omega_n$ centered at $z=p_n$ of radius $\delta$ and subtending $2\pi$ radians with initial and terminal endpoints on opposite sides of $\Gamma_n$.  

The contour $W_0$ is defined as the concatenation of:
\begin{itemize}
    \item all local steepest descent paths $S_{y_{j}}$ for real critical points $y_j$ with $j=0,\dots,2J$;
    \item noncritical trajectory arcs $L_j$ joining consecutive local steepest descent paths $S_{y_j}$ and $S_{y_{j+1}}$, $j=0,\dots,2J-1$ and two additional noncritical trajectory arcs $L_\infty^\pm$ with $L_\infty^-$ joining $S_{y_{2J}}$ with a point $z^-$ near the negative real line with $|z^-|=\delta^{-1}$ (and $\mathrm{Im}(z^-)>0$), and $L_\infty^+$ joining $S_{y_0}$ with a point $z^+$ near the positive real line with $|z^+|=\delta^{-1}$ (and $\mathrm{Im}(z^+)<0$);
    \item the orthogonal trajectories $S_\infty^\pm$ from $z^\pm$ to $z=\infty$ in the direction of steepest decrease of $\mathrm{Re}(E(z))$ ($S_\infty^+$ asymptotic to the ray $\arg(z)=-\pi/4$ and $S_\infty^-$ asymptotic to the ray $\arg(z)=3\pi/4$).
\end{itemize}
\label{def:W}
\end{definition}

If two different poles $p_j\neq p_k$ correspond to vertices in the polytree $\mathcal{T}$ having a common ancestor (allowing for a vertex to be its own ancestor by definition) other than the root, then by the above construction the branch cuts $\Gamma_j$ and $\Gamma_k$ will have arcs in common.  In this situation, we appeal to the reasoning indicated in Remark~\ref{rem:move-branch-cuts} to place the common arcs of the cuts directly against each other along with the arcs of their ``wrapping'' contours.  The order (left-to-right) of the abutting arcs is important to consider to avoid intersecting cuts so that the value of $h(z)$ is well-defined along each $W_n$.  A similar observation applies to $W_0$ which lies along some arcs of branch cuts $\Gamma_1,\dots,\Gamma_N$ if $J>0$.

\begin{proof}[Proof of Proposition~\ref{prop:W}]
    The construction of the contours $W_0,\dots,W_N$ assumed that $u_0$ and $(t,x)$ were in nonspecial configuration, so we begin by giving the proof in that case.

    First we establish the indicated $1-1$ correspondence.  For $n=1,\dots,N$, let $W_n$ be the path enclosing the branch cut $\Gamma_n$ emanating from the pole $p_n$.  Note that the first Stokes tree trajectory from $p_n$ back towards the root necessarily has a critical point $y$ as its other endpoint.  
    
    If $\mathrm{Im}(y)>0$ then $\Gamma_n$ contains the steepest descent arc $S_y$ and we assign $y$ to $W_n$.  
    The only other Stokes tree trajectory with endpoint $y$ is closer to the root, so $y$ cannot be assigned to any other contour $W_j$, $j\neq n$.  
    
    If $y\in\mathbb{R}$ and the local steepest descent path $S_y$ is an arc of $\Gamma_n$, then $y$ is one of the pair $(y_{2j},y_{2j-1})$ for some $j=1,\dots,J$ and we assign that pair to $W_n$.  If $y\in\mathbb{R}$ but $\Gamma_n$ does not contain $S_y$, then the oriented arcs of the Stokes tree form a right turn at $y$, so $p_n$ is a source and $y=y_{2j-2}$ for some $j=1,\dots,J$; in this case we assign the pair $(y_{2j},y_{2j-1})$ to $W_n$.  In both cases, the branch cut $\Gamma_n$ contains both steepest descent arcs $S_{y_{2j}}$ and $S_{y_{2j-1}}$ as well as an intervening noncritical trajectory arc $L\subset\mathbb{C}_-$ close to the real Stokes tree trajectory joining the pair of critical points; as any other poles $p_j$ with $j\neq n$ that are endpoints of  Stokes tree trajectories along the path back from $p_n$ are closer to the root, the pair $(y_{2j},y_{2j-1})$ cannot be assigned to any other $W_j$, $j\neq n$.  

    Now we consider how $\mathrm{Re}(E(z))$ varies along $W_n$.  First note that $\mathrm{Re}(E(z))$ is strictly monotone along the circular arc $\Omega_n$ centered at $p_n$, so $W_n$ is the concatenation of $W_n^-$, $\Omega_n$ and $W_n^+$ (where the components $W_n^\pm$ lie along the same arcs of $\Gamma_n$ but on opposite sides).  If $\Delta_n>0$ denotes the absolute value of the increment of $\mathrm{Re}(E(z))$ along $\Omega_n$, and if $z^\pm$ denote the points of $W_n^\pm$ that are adjacent to the same arbitrary point $z\in \Gamma_n$, then the superscript notation on $W_n^\pm$ indicates that $\mathrm{Re}(E(z^+))-\mathrm{Re}(E(z^-))=\Delta_n$. This implies that any maxima of $\mathrm{Re}(E(z))$ along $W_n$ necessarily occur along $W_n^+$.  Next we observe that $\mathrm{Re}(E(z))$ is constant along all arcs of $W_n^+$ lying against trajectory arcs $L$ that are components of $\Gamma_n$.  Also, $\mathrm{Re}(E(z))$ is strictly decreasing in the direction away from $p_n$ towards $z=\infty$ along the final arc of $W_n^+$ adjacent to $S_\infty$ and along the arcs adjacent to the circular arcs $A$ for any intervening poles $p_j$, $j\neq n$.  Since on each local steepest descent path $S_y$ that is part of the branch cut $\Gamma_n$ the value of $\mathrm{Re}(E(z))$ increases by $\delta>0$ from both endpoints to achieve a local maximum at $z=y$, the maximum value of $\mathrm{Re}(E(z))$ is achieved on $W_n^+$ exactly at those critical points $y$ lying between $\Omega_n$ and the first circular arc $A$ or $S_\infty$ on the way to $z=\infty$.  This is either a unique complex critical point or a sequence of pairs $(y_{2j},y_{2j-1})$, $j=1,\dots,J$ of real critical points.

    Finally, similar arguments prove that the maximum of $\mathrm{Re}(E(z))$ is achieved at each of the real critical points of $h(z)$.  This completes the proof of the proposition in the case that $u_0$ and $(t,x)$ are in nonspecial configuration.

    Now assume to the contrary, but that $(t,x)\not\in\mathcal{D}$ so that all critical points, real or complex, of $h(z)$ are simple. First suppose that $\mathrm{Re}(c_n)\neq 0$ for all $n=1,\dots,N$, in which case $c_n\not\in\ii\epsilon\mathbb{N}$ holds for all $n$ and all $\epsilon>0$.  The contours $W_1,\dots,W_N$ then are loops encircling the individual branch cuts, but to specify the contours we cannot use the same procedure because the extended level curves can join critical points to themselves or to other critical points.  However, since the inequalities $\mathrm{Im}(c_n)\neq 0$ for $n=1,\dots,N$, \eqref{eq:no-homoclinics}, and \eqref{eq:no-heteroclinics} are open conditions on the parameters, there is a nearby configuration of $u_0$ and $(t,x)$ that is arbitrarily close in Euclidean norm on the parameters $\{p_n\}_{n=1}^N$ and $\{c_n\}_{n=1}^N$ of $u_0$ and $(t,x)\in\mathbb{R}^2$ that is nonspecial.  Let $\widetilde{\Gamma}_1,\dots,\widetilde{\Gamma}_N$ denote the branch cuts and $\widetilde{W}_0,\dots,\widetilde{W}_N$ denote the corresponding contours for that configuration.  Since $\mathrm{Re}(E(z))$ is maximal for the nonspecial configuration only at the critical points, and since these are all simple and distinct, we can take $\Gamma_n=\widetilde{\Gamma}_n$ and $W_n=\widetilde{W}_n$ except in a small neighborhood of the critical points.  So, locally within a small disk centered near each critical point $\widetilde{y}$ we simply replace the local steepest descent path $\widetilde{S}_{\widetilde{y}}$ of the nearby configuration with that of the actual configuration, $S_y$ traversing $y$, tracing the boundary of the disk by some small angles to maintain connectivity.  The branch cuts and integration contours that followed the local steepest descent path for the nearby configuration are then replaced by corresponding cuts and contours following the actual local path to complete the definition of $\Gamma_1,\dots,\Gamma_N$ and $W_0,\dots,W_N$.

    If there is a value of $n$ such that $\mathrm{Re}(c_n)=0$ and $\mathrm{Im}(c_n)>0$, then an additional modification is needed for those values of $\epsilon$ where $c_n\in\ii\epsilon\mathbb{N}$, because the actual contour $W_n$ must terminate at the pole $p_n$ while for an arbitrarily close nonspecial configuration $\widetilde{W}_n$ is a loop surrounding $p_n$.  However, in this case to define $\Gamma_n$ and $W_n$ away from the corresponding critical point(s), it suffices to shrink the circle $\widetilde{\Omega}_n$ arc of $\widetilde{W}_n$ to a point (because $\ee^{-\ii h(z)/\epsilon}$ is certainly integrable at $z=p_n$ for the nearby configuration) and then simply to omit from $W_n$ the half of $\widetilde{W}_n$ on which the corresponding critical point (points) is (are) subdominant.  

    Since the contours include the actual local steepest descent paths and since an  inequality of the form $\mathrm{Re}(E(z))\le \mathrm{Re}(E(y))-\delta<0$ that held on the unchanged arcs of the contour not in small neighborhoods of any critical points is preserved under perturbation, the proof of Proposition~\ref{prop:W} is complete in the more general case that only $(t,x)\not\in\mathcal{D}$.
\end{proof}

\section{Steepest descent analysis}\label{sec:steepest-descent}
In this section, we complete the proof of Theorem~\ref{thm:u-app} in four broad strokes.  First, we apply the steepest descent method on the matrices $\mathbf{A}(t,x;\epsilon)$, $\mathbf{B}(t,x;\epsilon)$ by using the contours $\{W_n\}_{n=0}^N$ from Proposition~\ref{prop:W} and obtain the leading order behavior of $u(t,x;\epsilon)$ in terms of a ratio of determinants of $(N+1)\times(N+1)$ matrices $\mathbf{A}^{0}(t,x;\epsilon)$ and $\mathbf{B}^{0}(t,x;\epsilon)$, see Lemma~\ref{lem:u-asymp-fraction}.  We then establish a differential identity relating the determinants of $\mathbf{A}^{0}(t,x;\epsilon)$ and $\mathbf{B}^{0}(t,x;\epsilon)$ and use this relationship to prove Proposition~\ref{prop:ultimate}, which yields the leading order behavior of $u(t,x;\epsilon)$ only in terms of $\det\mathbf{B}^{0}(t,x;\epsilon)$.  Next, in Propositions \ref{prop:schur}, \ref{prop:super-ultimate}, we use the Schur complement formula to express $\det\mathbf{B}^0(t,x;\epsilon)$ in terms of a determinant of a $J\times J$ matrix $\widetilde{\mathbf{M}}(t,x;\epsilon)$.  Finally, in Lemmas \ref{lem:abstildeM}, \ref{lem:Phi}, we show that $\widetilde{\mathbf{M}}(t,x;\epsilon)=\mathbf{M}(t,x;\epsilon)$, which completes the proof of Theorem~\ref{thm:u-app}.

\subsection{Expanding the determinants}
To begin, we remind the reader of frequently used notation and mention two important identities.

\begin{remark}\label{rem:hprop}
We recall that $y_0(t,x)>\cdots>y_{2J(t,x)}(t,x)$ are the real critical points of $h(z;t,x)$ and $z_j(t,x)$, $z_j(t,x)^*$, $j=1,\ldots,N-J(t,x)$ are the non-real critical points of $h(z;t,x)$ with $\mathrm{Im}(z_j(t,x))>0$.  It follows from the definition~\eqref{eq:h-def} of $h(z;x,t)$ that
\begin{align}\label{eq:hpcrit}
    0=h'(y_k(t,x);t,x)=\frac{y_{k}(t,x)-x}{2t}+u_0(y_k(t,x)),
\end{align}
where $\prime$ denotes differentiation with respect to $z$.  An immediate consequence of this identity is
\begin{align}\label{eq:yjdiff}
    y_j(t,x)-y_k(t,x)=-2t(u_{j}^{\mathrm{B}}(t,x)-u_{k}^{\mathrm{B}}(t,x)),
\end{align}
where we recall that $u_0(y_k(t,x))=u_{k}^{\mathrm{B}}(t,x)$. Sometimes it will be convenient to denote $z_0(t,x)=y_0(t,x)$.  We will often omit the dependence on $(t,x)$ in $h(z;t,x)$, $y_k(t,x)$, $J(t,x)$, $u_{k}^{\mathrm{B}}(t,x)$ for convenience.
\end{remark}

The identities of Remark~\ref{rem:hprop} are used to derive the compact expressions for the local wavenumbers and frequencies given earlier in Remark~\ref{rem:local-wavenumber}. Indeed, we have
\begin{equation}
\begin{split}
    \kappa_j(t,x):=\frac{\partial\theta_j}{\partial x}(t,x)&=\frac{\partial h}{\partial y}(y_{2j-1}(t,x);t,x)\frac{\partial y_{2j-1}}{\partial x}(t,x)+\frac{\partial h}{\partial x}(y_{2j-1}(t,x);t,x)\\
    &\quad\quad-\frac{\partial h}{\partial y}(y_{2j}(t,x);t,x)\frac{\partial y_{2j}}{\partial x}(t,x)-\frac{\partial h}{\partial x}(y_{2j}(t,x);t,x)\\
    &=\frac{x-y_{2j-1}(t,x)}{2t}-\frac{x-y_{2j}(t,x)}{2t}\\
    &=u^\mathrm{B}_{2j-1}(t,x)-u^\mathrm{B}_{2j}(t,x),
    \end{split}
\label{eq:wavenumbers}
\end{equation}
\begin{equation}
\begin{split}
    \omega_j(t,x):=-\frac{\partial\theta_j}{\partial t}(t,x)&= -\frac{\partial h}{\partial y}(y_{2j-1}(t,x);t,x)\frac{\partial y_{2j-1}}{\partial t}(t,x)-\frac{\partial h}{\partial t}(y_{2j-1}(t,x);t,x)\\
    &\quad\quad+\frac{\partial h}{\partial y}(y_{2j}(t,x);t,x)\frac{\partial y_{2j}}{\partial t}(t,x)+\frac{\partial h}{\partial t}(y_{2j}(t,x);t,x)\\
    &=\frac{(y_{2j-1}(t,x)-x)^2}{4t^2}-\frac{(y_{2j}(t,x)-x)^2}{4t^2}\\
    &=u^\mathrm{B}_{2j-1}(t,x)^2-u^\mathrm{B}_{2j}(t,x)^2.
\end{split}
\label{eq:frequencies}
\end{equation}

We are now ready to apply the method of steepest descent to the matrices $\mathbf{A}^{0}(t,x;\epsilon)$, $\mathbf{B}^{0}(t,x;\epsilon)$ appearing in the formula  \eqref{eq:lambda-formula} for $u(t,x;\epsilon)$.

\begin{lemma}\label{lem:u-asymp-fraction}
Let $\epsilon>0$ be sufficiently small.  Then,
\begin{equation}\label{eq:u-asymp-fraction}
    u(t,x;\epsilon)=2\mathrm{Re}\left(\frac{\det(\mathbf{A}^0(t,x;\epsilon))}{\det(\mathbf{B}^0(t,x;\epsilon))}\right) + \mathcal{O}(\epsilon),
\end{equation}
where the error term is uniformly small and
\begin{multline}
\mathbf{A}^0(t,x;\epsilon):=\\
    \begin{bmatrix}
        u_0(y_1)+\alpha_1 u_0(y_2) & \displaystyle\frac{1}{y_1-p_1}+\frac{\alpha_1}{y_2-p_1} &\cdots & \displaystyle\frac{1}{y_1-p_N}+\frac{\alpha_1}{y_2-p_N}\\
        \vdots & \vdots & \ddots &\vdots\\
        u_0(y_{2J-1})+\alpha_J u_0(y_{2J}) & \displaystyle\frac{1}{y_{2J-1}-p_1}+\frac{\alpha_J}{y_{2J}-p_1}&\cdots &
        \displaystyle\frac{1}{y_{2J-1}-p_N}+\frac{\alpha_J}{y_{2J}-p_N}\\
        u_0(z_0) & \displaystyle \frac{1}{z_0-p_1} & \cdots & \displaystyle \frac{1}{z_0-p_N}\\
        \vdots & \vdots & \ddots & \vdots\\
        u_0(z_{N-J}) & \displaystyle\frac{1}{z_{N-J}-p_1} & \cdots &
        \displaystyle\frac{1}{z_{N-J}-p_N}
    \end{bmatrix},
\end{multline}
\begin{equation}
     \mathbf{B}^0(t,x;\epsilon):=
    \begin{bmatrix}
        1+\alpha_1& \displaystyle\frac{1}{y_1-p_1}+\frac{\alpha_1}{y_2-p_1} &\cdots & \displaystyle\frac{1}{y_1-p_N}+\frac{\alpha_1}{y_2-p_N}\\
        \vdots & \vdots & \ddots &\vdots\\
       1+\alpha_J  & \displaystyle\frac{1}{y_{2J-1}-p_1}+\frac{\alpha_J}{y_{2J}-p_1}&\cdots &
        \displaystyle\frac{1}{y_{2J-1}-p_N}+\frac{\alpha_J}{y_{2J}-p_N}\\
        1 & \displaystyle \frac{1}{z_0-p_1} & \cdots & \displaystyle \frac{1}{z_0-p_N}\\
        \vdots & \vdots & \ddots & \vdots\\
        1 & \displaystyle\frac{1}{z_{N-J}-p_1} & \cdots &
        \displaystyle\frac{1}{z_{N-J}-p_N}
    \end{bmatrix},
\end{equation}
and
\begin{equation}\label{eq:alpha-j}
    \alpha_j:=-\ii\sqrt{-\frac{h''(y_{2j-1})}{h''(y_{2j})}}\ee^{\ii\theta_j/\epsilon},\quad\theta_j:=h(y_{2j-1})-h(y_{2j})\in\mathbb{R}, \quad j=1,\dots J.
\end{equation}
\end{lemma}

We will show below that $\det(\mathbf{B}^0(t,x;\epsilon))$ is bounded away from zero uniformly for $(t,x)$ with $t>0$ in a neighborhood of any point for which $h(z)$ has simple critical points and for all $\epsilon>0$ (see identity~\eqref{eq:B-to-M} and Remark~\ref{rem:B-not-zero}).  From this it follows from Theorem~\ref{thm:inversion-formula} that~\eqref{eq:u-asymp-fraction} holds.

\begin{proof}
By factoring out a sign from the corresponding determinants, we can re-order the rows of $\mathbf{A}(t,x;\epsilon)$ and $\mathbf{B}(t,x;\epsilon)$ so that $W_j$ with $j=1,\dots,J$ (in row $j+1$) is the path that corresponds to the pair $(y_{2j},y_{2j-1})$ of consecutive real critical points.  Then by the same sequence of elementary row operations applied to both determinants, we cancel the contributions to integrals over $W_0$ from all real critical points except for $y_0$ and also ensure that the elements in row $j=2,\dots,J+1$ are corresponding to the paths that will produce contributions from only the pair $(y_{2j-2},y_{2j-3})$.  Finally, we move the first row (originally corresponding to the contour $W_0$) and insert it after the $J^\mathrm{th}$ row in each determinant. 

Then the application of Lemma~\ref{lem:steepest-descent} to the modified determinants shows that for some common nonzero factor $K(t,x;\epsilon)$,
\begin{equation}
    \frac{\det(\mathbf{A}(t,x;\epsilon))}{K(t,x;\epsilon)}=\det(\mathbf{A}^0(t,x;\epsilon)) + \mathcal{O}(\epsilon),
    \quad
    \frac{\det(\mathbf{B}(t,x;\epsilon))}{K(t,x;\epsilon)}=\det(\mathbf{B}^0(t,x;\epsilon))+\mathcal{O}(\epsilon).
\end{equation}
The factor $K(t,x;\epsilon)$ includes, among others, the products $\ee^{-\ii h(y_1)/\epsilon}\ee^{-\ii h(y_3)/\epsilon}\cdots\ee^{-\ii h(y_{2J-1})/\epsilon}$ and $\ee^{-\ii h(z_0)/\epsilon}\ee^{-\ii h(z_1)/\epsilon}\cdots\ee^{-\ii h(z_{N-J})/\epsilon}$ of $N+1$ exponential factors from~\eqref{eq:steepest-descent}, each of which is common to all of the elements in the same row. Here the only dependence on $\epsilon$ enters via the constants $\alpha_j$ in~\eqref{eq:alpha-j}, which comes from the ratio of formula~\eqref{eq:steepest-descent} at the points $y_{2j-1}$ and $y_{2j}$, with steepest descent angles $\phi=\mp\frac{\pi}{4}$, explaining the factor $-\ii$ in~\eqref{eq:alpha-j}.

Since by choice of branch cuts the function $h(z)$ is analytic along each integration contour, and since by construction the contour passing through the real critical points $y_{2j-1}$ and $y_{2j}$ itself remains close to the real line between them, we may use the formula~\eqref{eq:h-def} which is real-valued and analytic on $\mathbb{R}$ to define the values $h(y_{2j-1})$ and $h(y_{2j})$ in the definition of $\theta_j$.
\end{proof}

The following identity, which is a mild generalization of~\cite[eq.~(10)]{Gerard23}, will be useful in the proof of Propositions~\ref{prop:J0} and~\ref{prop:ultimate}.

\begin{lemma}\label{lem:global-identity}
The zeros and poles of $h'(z)$ satisfy the relationship
\begin{equation}\label{eq:global-identity}
    y_0 + \sum_{j=1}^J(y_{2j-1}+y_{2j}) + \sum_{m=1}^{N-J}(z_m+z_m^*) = x+\sum_{n=1}^N(p_n+p_n^*).
\end{equation}
\end{lemma}

\begin{proof}
Since $2th'(z)=z-x+2tu_0(z)$ is a rational function which has simple poles at $p_1,\dots,p_n$ and $p_1^*,\dots,p_n^*$ only, has $2N+1$ zeros $z_1,\dots,z_{N-J}, z_1^*,\dots,z_{N-J}^*$ and $y_0,\dots,y_{2J}$, and as $z\to\infty$, $2th'(z)=z+\mathcal{O}(1)$.  Therefore we can write $2th'(z)$ in factored form as
\begin{equation}\label{eq:hprimefactors}
    2th'(z)=\frac{\displaystyle\prod_{j=0}^{2J}(z-y_j)\cdot\prod_{m=1}^{N-J}(z-z_m)(z-z_m^*)}{\displaystyle\prod_{n=1}^N(z-p_n)(z-p_n^*)}.
\end{equation}
On the other hand, using just the definition of $u_0(z)$~\eqref{eq:rationalIC},
\begin{multline*}
    2th'(z)=\frac{(z-x)\displaystyle\prod_{n=1}^N(z-p_n)(z-p_n^*)}{\displaystyle\prod_{n=1}^N(z-p_n)(z-p_n^*)}+ \\
    2t\frac{\displaystyle\sum_{k=1}^Nc_k\displaystyle\prod_{\substack{n=1 \\ n\neq k}}^{N}(z-p_n)\cdot\displaystyle\prod_{m=1}^N(z-p_m^*) + \displaystyle\sum_{k=1}^Nc_k^*\displaystyle\prod_{n=1}^N(z-p_n)\cdot\displaystyle\prod_{\substack{m=1 \\ m\neq k}}^N(z-p_m^*)}{\displaystyle\prod_{n=1}^N(z-p_n)(z-p_n^*)}.
\end{multline*}
The identity~\eqref{eq:global-identity} follows from comparing coefficients of the $z^{2N}$ coefficients in the numerators of the previous two equations.
\end{proof}

With the help of Lemmas~\ref{lem:u-asymp-fraction} and~\ref{lem:global-identity}, we can now prove Theorem~\ref{thm:u-app} for $J(t,x)=0$.

\begin{proposition}\label{prop:J0}
Suppose $J(t,x)=0$ and $\epsilon>0$ is sufficiently small. Then,
\begin{equation}
u(t,x;\epsilon)=u_0^B(t,x)+\mathcal{O}(\epsilon),
\end{equation}
where the error term is uniform.
\end{proposition}

\begin{proof}
In the case $J=0$, the quantities $\alpha_j$ do not appear in the matrices $\mathbf{A}^0$ or $\mathbf{B}^0$, so they are independent of $\epsilon$ and here we will write $\mathbf{A}^0=\mathbf{A}^0(t,x)$ and $\mathbf{B}^0=\mathbf{B}^0(t,x)$.  These matrices have the simpler structure
\begin{equation}
    \mathbf{A}^0=\begin{bmatrix}u_0(z_0) &\displaystyle\frac{1}{z_0-p_1} &\cdots & \displaystyle\frac{1}{z_0-p_N}\\
    \vdots & \vdots & \ddots & \vdots\\
    u_0(z_N) & \displaystyle\frac{1}{z_N-p_1} &\cdots & \displaystyle\frac{1}{z_N-p_N}\end{bmatrix},\;
\mathbf{B}^0=\begin{bmatrix}
1 &\displaystyle\frac{1}{z_0-p_1} &\cdots & \displaystyle\frac{1}{z_0-p_N}\\
    \vdots & \vdots & \ddots & \vdots\\
  1 & \displaystyle\frac{1}{z_N-p_1} &\cdots & \displaystyle\frac{1}{z_N-p_N}    
\end{bmatrix}.
\end{equation}
It is convenient to replace the first column of $\mathbf{A}^0$ by using the critical point equation~\eqref{eq:hpcrit}. 
 Then,~\cite[Lemma 8]{Gerard23} yields
\begin{equation}
    \frac{\det(\mathbf{A}^0(t,x))}{\det(\mathbf{B}^0(t,x))}=\frac{1}{2t}\left(x-\sum_{m=0}^{N}z_m +\sum_{n=1}^Np_n\right).
\end{equation}
Since $z_0=y_0$ is the only real critical point, from~\eqref{eq:u-asymp-fraction} we obtain
\begin{align*}
    u(t,x;\epsilon)&=\frac{1}{2t}\left(2x-2y_0-\sum_{m=1}^N(z_m+z_m^*) +\sum_{n=1}^N(p_n+p_n^*)\right) + \mathcal{O}(\epsilon) \\
    &=\frac{x-y_0}{2t}+\mathcal{O}(\epsilon)=u_0^\mathrm{B}(t,x)+\mathcal{O}(\epsilon),
\end{align*}
where we have used Lemma~\ref{lem:global-identity} in the $J=0$ setting.
\end{proof}

\subsection{Differential identities and expressing \texorpdfstring{$u(t,x;\epsilon)$}{u} in terms of \texorpdfstring{$\det(\mathbf{B}^0)$}{detB0} only}
Now we consider $J(t,x)\geq1$.  First, the matrices $\mathbf{A}^{0}$, $\mathbf{B}^{0}$ only differ in their first column, so it should not be a surprise that their determinants are related.

\begin{lemma}\label{lem:detAtodetB}
Suppose $\epsilon>0$ is sufficiently small.  Then, 
\begin{equation}
    \det(\mathbf{A}^0)=\ii\epsilon\frac{\partial}{\partial x}\det(\mathbf{B}^0)+\frac{1}{2t}\left(x+\sum_{n=1}^Np_n-\sum_{j=1}^Jy_{2j-1}-\sum_{m=0}^{N-J}z_m\right)\det(\mathbf{B}^0)+\mathcal{O}(\epsilon),
\end{equation}
where the error term is uniform.
\end{lemma}

\begin{proof}
The Lemma will follow from the differentiation of $\det\mathbf{B}^{0}$, so to prepare for the computation we first prove the identities
\begin{align}\label{eq:derivs-columns}
\begin{split}
    \epsilon\frac{\partial}{\partial x}\vec{B}_{1}^{0}&=\ii\widetilde{\mathbf{D}}\vec{B}_{1}^{0}-\ii\vec{A}_{1}^{0}+\mathcal{O}(\epsilon), \\
    \epsilon\frac{\partial}{\partial x}\vec{B}_{k}&=\ii\widetilde{\mathbf{D}}\vec{B}_{k}^{0}-\ii\frac{x-p_{k-1}}{2t}\vec{A}_{k}^{0}+\frac{\ii}{2t}\vec{B}_{1}^{0}+\mathcal{O}(\epsilon), \quad k=2,\ldots,N+1,
\end{split}
\end{align}
where the notation $\vec{C}_k$ refers to the $k^\mathrm{th}$ column of the matrix $\mathbf{C}$ and
\begin{align}
    \widetilde{\mathbf{D}}:=\mathrm{diag}(u_0(y_1),u_0(y_3),\ldots,u_0(y_{2J-1}),u_0(z_0),\ldots,u_0(z_{N-J})).
\end{align}
Note that error terms of order $\epsilon$ will be generated by the action of the differential operator $\epsilon\partial_x$ unless it acts on one of the exponentials $\ee^{\ii\theta_j(t,x)/\epsilon}$, in which case~\eqref{eq:wavenumbers} applies.  Consequently, if $j=1,\dots,J$, we have from~\eqref{eq:alpha-j} that
\begin{equation}
\begin{split}
    \epsilon\frac{\partial B_{j1}^0}{\partial x}&= -\ii \alpha_j\frac{\partial}{\partial x}\theta_j(t,x) + \mathcal{O}(\epsilon) \\
    &=-\ii (u_0(y_{2j-1})+u_0(y_{2j})\alpha_j)+\ii u_0(y_{2j-1})(1+\alpha_j) + \mathcal{O}(\epsilon) \\
    &=-\ii A_{j1}^0 + \ii u_0(y_{2j-1})B_{j1}^0 + \mathcal{O}(\epsilon).
\end{split}
\end{equation}
An even simpler exact calculation for $j=J+1,\dots,N+1$ gives
\begin{equation}
    \epsilon\frac{\partial B_{j1}^0}{\partial x}=\epsilon\frac{\partial}{\partial x}(1)=0=-\ii u_0(z_{j-J-1}) + \ii u_0(z_{j-J-1})\cdot 1 = -\ii A_{j1}^0 +\ii u_0(z_{j-J-1})B_{j1}^0,
\end{equation}
and therefore the first equality of~\eqref{eq:derivs-columns} is proven.  Now let $k=2,3,\dots,N+1$.  If $j=1,\dots,J$, then
\begin{equation}
\begin{split}
    \epsilon\frac{\partial B_{jk}^0}{\partial x}&=\frac{\ii \alpha_j}{y_{2j}-p_{k-1}}\frac{\partial}{\partial x}\theta_j(t,x) + \mathcal{O}(\epsilon)
    =\ii u_0(y_{2j-1})\frac{\alpha_j}{y_{2j}-p_{k-1}}-\ii u_0(y_{2j})\frac{\alpha_j}{y_{2j}-p_{k-1}} + \mathcal{O}(\epsilon).
\end{split}
\end{equation}
But because $x-y=2tu_0(y)$ holds for $y=y_{2j-1}$ and $y=y_{2j}$,
\begin{equation}
\begin{split}
    \epsilon\frac{\partial B_{jk}^0}{\partial x}
    &=\ii u_0(y_{2j-1})B_{jk}^0-\frac{\ii u_0(y_{2j-1})}{y_{2j-1}-p_{k-1}}-\frac{\ii\alpha_j u_0(y_{2j})}{y_{2j}-p_{k-1}} + \mathcal{O}(\epsilon)\\
    &=\ii u_0(y_{2j-1})B_{jk}^0 -\frac{\ii}{2t}\left(\frac{x-y_{2j-1}}{y_{2j-1}-p_{k-1}}+\alpha_j\frac{x-y_{2j}}{y_{2j}-p_{k-1}}\right) + \mathcal{O}(\epsilon)\\
    &=\ii u_0(y_{2j-1})B_{jk}^0 -\ii\frac{x-p_{k-1}}{2t}A_{jk}^0+\frac{\ii}{2t}B_{j1}^0 + \mathcal{O}(\epsilon).
\end{split}
\end{equation}
Likewise, if $j=J+1,\dots,N+1$, then 
\begin{equation}
    \epsilon\frac{\partial B_{jk}^0}{\partial x} 	
    =\epsilon\frac{\partial}{\partial x}\left(\frac{1}{z_{j-J-1}-p_{k-1}}\right)
    =\mathcal{O}(\epsilon),
\end{equation}
however it is also true for these indices that
\begin{equation}
\begin{split}
    -\ii\frac{x-p_{k-1}}{2t}A_{jk}^0 +\frac{\ii}{2t}B_{j1}^0 +\ii u_0(z_{j-J-1})B_{jk}^0
    &=\frac{\ii u_0(z_{j-J-1})}{z_{j-J-1}-p_{k-1}}-\frac{\ii}{2t}\frac{x-p_{k-1}}{z_{j-J-1}-p_{k-1}}+\frac{\ii}{2t}\\
    &=\frac{\ii}{2t}\left(\frac{x-z_{j-J-1}}{z_{j-J-1}-p_{k-1}}-\frac{x-p_{k-1}}{z_{j-J-1}-p_{k-1}} + 1\right)\\
    &=0,
\end{split}
\end{equation}
which completes the verification of~\eqref{eq:derivs-columns}.  Now, define the matrix $\mathbf{D}$ as
\begin{align}
    \mathbf{D}:=\mathrm{diag}\left(\ee^{-\ii h(y_1)/\epsilon},\ee^{-\ii h(y_3)/\epsilon},\ldots,\ee^{-\ii h(y_{2J-1})/\epsilon},\ee^{-\ii h(z_0)/\epsilon},\ldots,\ee^{-\ii h(z_{N-J})/\epsilon}\right).
\end{align}
Via chain rule and Remark~\ref{rem:hprop},
\begin{align}\label{eq:D-deriv}
    \epsilon\frac{\partial}{\partial x}\mathbf{D}=-\ii\widetilde{\mathbf{D}}\mathbf{D}.
\end{align}
Calculating the derivative by summing determinants with differentiated columns using the fact that $\mathbf{A}^0$ and $\mathbf{B}^0$ only differ in the first column, one applies~\eqref{eq:derivs-columns} and~\eqref{eq:D-deriv} to obtain
\begin{equation}\label{eq:partialxDB1}
    \epsilon\frac{\partial}{\partial x}\det(\mathbf{D}\mathbf{B}^0) = \det(\mathbf{D})\left[-\frac{\ii}{2t}\left(Nx-\sum_{n=1}^Np_n\right)\det(\mathbf{B}^0)-\ii\det(\mathbf{A}^0)+\mathcal{O}(\epsilon)\right].
\end{equation}
On the other hand, using~\eqref{eq:D-deriv} without differentiating explicitly $\mathbf{B}^0$, we have the exact result
\begin{equation}\label{eq:partialxDB2}
\begin{split}
    \epsilon\frac{\partial}{\partial x}\det(\mathbf{D}\mathbf{B}^0)&=\det(\mathbf{D})\left[-\ii\left(\sum_{j=1}^Ju_0(y_{2j-1})+\sum_{m=0}^{N-J}u_0(z_m)\right)\det(\mathbf{B}^0) + \epsilon\frac{\partial}{\partial x}\det(\mathbf{B}^0)\right]\\
    &=\det(\mathbf{D})\left[-\ii\left(\sum_{j=1}^J\frac{x-y_{2j-1}}{2t}+\sum_{m=0}^{N-J}\frac{x-z_m}{2t}\right)\det(\mathbf{B}^0)+\epsilon\frac{\partial}{\partial x}\det(\mathbf{B}^0)\right]\\
    &=\det(\mathbf{D})\left[-\frac{\ii}{2t}\left((N+1)x-\sum_{j=1}^Jy_{2j-1}-\sum_{m=0}^{N-J}z_m\right)\det(\mathbf{B}^0)+\epsilon\frac{\partial}{\partial x}\det(\mathbf{B}^0)\right].
\end{split}
\end{equation}
Since $\det(\mathbf{D})\neq 0$, the lemma follows by comparing~\eqref{eq:partialxDB1} and~\eqref{eq:partialxDB2}.
\end{proof}

Using the identity of Lemma~\ref{lem:detAtodetB} in Lemma~\ref{lem:u-asymp-fraction}, and for any branch of the complex logarithm,
\begin{multline}
    2\mathrm{Re}\left(\frac{\det(\mathbf{A}^0)}{\det(\mathbf{B}^0)}\right)=
    \frac{1}{2t}\left(2x-2y_0-2\sum_{j=1}^Jy_{2j-1} +\sum_{n=1}^N(p_n+p_n^*) -\sum_{m=1}^{N-J}(z_m+z_m^*)\right)\\
    -2\epsilon\mathrm{Im}\left(\frac{\partial}{\partial x}\log(\det(\mathbf{B}^0))\right) + \mathcal{O}(\epsilon).
\label{eq:penultimate}
\end{multline}
Applying the identity~\eqref{eq:global-identity} to eliminate the sum over poles in~\eqref{eq:penultimate} yields
\begin{equation}
\begin{split}
    2\mathrm{Re}\left(\frac{\det(\mathbf{A}^0)}{\det(\mathbf{B}^0)}\right)&=\frac{1}{2t}\left(x-y_0+\sum_{j=1}^J(y_{2j}-y_{2j-1})\right) -2\epsilon\mathrm{Im}\left(\frac{\partial}{\partial x}\log(\det(\mathbf{B}^0))\right)+\mathcal{O}(\epsilon)\\
    &=u_0(y_0)+\sum_{j=1}^J(u_0(y_{2j-1})-u_0(y_{2j}))
    -2\epsilon\mathrm{Im}\left(\frac{\partial}{\partial x}\log(\det(\mathbf{B}^0))\right)+\mathcal{O}(\epsilon).
\end{split}
\end{equation}
Thus, we have proven the following Proposition.

\begin{proposition}\label{prop:ultimate}
Suppose $J(t,x)\geq1$ and $\epsilon>0$ is sufficiently small.  Then,
\begin{align}\label{eq:ultimate}
    u(t,x;\epsilon)=u_0^\mathrm{B}(t,x)+\sum_{j=1}^J(u^\mathrm{B}_{2j-1}(t,x)-u^\mathrm{B}_{2j}(t,x))-2\epsilon\mathrm{Im}\left(\frac{\partial}{\partial x}\log \det(\mathbf{B}^0(t,x;\epsilon))\right)+\mathcal{O}(\epsilon),
\end{align}
where the error term is uniform.
\end{proposition}

\begin{remark}
According to \cite[Theorem 1.3]{RationalBO}, an equivalent expression of $u(t,x;\epsilon)$ given in \eqref{eq:lambda-formula} is
\begin{align}\label{eq:u-tildeB}
    u(t,x;\epsilon)=-2\epsilon\mathrm{Im}\left(\frac{\partial}{\partial x}\log\det\widetilde{\mathbf{B}}(t,x;\epsilon)\right),
\end{align}
where the first column of $\mathbf{B}(t,x;\epsilon)$ and $\widetilde{\mathbf{B}}(t,x;\epsilon)$ are equal and $(\widetilde{\mathbf{B}})_{jk}=\ee^{\ii(x-p_{k-1})^2/(4t\epsilon)}B_{jk}$ for $k=2,\dots,N+1$, where we recall that $\mathbf{B}(t,x;\epsilon)$ was defined in \eqref{eq:Bmatrix}. One can apply the steepest descent method to $\widetilde{\mathbf{B}}(t,x;\epsilon)$ and formally differentiate the result to more directly obtain the leading terms in Proposition~\ref{prop:ultimate}, but this approach is not rigorous because the topology of the estimate of the error term is not sufficient to control its derivative.  In general, great care is required to properly mix differentiation and asymptotics.
\end{remark}

\subsection{Schur complement reduction to a \texorpdfstring{$J\times J$}{JxJ} determinant}
Next, we strive to express \eqref{eq:ultimate} in a form resembling the multiphase solution formula of Dobrokhotov and Krichever \eqref{eq:DK-exact} for $J\geq1$.  We can write the determinant of the $(N+1)\times(N+1)$ matrix $\mathbf{B}^{0}(t,x;\epsilon)$ in terms of a $J\times J$ determinant via the Schur complement formula.  Explicitly,
\begin{equation}\label{eq:detBSchur}
    \frac{\det(\mathbf{B}^0(t,x;\epsilon))}{\det(\mathbf{C})}=\det(\mathbf{K}-\mathbf{H}\mathbf{C}^{-1}\mathbf{J}),
\end{equation}
where $\mathbf{K}$, $\mathbf{H}$, $\mathbf{J}$, $\mathbf{C}$ are matrices of size $J\times J$, $J\times(N-J+1)$, $(N-J+1)\times J$, $(N-J+1)\times(N-J+1)$, respectively, and are defined as
\begin{equation*}
 \mathbf{K}   :=\begin{bmatrix}
        1+\alpha_1 & \displaystyle\frac{1}{y_1-p_1}+\frac{\alpha_1}{y_2-p_1} & \cdots & \displaystyle\frac{1}{y_1-p_{J-1}}+\frac{\alpha_1}{y_2-p_{J-1}}\\
        \vdots & \vdots &\ddots & \vdots\\
        1 + \alpha_J & \displaystyle\frac{1}{y_{2J-1}-p_1}+\frac{\alpha_J}{y_{2J}-p_1} &\cdots &\displaystyle\frac{1}{y_{2J-1}-p_{J-1}}+\frac{\alpha_J}{y_{2J}-p_{J-1}}
    \end{bmatrix},
\end{equation*}
\begin{equation*}
    \mathbf{H}:=\begin{bmatrix}
        \displaystyle\frac{1}{y_1-p_J}+\frac{\alpha_1}{y_2-p_J} & \cdots & \displaystyle\frac{1}{y_1-p_N} +\frac{\alpha_1}{y_2-p_N}\\
        \vdots & \ddots & \vdots\\
        \displaystyle\frac{1}{y_{2J-1}-p_J}+\frac{\alpha_J}{y_{2J}-p_J} & \cdots & \displaystyle\frac{1}{y_{2J-1}-p_N} + \frac{\alpha_J}{y_{2J}-p_N}
    \end{bmatrix},
\end{equation*}
\begin{equation*}
 \mathbf{J}   
    :=\begin{bmatrix}
        1 & \displaystyle\frac{1}{z_0-p_1} & \cdots & \displaystyle\frac{1}{z_0-p_{J-1}}\\
        \vdots & \vdots & \ddots & \vdots\\
        1 & \displaystyle\frac{1}{z_{N-J}-p_1} & \cdots & \displaystyle\frac{1}{z_{N-J}-p_{J-1}}
    \end{bmatrix}, \quad
    \mathbf{C}:=\begin{bmatrix}
    \displaystyle\frac{1}{z_0-p_J} & \cdots & \displaystyle\frac{1}{z_0-p_N}\\
    \vdots & \ddots & \vdots\\
    \displaystyle\frac{1}{z_{N-J}-p_J}&\cdots &\displaystyle\frac{1}{z_{N-J}-p_N}\end{bmatrix}.
\end{equation*}
It is well-known that the Cauchy matrix $\mathbf{C}$ is explicitly invertible.

\begin{proposition}\label{prop:schur}
We have the identity
\begin{multline}\label{eq:detB-2}
    \frac{\det(\mathbf{B}^0(t,x;\epsilon))}{\det(\mathbf{C})} \\
    =\begin{vmatrix}
        W_1+W_2\alpha_1 & \displaystyle\frac{Q_2W_1}{y_1-p_1} +\frac{Q_2W_2\alpha_1}{y_2-p_1} &\cdots&
        \displaystyle\frac{Q_JW_1}{y_1-p_{J-1}}+\frac{Q_JW_2\alpha_1}{y_2-p_{J-1}}\\
        \vdots & \vdots & \ddots & \vdots\\
        W_{2J-1}+W_{2J}\alpha_J & 
        \displaystyle\frac{Q_2W_{2J-1}}{y_{2J-1}-p_{1}}+\frac{Q_2W_{2J}\alpha_J}{y_{2J}-p_1} & \cdots & 
        \displaystyle\frac{Q_JW_{2J-1}}{y_{2J-1}-p_{J-1}}+\frac{Q_JW_{2J}\alpha_J}{y_{2J}-p_{J-1}}
    \end{vmatrix},
\end{multline}
where, for $j=1,\dots,2J$ and $l=1,\dots, J$,
\begin{equation}\label{eq:const-Ql}
    W_j:=\frac{\displaystyle\prod_{m=0}^{N-J}(y_j-z_m)}{\displaystyle\prod_{n=J}^N(y_j-p_n)},
\quad
    Q_l:=\frac{\displaystyle \prod_{n=J}^{N}(p_{l-1}-p_n)}{\displaystyle \prod_{m=0}^{N-J}(p_{l-1}-z_m)}.
\end{equation}
\end{proposition}

\begin{proof}
From simple matrix multiplication,
\begin{align}\label{eq:B0schur expand}
    (\mathbf{K}-\mathbf{H}\mathbf{C}^{-1}\mathbf{J})_{il}=g_l(y_{2i-1})+\alpha_ig_l(y_{2i}), \quad i,j=1,\ldots,J,
\end{align}
where
\begin{align}\label{def:gl}
    g_{l}(y):=\begin{cases}
        \displaystyle{1-\sum_{j,k=1}^{N-J+1}\frac{(\mathbf{C}^{-1})_{jk}}{y-p_{J+j-1}}}, &l=1, \\
        \displaystyle{\frac{1}{y-p_{l-1}}-\sum_{j,k=1}^{N-J+1}\frac{(\mathbf{C}^{-1})_{jk}}{(z_{k-1}-p_{l-1})(y-p_{J+j-1})}},& l=2,\ldots,J.
    \end{cases}
\end{align}
We will show that $g_l(y)$ can be simplified to
\begin{align}\label{eq:glsimplify}
    g_{l}(y)=\begin{cases}
    \displaystyle{\frac{\displaystyle\prod_{m=0}^{N-J}(y-z_m)}{\displaystyle\prod_{n=J}^N(y-p_n)}}, &l=1, \\
    \displaystyle{\frac{Q_l}{y-p_{l-1}}\frac{\displaystyle\prod_{m=0}^{N-J}(y-z_m)}{\displaystyle\prod_{n=J}^{N}(y-p_n)}}, &l=2,\ldots,J.
    \end{cases}
\end{align}
The Proposition follows from this identity and \eqref{eq:detBSchur}, \eqref{eq:const-Ql}, \eqref{eq:B0schur expand}.  The derivation of \eqref{eq:glsimplify} is similar for $l=1$ and $l=2,\ldots,J$, so for brevity, we omit the proof for $l=1$.  For $1\le j,k\le N-J+1$ the elements of the inverse matrix $\mathbf{C}$ are
\begin{equation}
    (\mathbf{C}^{-1})_{jk}=\frac{1}{p_{J+j-1}-z_{k-1}}\cdot\frac{\displaystyle\prod_{m=1}^{N-J+1}(p_{J+j-1}-z_{m-1})\prod_{n=1}^{N-J+1}(z_{k-1}-p_{J+n-1})}{\displaystyle\mathop{\prod_{m=1}^{N-J+1}}_{m\neq k}(z_{k-1}-z_{m-1})\mathop{\prod_{n=1}^{N-J+1}}_{n\neq j}(p_{J+j-1}-p_{J+n-1})}.
\end{equation}
We can explicitly evaluate the double sum as follows.  By reindexing the sums, we have
\begin{multline}
    \sum_{j,k=1}^{N-J+1}
\frac{(\mathbf{C}^{-1})_{jk}}{(z_{k-1}-p_{l-1})(y-p_{J+j-1})}\\
	=\sum_{j=J}^N\sum_{k=0}^{N-J}
	\frac{1}{(z_k-p_{l-1})(y-p_{j})(p_j-z_k)}
	\frac{\displaystyle\prod_{m=0}^{N-J}(p_j-z_m)\cdot\prod_{n=J}^N(z_k-p_n)}{\displaystyle\mathop{\prod_{m=0}^{N-J}}_{m\neq k}(z_k-z_m)\cdot\mathop{\prod_{n=J}^N}_{n\neq j}(p_j-p_n)}\\
	=-\sum_{j=J}^N\frac{1}{y-p_j}\frac{\displaystyle\prod_{m=0}^{N-J}(p_j-z_m)}
	{\displaystyle\mathop{\prod_{n=J}^N}_{n\neq j}(p_j-p_n)}\sum_{k=0}^{N-J}\frac{1}{z_k-p_{l-1}}\frac{\displaystyle\prod_{\substack{ n=J\\ n\neq j}}^N(z_k-p_n)}
	{\displaystyle\mathop{\prod_{m=0}^{N-J}}_{m\neq k}(z_k-z_m)}.
	\label{eq:double-sum}
\end{multline}
The inner sum can be written as the sum of residues of a rational function. Indeed, by temporarily denoting
\begin{equation}
    f_j(z):=\frac{1}{z-p_{l-1}}\prod_{ \substack{n=J \\ n\neq j}}^{N}(z-p_n)\cdot{\displaystyle\prod_{m=0}^{N-J}(z-z_m)^{-1}},
\end{equation} 
it follows by the residue theorem that
\begin{equation}
    \sum_{k=0}^{N-J}\frac{1}{z_k-p_{l-1}}\displaystyle\prod_{\substack{ n=J \\ n\neq j}}^N(z_k-p_n)\cdot\displaystyle\mathop{\prod_{m=0}^{N-J}}_{m\neq k}(z_k-z_m)^{-1}
	=\sum_{k=0}^{N-J}\mathop{\mathrm{Res}}_{z=z_k}f_j
	=\frac{1}{2\pi\ii}\oint_L f_j(z)\,\dd z,
\end{equation}
where $L$ is a positively-oriented loop that encloses the points $z_0,\dots,z_{N-J}$ but excludes the other pole of the integrand at $z=p_{l-1}$ (since $2\le l\le J$ this is the only singularity of the integrand outside of $L$).  Noting that the integrand is $\mathcal{O}(z^{-2})$ as $z\to\infty$, the sum of residues above is equivalent to the single residue at $z=p_{l-1}$, so we have the identity
\begin{multline}
    \sum_{k=0}^{N-J}\frac{1}{z_k-p_{l-1}}\displaystyle\prod_{\substack{ n=J \\ n\neq j}}^N(z_k-p_n)\cdot\displaystyle\mathop{\prod_{m=0}^{N-J}}_{m\neq k}(z_k-z_m)^{-1} \\
	=-\mathop{\mathrm{Res}}_{z=p_{l-1}} f_j
	=-\displaystyle\prod_{\substack{ n=J\\ n \neq j}}^{N}(p_{l-1}-p_n)\cdot\displaystyle\prod_{m=0}^{N-J}(p_{l-1}-z_m)^{-1}.   
\end{multline}
Using this identity in~\eqref{eq:double-sum}, we get
\begin{equation}
\sum_{j,k=1}^{N-J+1}
\frac{(\mathbf{C}^{-1})_{jk}}{(z_{k-1}-p_{l-1})(y-p_{J+j-1})}=   \frac{\displaystyle\prod_{n=J}^{N}(p_{l-1}-p_n)}{\displaystyle\prod_{m=0}^{N-J}(p_{l-1}-z_m)}\sum_{j=J}^N\frac{1}{(y-p_j)(p_{l-1}-p_j)} 
\frac{\displaystyle\prod_{m=0}^{N-J}(p_j-z_m)}{\displaystyle\mathop{\prod_{n=J}^N}_{n\neq j}(p_j-p_n)}.
\end{equation}
The remaining sum can be evaluated in a similar manner.  By temporarily denoting
\begin{equation}
f(z)
	:=\frac{1}{(y-z)(p_{l-1}-z)}\displaystyle\prod_{m=0}^{N-J}(z-z_m)\cdot\displaystyle\prod_{n=J}^N(z-p_n)^{-1},
\end{equation}
we see that
\begin{equation}
    \sum_{j=J}^N\frac{1}{(y-p_j)(p_{l-1}-p_j)}\displaystyle\prod_{m=0}^{N-J}(p_j-z_m)\cdot\displaystyle\mathop{\prod_{n=J}^N}_{n\neq j}(p_j-p_n)^{-1}
	=\sum_{j=J}^N\mathop{\mathrm{Res}}_{z=p_j}f
	=\frac{1}{2\pi\ii}\oint_L f(z)\,\dd z,
\end{equation}
where $L$ is a contour with positive orientation that encloses the points $p_J,\dots,p_{N}$ but neither $y$ nor $p_{l-1}$ (the index $l-1$ is less than $J$).  As before, $f(z)=\mathcal{O}(z^2)$ as $z\to\infty$, so the sum of residues above is expressed in terms of residues at $z=y,p_{l-1}$ and we obtain the identity
\begin{multline}
\sum_{j=J}^N\frac{1}{(y-p_j)(p_{l-1}-p_j)} 
\frac{\displaystyle\prod_{m=0}^{N-J}(p_j-z_m)}{\displaystyle\mathop{\prod_{n=J}^N}_{n\neq j}(p_j-p_n)}=\frac{-1}{y-p_{l-1}}\frac{\displaystyle\prod_{m=0}^{N-J}(y-z_m)}{\displaystyle\prod_{n=J}^{N}(y-p_n)}+\frac{1}{y-p_{l-1}}\frac{\displaystyle\prod_{m=0}^{N-J}(p_{l-1}-z_m)}{\displaystyle\prod_{n=J}^{N}(p_{l-1}-p_n)}.    
\end{multline}
We have shown that
\begin{equation}
\sum_{j,k=1}^{N-J+1}
\frac{(\mathbf{C}^{-1})_{jk}}{(z_{k-1}-p_{l-1})(y-p_{J+j-1})}
=\frac{1}{y-p_{l-1}}-\frac{1}{y-p_{l-1}}\frac{\displaystyle\prod_{m=0}^{N-J}(y-z_m)\cdot\prod_{n=J}^{N}(p_{l-1}-p_n)}{\displaystyle\prod_{m=0}^{N-J}(p_{l-1}-z_m)\cdot\prod_{n=J}^{N}(y-p_n)},
\end{equation}
and it is now clear that \eqref{def:gl} and \eqref{eq:glsimplify} are equivalent, completing the proof.
\end{proof}

Now we use the previous proposition to reduce the size of the determinant appearing in \eqref{eq:ultimate}, while taking advantage of the operator $\partial_x\log$.

\begin{proposition}\label{prop:super-ultimate}
Suppose $J(t,x)\geq1$ and $\epsilon>0$ is sufficiently small.  Then,
\begin{align}\label{eq:super-ultimate}
    \epsilon\mathrm{Im}\left(\frac{\partial}{\partial x}\log \det(\mathbf{B}^0(t,x;\epsilon))\right)=\epsilon\mathrm{Im}\left(\frac{\partial}{\partial x}\log\det(\widetilde{\mathbf{M}}(t,x;\epsilon))\right)+\mathcal{O}(\epsilon),
\end{align}
where the error term is uniform and the matrix $\widetilde{\mathbf{M}}(t,x;\epsilon)$ is defined as
\begin{align}\label{eq:tildeM}
    \left(\widetilde{\mathbf{M}}(t,x;\epsilon)\right)_{jk}:=\widetilde{m}_j(t,x;\epsilon)\delta_{jk}+\frac{1}{u^\mathrm{B}_{2j-1}(t,x)-u^\mathrm{B}_{2k}(t,x)},
\end{align}
for $j,k=1,\ldots,J$, with
\begin{align}\label{eq:DK-cj}
    \widetilde{m}_j(t,x;\epsilon)=\frac{\alpha_jV_{2j}}{V_{2j-1}}\frac{\displaystyle\mathop{\prod_{i=1}^J}_{i\neq j}(u_{2j}^\mathrm{B}(t,x)-u_{2i}^\mathrm{B}(t,x))}{\displaystyle\prod_{i=1}^J(u_{2j-1}^\mathrm{B}(t,x)-u_{2i}^\mathrm{B}(t,x))},
\end{align}
and
\begin{equation}
 V_j:=\frac{\displaystyle\prod_{m=0}^{N-J}(y_j-z_m)}{\displaystyle\prod_{n=1}^N(y_j-p_n)},
 \quad j=1,\dots, 2J.
\end{equation}
\end{proposition}

\begin{proof}
Recalling the definition of $W_j$ in \eqref{eq:const-Ql}, we can write
\begin{equation}
    W_j
    = V_jP(y_j),
	\quad\text{where}\quad
     P(y):=\prod_{n=1}^{J-1}(y-p_n).
\label{eq:VjP-define}
\end{equation}
Unlike $W_j$, the factor $V_j$ is symmetrical in the poles $p_1,\dots,p_N$.  Let us denote by $P_n(y)$ the polynomials of degree $J-2$ given by
\begin{equation}
    P_n(y):=\mathop{\prod_{m=1}^{J-1}}_{m\neq n} (y-p_m),\quad n=1,\dots,J-1.
\end{equation}
Since these polynomials are proportional by nonzero factors to the Lagrange interpolating polynomials for the distinct points $\{p_1,\dots,p_{J-1}\}$, they are linearly independent, and hence 
\begin{equation}
    \mathrm{span}\left\{P_1(y),\dots,P_{J-1}(y)\right\}=\mathrm{span}\left\{1,y,\dots,y^{J-3},y^{J-2}\right\}.
\end{equation}
Therefore, there are coefficients forming a square matrix $\mathbf{S}$ of dimension $(J-1)\times (J-1)$ with $\det(\mathbf{S})\neq 0$ such that 
\begin{equation}
    P_n(y)=\sum_{j=1}^{J-1}s_{nj}y^{j-1},\quad n=1,\dots,J-1.
\end{equation}
Denoting also
\begin{equation}
    Q:=Q_2Q_3\cdots Q_J\neq 0,
\end{equation}
it follows by column operations applied to the last $J-1$ columns of the determinant on the right-hand side of
\eqref{eq:detB-2} that
\begin{multline}\label{eq:detB-2.5}
    \frac{\det(\mathbf{B}^0(t,x;\epsilon))}{\det(\mathbf{C})Q\det(\mathbf{S})}  \\
    =\begin{vmatrix}V_1 P(y_1) +V_2 P(y_2)\alpha_1 & V_1 + V_2\alpha_1&\cdots & V_1y_1^{J-2}+V_2y_2^{J-2}\alpha_1\\
    \vdots & \vdots & \ddots & \vdots\\
    V_{2J-1}P(y_{2J-1})+V_{2J}P(y_{2J})\alpha_J & V_{2J-1} + V_{2J}\alpha_J & \cdots & 
    V_{2J-1}y_{2J-1}^{J-2} + V_{2J}y_{2J}^{J-2}\alpha_J
    \end{vmatrix}
    \\
    =(-1)^J\begin{vmatrix}
        V_1+V_2\alpha_1 & \cdots & V_1y_1^{J-2}+V_2y_2^{J-2}\alpha_1 & V_1P(y_1)+V_2P(y_2)\alpha_1\\
        \vdots & \ddots & \vdots & \vdots\\
        V_{2J-1}+V_{2J}\alpha_J & \cdots & V_{2J-1}y_{2J-1}^{J-2}+V_{2J}y_{2J}^{J-2}\alpha_J & V_{2J-1}P(y_{2J-1})+V_{2J}P(y_{2J})\alpha_J
    \end{vmatrix}.
\end{multline}
where for the second equality we used the left-shift cyclic permutation matrix to move the first column to the last.
Since $P(y)$ is a monic polynomial of degree $J-1$, replacing the last column by itself plus a linear combination of the remaining columns yields
\begin{multline}
    \frac{\det(\mathbf{B}^0(t,x;\epsilon))}{\det(\mathbf{C})Q\det(\mathbf{S})}\\
    =(-1)^J\begin{vmatrix} V_1 + V_2\alpha_1 &V_1y_1+V_2y_2\alpha_1 &\cdots  & V_1 y_1^{J-1} +V_2 y_2^{J-1}\alpha_1\\
    \vdots &\vdots & \ddots  & \vdots\\
    V_{2J-1} + V_{2J}\alpha_J &
    V_{2J-1}y_{2J-1}+V_{2J}y_{2J}\alpha_J& \cdots   & 
    V_{2J-1}y_{2J-1}^{J-1}+V_{2J}y_{2J}^{J-1}\alpha_J 
    \end{vmatrix}.
\label{eq:detB-3}
\end{multline}
Let us define diagonal matrices by
\begin{equation}
    \begin{split}
        \mathbf{D}_\mathrm{o}&:=\mathrm{diag}(V_1,V_3,\dots,V_{2J-1}), \\
        \mathbf{D}_\mathrm{e}&:=\mathrm{diag}(V_2,V_4,\dots,V_{2J}), \\
        \mathbf{D}_\alpha&:=\mathrm{diag}(\alpha_1,\dots,\alpha_J),
    \end{split}
\end{equation}
and two Vandermonde matrices by
\begin{equation}
    \mathbf{V}_\mathrm{o}:=\mathbf{\Delta}(y_1,y_3,\dots,y_{2J-1})=\begin{bmatrix}1 & y_1 &\cdots &y_1^{J-1}\\
    \vdots &\vdots & \ddots & \vdots\\
    1& y_{2J-1} & \cdots & y_{2J-1}^{J-1}
    \end{bmatrix}
\end{equation}   
and similarly $\mathbf{V}_\mathrm{e}:=\mathbf{\Delta}(y_2,y_4,\dots,y_{2J})$.
Then, we can write
\begin{equation}
    (-1)^J \frac{\det(\mathbf{B}^0(t,x;\epsilon))}{\det(\mathbf{C})Q\det(\mathbf{S})}
    \\
    =\det\left(\mathbf{D}_\mathrm{o}\mathbf{V}_\mathrm{o} + \mathbf{D}_\mathrm{e}\mathbf{D}_\alpha\mathbf{V}_\mathrm{e}\right).
    \label{eq:B-to-V}
\end{equation}
Now the matrix $\mathbf{V}_\mathrm{e}$ has an explicit inverse and nonzero determinant, so
\begin{equation}
\det(\mathbf{D}_\mathrm{o}\mathbf{V}_\mathrm{o} + \mathbf{D}_\mathrm{e}\mathbf{D}_\alpha\mathbf{V}_\mathrm{e})=\det(\mathbf{D}_\mathrm{e}\mathbf{D}_\alpha +\mathbf{D}_\mathrm{o}\mathbf{V}_\mathrm{o}\mathbf{V}_\mathrm{e}^{-1})\det(\mathbf{V}_\mathrm{e}),\quad\det(\mathbf{V}_\mathrm{e})=\prod_{1\le i<j\le J}(y_{2j}-y_{2i}).
\end{equation}
The inverse of $\mathbf{V}_\mathrm{e}$ is constructed from Lagrange interpolating polynomials for the points $(y_2,y_4,\dots,y_{2J})$.  It is then easy to show that the matrix elements of $\mathbf{V}_\mathrm{o}\mathbf{V}_\mathrm{e}^{-1}$ are given explicitly by
\begin{equation}
\left(\mathbf{V}_\mathrm{o}\mathbf{V}_\mathrm{e}^{-1}\right)_{jk}=\mathop{\prod_{i=1}^J}_{i\neq k}\frac{y_{2j-1}-y_{2i}}{y_{2k}-y_{2i}}=
\prod_{i=1}^J(y_{2j-1}-y_{2i})\cdot\frac{1}{y_{2j-1}-y_{2k}}\cdot\mathop{\prod_{i=1}^J}_{i\neq k}\frac{1}{y_{2k}-y_{2i}}.
\end{equation}
The left-hand factor depends on the index $j$ only, while the right-hand factor depends on the index $k$ only.  We define related diagonal matrices by
\begin{equation}
\begin{split}
\widetilde{\mathbf{D}}_\mathrm{o}&:=\mathrm{diag}\left(V_1\prod_{i=1}^J(y_1-y_{2i}),V_3\prod_{i=1}^J(y_3-y_{2i}),\dots,V_{2J-1}\prod_{i=1}^J(y_{2J-1}-y_{2i})\right), \\
\widetilde{\mathbf{D}}_\mathrm{e}&:=\mathrm{diag}\left(V_2\mathop{\prod_{i=1}^J}_{i\neq 1}(y_{2}-y_{2i}),V_4\mathop{\prod_{i=1}^J}_{i\neq 2}(y_4-y_{2i}),\dots,V_{2J}\mathop{\prod_{i=1}^J}_{i\neq J}(y_{2J}-y_{2i})\right),
\end{split}
\end{equation}
and a $J\times J$ Cauchy matrix $\widetilde{\mathbf{C}}$ with elements
\begin{equation}
    \widetilde{C}_{jk}:=\frac{1}{y_{2j-1}-y_{2k}},\quad 1\le j,k\le J.
\end{equation}
Notice that
\begin{align}
    \widetilde{\mathbf{D}}_\mathrm{o}^{-1}\mathbf{D}_\alpha\widetilde{\mathbf{D}}_\mathrm{e}+\widetilde{\mathbf{C}}=\frac{-1}{2t}\widetilde{\mathbf{M}}(t,x;\epsilon)
\end{align}
by use of the identity $y_j-y_k=-2t(u_{j}^{\mathrm{B}}-u_{k}^{\mathrm{B}})$.  Now we see that
\begin{equation}
\begin{split}
\det(\mathbf{D}_\mathrm{o}\mathbf{V}_\mathrm{o} + \mathbf{D}_\mathrm{e}\mathbf{D}_\alpha\mathbf{V}_\mathrm{e})&=\det(\mathbf{D}_\mathrm{e}\mathbf{D}_\alpha + \widetilde{\mathbf{D}}_\mathrm{o}\widetilde{\mathbf{C}}\widetilde{\mathbf{D}}_\mathrm{e}^{-1}\mathbf{D}_\mathrm{e})\det(\mathbf{V}_\mathrm{e})\\
&=\det(\widetilde{\mathbf{D}}_\mathrm{o}(\widetilde{\mathbf{D}}_\mathrm{o}^{-1}\mathbf{D}_\alpha\widetilde{\mathbf{D}}_\mathrm{e} + \widetilde{\mathbf{C}})\widetilde{\mathbf{D}}_\mathrm{e}^{-1}\mathbf{D}_\mathrm{e})\det(\mathbf{V}_\mathrm{e})\\
&=(-2t)^{-J}\det(\widetilde{\mathbf{M}})\det(\widetilde{\mathbf{D}}_\mathrm{o})\det(\widetilde{\mathbf{D}}_\mathrm{e}^{-1})\det(\mathbf{D}_\mathrm{e})\det(\mathbf{V}_\mathrm{e}).
\end{split}
\end{equation}
Therefore, going back to~\eqref{eq:B-to-V} we have
\begin{equation}\label{eq:B-to-M}
    \det(\mathbf{B}^0(t,x;\epsilon))=(2t)^{-J}\det(\mathbf{C})Q\det(\mathbf{S})\frac{\displaystyle\prod_{j=1}^JV_{2j-1}\prod_{i=1}^J(y_{2j-1}-y_{2i})}{\displaystyle\prod_{1\le j<i\le J}(y_{2j}-y_{2i})}\det(\widetilde{\mathbf{M}}).
\end{equation}
Using the fact that the factors other than $\det(\widetilde{\mathbf{M}})$ are independent of~$\epsilon$ and differentiable, we see that~\eqref{eq:super-ultimate} holds by applying $\epsilon\mathrm{Im}\partial x\log(\diamond)$ to both sides of the above equation.
\end{proof}

\subsection{Simplification of diagonal elements of \texorpdfstring{$\widetilde{\mathbf{M}}$}{tildeM}}
The final step required to prove Theorem~\ref{thm:u-app} is to show that $\widetilde{\mathbf{M}}(t,x;\epsilon)=\mathbf{M}(t,x;\epsilon)$.  By comparing \eqref{eq:tildeM} and \eqref{eq:Mjk}, it is obvious that we only need to show the equality of their diagonal elements.  This is achieved in the following two lemmas.

\begin{lemma}\label{lem:abstildeM}
We have the identity
\begin{align}
    \left|\widetilde{m}_j(t,x;\epsilon)\right|^2=\left|\gamma_j(t,x)\right|^2, \quad j=1,\ldots,J.
\end{align}
\end{lemma}

\begin{proof}
We recall from~\eqref{eq:DK-cj} that
\begin{equation}\label{eq:absgamj2}
    \left|\widetilde{m}_j(t,x;\epsilon)\right|^2 = \frac{|\alpha_j|^2|V_{2j}|^2}{|V_{2j-1}|^2}\cdot\frac{\displaystyle\mathop{\prod_{i=1}^J}_{i\neq j}(u_{2j}^\mathrm{B}(t,x)-u_{2i}^\mathrm{B}(t,x))^2}{\displaystyle\prod_{i=1}^J(u_{2j-1}^\mathrm{B}(t,x)-u_{2i}^\mathrm{B}(t,x))^2}.
\end{equation}
Using~\eqref{eq:const-Ql} and~\eqref{eq:VjP-define}, we find
\begin{equation}
    \frac{|V_{2j}|^2}{|V_{2j-1}|^2}=\frac{\displaystyle (y_{2j}-y_0)^2\prod_{m=1}^{N-J}(y_{2j}-z_m)(y_{2j}-z_m^*)\cdot\prod_{n=1}^N(y_{2j-1}-p_n)(y_{2j-1}-p_n^*)}{\displaystyle(y_{2j-1}-y_0)^2\prod_{m=1}^{N-J}(y_{2j-1}-z_m)(y_{2j-1}-z_m^*)\cdot\prod_{n=1}^N(y_{2j}-p_n)(y_{2j}-p_n^*)}.
    \label{eq:V2j/V2j-1}
\end{equation}
Recall, from~\eqref{eq:alpha-j}, that
\begin{equation}
    |\alpha_j|^2 = -\frac{h''(y_{2j-1})}{h''(y_{2j})}.
\end{equation}
Using the factorization of $h'(z)$ in~\eqref{eq:hprimefactors},
\begin{equation}
    |\alpha_j|^2=-\frac{\displaystyle\mathop{\prod_{k=0}^{2J}}_{k\neq 2j-1}(y_{2j-1}-y_k)\cdot\prod_{m=1}^{N-J}(y_{2j-1}-z_m)(y_{2j-1}-z_m^*)\cdot\prod_{n=1}^N(y_{2j}-p_n)(y_{2j}-p_n^*)}{\displaystyle\mathop{\prod_{k=0}^{2J}}_{k\neq 2j}(y_{2j}-y_k)\cdot\prod_{m=1}^{N-J}(y_{2j}-z_m)(y_{2j}-z_m^*)\cdot\prod_{n=1}^N(y_{2j-1}-p_n)(y_{2j-1}-p_n^*)}.
\end{equation}
So combining with~\eqref{eq:V2j/V2j-1} and using the identity~\eqref{eq:yjdiff},
\begin{align*}
    \frac{|\alpha_j|^2|V_{2j}|^2}{|V_{2j-1}|^2}&=
    -\frac{\displaystyle (y_{2j}-y_0)\mathop{\prod_{k=1}^{2J}}_{k\neq 2j-1}(y_{2j-1}-y_k)}{\displaystyle (y_{2j-1}-y_0)\mathop{\prod_{k=1}^{2J}}_{k\neq 2j}(y_{2j}-y_k)}=-\frac{\displaystyle(u_{2j}^\mathrm{B}-u_0^\mathrm{B})\mathop{\prod_{k=1}^{2J}}_{k\neq 2j-1}(u_{2j-1}^\mathrm{B}-u_k^\mathrm{B})}{\displaystyle (u^\mathrm{B}_{2j-1}-u_0^\mathrm{B})\mathop{\prod_{k=1}^{2J}}_{k\neq 2j}(u_{2j}^\mathrm{B}-u_k^\mathrm{B})}.
\end{align*}
Returning to~\eqref{eq:absgamj2}, we have shown that
\begin{align}
    \left|\widetilde{m}_j(t,x;\epsilon)\right|^2=-\frac{\displaystyle (u^\mathrm{B}_{2j-1}-u_0^\mathrm{B})\mathop{\prod_{i=1}^{J}}_{i\neq j}(u_{2j}^\mathrm{B}-u_{2i}^\mathrm{B})\mathop{\prod_{k=1}^{J}}_{k\neq j}(u_{2j-1}^\mathrm{B}-u_{2k-1}^\mathrm{B})}{\displaystyle(u_{2j}^\mathrm{B}-u_0^\mathrm{B})\mathop{\prod_{i=1}^{J}}(u_{2j-1}^\mathrm{B}-u_{2i}^\mathrm{B})\mathop{\prod_{k=1}^{J}}(u_{2j}^\mathrm{B}-u_{2k-1}^\mathrm{B})},
\end{align}
which is of the form~\eqref{eq:mod-gamma-exact} with the identification $R_k=u_k^{\mathrm{B}}(t,x)$ for $k=0,1,\ldots,2J$.
\end{proof}

We next express $\varphi_j(t,x):=\arg(\widetilde{m}_j(t,x;\epsilon))-\theta_j(t,x)/\epsilon$ explicitly in terms of the branches of the solution of the inviscid Burgers equation with initial data $u_0$.  We observe that the ratio of products involving the functions $u_k^\mathrm{B}(t,x)$ appearing in~\eqref{eq:tildeM} is real-valued and negative.  From~\eqref{eq:alpha-j} and~\eqref{eq:VjP-define}, and using that $z_0=y_0\in\mathbb{R}$, we then get that
\begin{equation}
    \varphi_j(t,x)=\frac{\pi}{2} +\Phi(y_{2j-1}(t,x))-\Phi(y_{2j}(t,x)),
\label{eq:varphi-j}
\end{equation}
wherein
\begin{equation}
    \Phi(y):=\sum_{n=1}^N\arg(y-p_n)-\sum_{m=1}^{N-J}\arg(y-z_m),\quad y\in\mathbb{R}.
    \label{eq:Phi-define}
\end{equation}
Here we consistently choose the principal branch, so each of the arguments lies in the range $(-\pi,0)$ because the points $\{p_n\}_{n=1}^N$ and $\{z_m\}_{m=1}^{N-J}$ lie in the open upper half-plane.

\begin{lemma}\label{lem:Phi}
We have the identity
\begin{equation}\label{eq:lemPhi}
    \Phi(y)=-J\frac{\pi}{2}+\frac{1}{2\pi}\int_0^{+\infty}\ln\left(\frac{g(y-s)}{g(y+s)}\right)\, \frac{\dd s}{s},\quad y\in\mathbb{R},
\end{equation}
where
\begin{equation}
    g(y)=\frac{y+2tu_0(y)-x}{\displaystyle{\prod_{k=0}^{2J}(y+2tu_k^\mathrm{B}(t,x)-x)}}.
\end{equation}
\end{lemma}

\begin{proof}
Since for $k=0,\dots,2J$, $u_k^\mathrm{B}(t,x)=u_0(y_k(t,x))$ where $y_k(t,x)$ is a real critical point of $h(z)$, i.e., a solution of $z + 2tu_0(z)-x=0$, $g(y)$ can be written in the alternate form
\begin{equation}\label{eq:g-rewrite}
    g(y)
    =\frac{2th'(y)}{\displaystyle\prod_{k=0}^{2J}(y-y_k(t,x))}
    =\frac{\displaystyle\prod_{m=1}^{N-J}(y-z_m)(y-z_m^*)}{\displaystyle\prod_{n=1}^N(y-p_n)(y-p_n^*)}.
\end{equation}
The second equality was obtained via \eqref{eq:hprimefactors}.  Clearly $g(y)$ is strictly positive for $y\in\mathbb{R}$.  The same is true of the product $(y^2+1)^Jg(y)$, which has the additional property that $(y^2+1)^Jg(y)=1+\mathcal{O}(y^{-1})$ as $y\to\infty$.  It follows that
\begin{equation}
L(y):=\ln((y^2+1)^J g(y))
\end{equation}
is a well-defined real-valued function of $y\in\mathbb{R}$ that is in $L^2(\mathbb{R})$. Therefore, it makes sense to apply the complementary Cauchy-Szeg\H{o} projection $\mathrm{Id}-\Pi$, where $\mathrm{Id}$ is the identity operator and $\Pi$ denotes the orthogonal projection from $L^2(\mathbb{R})$ onto the Hardy subspace of functions analytic in the upper half-plane.  Given the form of $g(y)$ as a rational function from~\eqref{eq:g-rewrite}, it is easy to do this explicitly by taking the terms from $L$ that are analytic in the lower half-plane
\begin{equation}   
({\rm Id}-\Pi )L(y)=J\log (y-\ii)+\sum_{m=1}^{N-J} \log (y-z_m)-\sum_{n=1}^N \log (y-p_n),\quad y\in\mathbb{R},
\end{equation}
where the complex logarithms denote the principal branches.  The function $\Phi(y)$ defined in~\eqref{eq:Phi-define} can then be recognized from the imaginary part of the right-hand side, leading to the formula
\begin{equation}
    \Phi(y)=J\arg(y-\ii)-\mathrm{Im}((\mathrm{Id}-\Pi)L(y)).
\end{equation}
Since $\mathrm{Im}(\log(y-\ii))=\arg(y-i)\in (-\pi,0)$ for $y\in\mathbb{R}$, we can write $\arg(y-\ii)=\arctan(y)-\pi/2$.

It remains to simplify the integral term.  We now use the fact that not only is $L\in L^2(\mathbb{R})$ but also $L\in C^1(\mathbb{R})$.  We note that for every real-valued function $f\in L^2(\mathbb{R})\cap C^1(\mathbb{R})$ we have
\begin{equation}
\mathrm{Im}((\mathrm{Id}-\Pi )f(y))=\frac{1}{2\pi}\int_0^{+\infty} \frac{f(y+s)-f(y-s)}{s}\, \dd s,
\end{equation}
see for instance~\cite[Lemma 3.5]{Badreddine24zero} for a proof.
We obtain
\begin{equation}\label{eq:Phi-rewrite}
\Phi (y)=J\left (\arctan (y)-\frac \pi 2\right )+\frac{J}{2\pi}\int_0^{+\infty}\ln\left(\frac{(y-s)^2+1}{(y+s)^2+1}\right)\, \frac{\dd s}{s}  +\frac{1}{2\pi}\int_0^{+\infty}\ln\left( \frac{g(y-s)}{g(y+s)}\right)\, \frac{\dd s}{s}.
\end{equation}
It is easy to check that 
\[
f_0(y):=\frac{1}{2\pi}\int_0^{+\infty}\ln\left(\frac{(y-s)^2+1}{(y+s)^2+1}\right)\, \frac{\dd s}{s}
\]
satisfies $f_0(0)=0$ and 
\begin{align}
\begin{split}
f_0'(y)=&\frac{1}{\pi}\int_0^{+\infty} \left[\frac{y-s}{(y-s)^2+1}- \frac{y+s}{(y+s)^2+1}\right]\, \frac{\dd s}{s}\\
=&\frac{1}{\pi}\int_{\mathbb{R}} \frac{(y^2-s^2-1)\, \dd s}{((y-s)^2+1)((y+s)^2+1)}=-\frac{1}{1+y^2},
\end{split}
\end{align}
so that $f_0(y)=-\arctan(y)$.
Using this in~\eqref{eq:Phi-rewrite} gives~\eqref{eq:lemPhi} and completes the proof.
\end{proof}

According to Lemmas \ref{lem:abstildeM}, \ref{lem:Phi} and \eqref{eq:gamma}, \eqref{eq:Mjk},
\begin{align}\label{eq:tildeM=M}
    \widetilde{m}_j=|\widetilde{m}_j|\ee^{\ii\arg\widetilde{m}_j}=|\gamma_j|\ee^{\ii\phi_j+\ii\theta_j/\epsilon}=\gamma_j\ee^{\ii\theta_j/\epsilon},
\end{align}
so it is clear that $\widetilde{\mathbf{M}}(t,x;\epsilon)=\mathbf{M}(t,x;\epsilon)$ by comparing \eqref{eq:tildeM} and \eqref{eq:Mjk}.  Finally, we conclude that Theorem \ref{thm:u-app} holds for $J(t,x)\geq1$ by Propositions \ref{prop:ultimate}, \ref{prop:super-ultimate}, and the equality $\widetilde{\mathbf{M}}(t,x;\epsilon)=\mathbf{M}(t,x;\epsilon)$.

\begin{remark}[Nonvanishing]
Finally, we note that thanks to~\cite[Lemma 1.1]{DobrokhotovK91}, $\mathbf{M}(t,x;\epsilon)$ is nonsingular, which implies that $\det(\mathbf{B}^0(t,x;\epsilon))$ is bounded away from zero locally for $(t,x)\in\mathbb{R}^2$ for which $J(t,x)=1,2,\dots$ is fixed, as was required to obtain the formula~\eqref{eq:ultimate} to begin with.
\label{rem:B-not-zero}
\end{remark}

\section{Convergence in \texorpdfstring{$L^2$}{L2}}\label{section:CV-L2}

To prove Corollary~\ref{cor:CV-L2}, we must demonstrate that
\begin{align}
    \|u(t,\diamond;\epsilon)-u^\mathrm{ZD}(t,\diamond;\epsilon)\|_{L^2(\mathbb{R})}^2\to0 \quad \text{as } \epsilon\to0.
\end{align}
It is clear to see that
\begin{multline}\label{eq:normu-uZD}
\|u(t,\diamond;\epsilon)-u^\mathrm{ZD}(t,\diamond;\epsilon)\|_{L^2(\mathbb{R})}^2
	=\|u(t,\diamond;\epsilon)\|_{L^2(\mathbb{R})}^2-\|u^\mathrm{ZD}(t,\diamond;\epsilon)\|_{L^2(\mathbb{R})}^2\\
	-2\left\langle u(t,\diamond;\epsilon)-u^\mathrm{ZD}(t,\diamond;\epsilon), u^\mathrm{ZD}(t,\diamond;\epsilon)\right\rangle.
\end{multline}
Since $t\mapsto u(t,\diamond;\epsilon)$ is the solution to the Benjamin-Ono equation~\eqref{eq:BO}, which conserves the $L^2$ norm, we have $\|u(t,\diamond;\epsilon)\|_{L^2(\mathbb{R})}^2=\|u_0\|_{L^2(\mathbb{R})}^2$.  Thus, we prove Corollary~\ref{cor:CV-L2} by showing that
\begin{align}\label{eq:uzdnormtou0norm}
    \|u^{\mathrm{ZD}}(t,\diamond;\epsilon)\|_{L^2(\mathbb{R})}^2\to\|u_0\|_{L^2(\mathbb{R})}^2 \quad \text{as } \epsilon\to0,
\end{align}
and
\begin{align}\label{eq:uuzdipto0}
    \left\langle u(t,\diamond;\epsilon)-u^\mathrm{ZD}(t,\diamond;\epsilon), u^\mathrm{ZD}(t,\diamond;\epsilon)\right\rangle\to0 \quad \text{as }\epsilon\to0.
\end{align}
We prove the statements~\eqref{eq:uzdnormtou0norm}, \eqref{eq:uuzdipto0} in Propositions~\ref{lem:norm-uapp}, \ref{prop:uuzdipto0}, respectively.  We begin with three preparatory lemmas.

\begin{lemma}\label{lem:uapp-bound}
If $t\ge 0$ is such that  $J(t,x)\leq 1$ for every generic $x$ where $J(t,x)$ is well-defined, then
\begin{equation}
|u^\mathrm{ZD}(t,x;\epsilon)|
	\leq 9\|u_0\|_{L^\infty(\mathbb{R})}.
 \label{eq:J1L-infty}
\end{equation}
\end{lemma}

\begin{proof}
Because $J(t,x)=0$ implies that $u^\mathrm{ZD}(t,x;\epsilon)=u_0^\mathrm{B}(t,x)$ of the inviscid Burgers equation with initial data $u_0$, the use of~\eqref{eq:uB-yj} yields the desired inequality~\eqref{eq:J1L-infty} without the factor of $9$. Therefore, it is enough to consider a point at which $J(t,x)=1$ and hence $u^\mathrm{B}(t,\cdot)$ has three branches. We note that since $|\cos(\cdot)|\leq 1$, we have
\begin{equation}
|u^\mathrm{ZD}(t,x;\epsilon)|
	\leq |u^\mathrm{B}_0(t,x)|+(u^\mathrm{B}_2(t,x)-u^\mathrm{B}_1(t,x))\frac{1-r(t,x)^2}{(1-r(t,x))^2}.
\end{equation}
Using the identity $(1-r^2)/(1-r)^2=(1+r)^2/(1-r^2)$ and~\eqref{eq:r} which also implies that
\begin{equation}
    1-r(t,x)^2=\frac{u_2^\mathrm{B}(t,x)-u_1^\mathrm{B}(t,x)}{u_2^\mathrm{B}(t,x)-u_0^\mathrm{B}(t,x)},
\end{equation}
we obtain
\begin{equation}
\begin{split}
|u^\mathrm{ZD}(t,x;\epsilon)|&
	\leq |u^\mathrm{B}_0(t,x)|+(u^\mathrm{B}_2(t,x)-u^\mathrm{B}_0(t,x))\left(1+\sqrt{\frac{u^\mathrm{B}_1(t,x)-u^\mathrm{B}_0(t,x)}{u^\mathrm{B}_2(t,x)-u^\mathrm{B}_0(t,x)}}\right)^2\\
 &=|u_0^\mathrm{B}(t,x)| + \left(\sqrt{u_2^\mathrm{B}(t,x)-u_0^\mathrm{B}(t,x)}+\sqrt{u_1^\mathrm{B}(t,x)-u_0^\mathrm{B}(t,x)}\right)^2\\
 &\le |u_0^\mathrm{B}(t,x)| + \left(\sqrt{|u_2^\mathrm{B}(t,x)|+|u_0^\mathrm{B}(t,x)|} +\sqrt{|u_1^\mathrm{B}(t,x)|+|u_0^\mathrm{B}(t,x)|}\right)^2.
\end{split}
\end{equation}
The desired inequality~\eqref{eq:J1L-infty} then follows from~\eqref{eq:uB-yj}.
\end{proof}

We have $J(t,x)=0$ except for $x$ in a finite union of pairwise disjoint closed intervals that we denote by $[X^-_k(t),X^+_k(t)]$ for $k=1,\dots,K$, with $X^-_1(t)\le X^+_1(t)<X^-_2(t)\le X^+_2(t)<\cdots<X^-_K(t)\le X^+_K(t)$, and $J(t,x)=1$ holds on each open interval $(X^-_k(t),X^+_k(t))$ (possibly empty if $X^-_k(t)=X^+_k(t)$).  We have $K=0$ if and only if $J(t,x)=0$ unambiguously for all $x\in\mathbb{R}$. Let $\delta>0$ and define the set 
\begin{equation}\label{eq:Sdelta}
S_\delta(t):=\mathbb{R}\setminus\bigcup_{k=1}^K\left((X_k^-(t)-\delta,X_k^-(t)+\delta)\cup (X_k^+(t)-\delta,X_k^+(t)+\delta)\right).
\end{equation}
Now, $S_\delta(t)$ is a finite union of pairwise disjoint closed intervals on each of which $J(t,x)=0$ or $J(t,x)=1$, and exactly two of the intervals on which $J(t,x)=0$ are infinite. We set $S_\delta(t)=S_\delta^0(t)\sqcup S_\delta^1(t)$, where $S_\delta^0(t)$ and $S_\delta^1(t)$ denote the unions of intervals in $S_\delta(t)$ on which $J(t,x)=0$ and $J(t,x)=1$ respectively. Therefore
\begin{equation}
\|u^\mathrm{ZD}(t,\diamond;\epsilon)\|_{L^2(\mathbb{R})}^2
	=\|u^\mathrm{ZD}(t,\diamond;\epsilon)\|_{L^2(S_\delta^0)}^2+\|u^\mathrm{ZD}(t,\diamond;\epsilon)\|_{L^2(S_\delta^1)}^2+r_\delta,
 \label{eq:carve-out-delta}
\end{equation}
where $r_\delta\leq C\delta$ tends to $0$ as $\delta\to 0$ independently of $\epsilon$ thanks to Lemma~\ref{lem:uapp-bound}. When $x$ belongs to $S_\delta^1(t)$, one can check that $r(t,x)$ is bounded away from $0$ and $1$: $0<\eta_\delta\leq r(t,x)<1-\eta_\delta<1$ for some small $\eta_\delta>0$.  For such $x$, we write $u^\mathrm{ZD}(t,x;\epsilon)^2$ in the form
\begin{equation}
    u^\mathrm{ZD}(t,x;\epsilon)^2=u_0^\mathrm{B}(t,x)^2 +2F^{(1)}(t,x,\epsilon^{-1}\theta(t,x)) + F^{(2)}(t,x,\epsilon^{-1}\theta(t,x)),\quad x\in S^1_\delta(t),
\label{eq:uZDsquared}
\end{equation}
where 
\begin{equation}
F^{(p)}(t,x,y):=u^\mathrm{B}_0(t,x)^{2-p}(u^\mathrm{B}_2(t,x)-u^\mathrm{B}_1(t,x))^pU_{r(t,x)}(y+\varphi(t,x))^p,\quad p=1,2
\end{equation}
are functions that are $2\pi$-periodic in $y$ and represent the terms linear and quadratic in $U_r$.

\begin{lemma}\label{lem:averaging}
Let $\delta>0$ be fixed and sufficiently small.  Then, on each maximal subinterval $[X^-_k(t)+\delta,X^+_k(t)-\delta]\subset S^1_\delta(t)$,
\begin{equation}\label{eq:zero-fourier}
\int_{X^-_k(t)+\delta}^{X^+_k(t)-\delta}F^{(p)}(t,x,\epsilon^{-1}\theta(t,x))\dd x
	=\int_{X^-_k(t)+\delta}^{X^+_k(t)-\delta}\langle F^{(p)}\rangle(t,x)\,\dd x
	+\mathcal{O}_\delta(\epsilon),\quad\epsilon\to 0,
\end{equation}
where $\langle F^{(p)}\rangle(t,x)$ denotes the period average over $y$:
\begin{equation}
    \langle F^{(p)}\rangle(t,x):=\frac{1}{2\pi}\int_0^{2\pi} F^{(p)}(t,x,y)\,\dd y.
\end{equation}
\end{lemma}

\begin{proof}
We expand $F^{(p)}(t,x,y)$ into a Fourier series
\begin{equation}
F^{(p)}(t,x,y)
	=\sum_{n\in\mathbb{Z}}F^{(p)}_n(t,x)\ee^{\ii ny},\quad p=1,2,
\end{equation}
where
\begin{equation}
F^{(p)}_n(t,x)
	=c_n^{(p)}(r(t,x))u^\mathrm{B}_0(t,x)^{2-p}(u^\mathrm{B}_2(t,x)-u^\mathrm{B}_1(t,x))^{p},
\end{equation}
and
\begin{align}
    c_n^{(p)}(r)=\frac{1}{2\pi}\int_{0}^{2\pi}\left(\frac{1-r^2}{1-2\cos(y)+r^2}\right)^{p}\ee^{-\ii ny}\dd y.
\end{align}
Using the change of variable $z=\ee^{-\ii\mathrm{sgn}(n)y}$, where $\mathrm{sgn}(n)=1$ for $n\geq0$ and $\mathrm{sgn}(n)=-1$ for $n<0$,
\begin{equation}\label{eq:cnr}
    c_{n}^{(p)}(r)=\frac{(r^{-1}-r)^{p}}{2\pi\ii}\oint_{|z|=1}\frac{z^{|n|+p-1}}{(z-r)^{p}(r^{-1}-z)^{p}}\dd z=r^{|n|}\left(|n|+\frac{1+r^2}{1-r^2}\right)^{p-1}
\end{equation}
by residue calculus, since $0<r(t,x)<1$. In both cases, we have
\begin{equation}
    \sum_{n\in\mathbb{Z}} |F^{(p)}_n(t,\diamond)|\in L^1([X_k^-(t)+\delta,X_k^+(t)-\delta]),
\end{equation}
and consequently one can exchange the order of sum and integral:
\begin{equation}\label{eq:sumfourier}
\int_{X_k^-(t)+\delta}^{X_k^+(t)-\delta} F^{(p)}(t,x,\epsilon^{-1}\theta(t,x))\,\dd x
	=\sum_{n\in\mathbb{Z}}\int_{X_k^-(t)+\delta}^{X_k^+(t)-\delta} F^{(p)}_n(t,x)\ee^{\ii n\theta(t,x)/\epsilon}\,\dd x.
\end{equation}
Integrating by parts, for every $n\neq 0$,
\begin{multline*}
\int_{X_k^-(t)+\delta}^{X_k^+(t)-\delta} F^{(p)}_n(t,x)\ee^{\ii n\theta(t,x)/\epsilon}\,\dd x\\
	=\epsilon\left[\frac{F^{(p)}_n(t,x)}{\ii n\theta_x(t,x)}\ee^{\ii n\theta(t,x)/\epsilon}\right]_{X_k^-(t)+\delta}^{X_k^+(t)-\delta}
	-\epsilon\int_{X_k^-(t)+\delta}^{X_k^+(t)-\delta} \frac{\dd}{\dd x}\left(\frac{F^{(p)}_n(t,x)}{\ii n\theta_x(t,x)}\right)\ee^{\ii n\theta(t,x)/\epsilon}\,\dd x.
\end{multline*}
Notice $\theta_x(t,x)=u_1^\mathrm{B}(t,x)-u_2^\mathrm{B}(t,x)$ is bounded away from zero on $x\in[X_k^-(t)+\delta,X_k^+(t)-\delta]$ and from the formula~\eqref{eq:r} of $r(t,x)$, one can check that $\partial_x r(t,x)\to +\infty$ as $x\to X_k^-(t)$ and $x\to X_k^+(t)$, so 
that $\partial_x r(t,x)\geq c>0$ on $(X_k^-(t),X_k^+(t))$.  Now, using~\eqref{eq:cnr}, it can be shown that
\begin{align}
    \left|\int_{X_k^-(t)+\delta}^{X_k^+(t)-\delta} F^{(p)}_n(t,x)\ee^{\ii n\theta(t,x)/\epsilon}\,\dd x\right|\leq \epsilon \widetilde{K}_{\delta}(t)r(t,x_{k}(\delta))^{|n|-1},
\end{align}
where $x_k(\delta)\in[X_k^-(t)+\delta,X_k^+(t)-\delta]$ maximizes $r(t,x)$ and $\widetilde{K}_{\delta}(t)$ is independent of $x$, $\epsilon$, and $n$.  We conclude that the sum on the right-hand side of~\eqref{eq:sumfourier} is absolutely convergent and thus
\begin{equation}
\mathop{\sum_{n\in\mathbb{Z}}}_{n\neq 0}\int_{X_k^-(t)+\delta}^{X_k^+(t)-\delta} F^{(p)}_n(t,x)\ee^{\ii n\theta(t,x)/\epsilon}\,\dd x
	=\mathcal{O}_\delta(\epsilon),\quad\epsilon\to 0.
\end{equation}
Finally, we get
\begin{equation}
\int_{X_k^-(t)+\delta}^{X_k^+(t)-\delta} F^{(p)}(t,x,\epsilon^{-1}\theta(t,x))\,\dd x
	= \int_{X_k^-(t)+\delta}^{X_k^+(t)-\delta} F^{(p)}_0(t,x)\,\dd x+\mathcal{O}_\delta(\epsilon),\quad\epsilon\to 0,\quad p=1,2,
\end{equation}
which amounts to~\eqref{eq:zero-fourier}.
\end{proof}

Now we use Lemma~\ref{lem:averaging},~\eqref{eq:uZDsquared}, and~\eqref{eq:cnr} to obtain
\begin{equation}\label{eq:int-over-Sdelta1}
    \|u^\mathrm{ZD}(t,\diamond;\epsilon)\|_{L^2(S_\delta^1)}^2=\int_{S_\delta^1(t)}\left[u_2^\mathrm{B}(t,x)^2-u_1^\mathrm{B}(t,x)^2+u_0^\mathrm{B}(t,x)^2\right]\,\dd x + \mathcal{O}_\delta(\epsilon).
\end{equation}
On the complementary part $S^0_\delta(t)$ of $S_\delta(t)$, we have $J(t,x)=0$ and hence $u^{\mathrm{ZD}}(t,x;\epsilon)=u_0^\mathrm{B}(t,x)$, where $u_0^\mathrm{B}(t,x)$ is the unique branch of the solution of Burgers' equation.  Therefore, we have the exact equality
\begin{equation}
    \|u^\mathrm{ZD}(t,\diamond;\epsilon)\|_{L^2(S_\delta^0)}^2= \int_{S_\delta^0(t)}u_0^\mathrm{B}(t,x)^2\,\dd x.
\label{eq:int-over-Sdelta0}
\end{equation}
Since all solution branches of the inviscid Burgers equation with initial data $u_0$ are bounded in absolute value by $\|u_0\|_{L^\infty}$, we can replace $S_\delta^1(t)$ and $S_\delta^0(t)$ in~\eqref{eq:int-over-Sdelta1} and~\eqref{eq:int-over-Sdelta0} respectively by the corresponding unions of intervals with $\delta=0$, at the cost of an error term independent of $\epsilon$ that tends to zero with $\delta$.  Absorbing this error into $r_\delta$ from~\eqref{eq:carve-out-delta} we obtain
\begin{equation}
    \|u^\mathrm{ZD}(t,\diamond;\epsilon)\|_{L^2(\mathbb{R})}^2=I(t)
        + \mathcal{O}_\delta(\epsilon) + r_\delta,
\end{equation}
where $I(t)$ denotes the sum of integrals
\begin{equation}
    I(t):=\int_{S_0^0(t)} u_0^\mathrm{B}(t,x)^2\,\dd x + \int_{S_0^1(t)}\left[u_2^\mathrm{B}(t,x)^2-u_1^\mathrm{B}(t,x)^2+u_0^\mathrm{B}(t,x)^2\right]\,\dd x 
\end{equation}
which is independent of both $\epsilon$ and $t$. It therefore follows that for each fixed $\delta>0$
\begin{equation}
\limsup_{\epsilon\to 0} \left|\|u^\mathrm{ZD}(t,\diamond;\epsilon)\|_{L^2(\mathbb{R})}^2
	-I(t)\right|\leq r_\delta,
\end{equation}
and then since $\delta$ can be taken to be arbitrarily small forcing $r_\delta\to 0$ we have proved that
\begin{equation}\label{eq:normuZDI}
    \lim_{\epsilon\to 0} \|u^\mathrm{ZD}(t,\diamond;\epsilon)\|_{L^2(\mathbb{R})}^2=I(t).
\end{equation}
We now state our final lemma.

\begin{lemma}\label{lem:I}
We have the identity $I(t)=\|u_0\|_{L^2(\mathbb{R})}^2$ (so in fact $I(t)$ is independent of $t$). 
\end{lemma}

\begin{proof}
Recalling the points $X_1^-(t)\le X_1^+(t)<X_2^-(t)\le X_2^+(t)<\cdots <X_K^-(t)\le X_K^+(t)$, we first observe that if $X_k^-(t)<X_k^+(t)$, two of the three branches $u_0^\mathrm{B}(t,x)<u_1^\mathrm{B}(t,x)<u_2^\mathrm{B}(t,x)$ of the inviscid Burgers solution defined for $x\in (X_k^-(t),X_k^+(t))$ are related to the single branch $u_0^\mathrm{B}(t,x)$ defined separately on the two disjoint intervals of $S_0^0(t)$ immediately to the left and right. It will be convenient to interpret $X_{k,L}^{\pm}$ as a left sided limit, i.e. $f(X_{k,L}^{\pm})$ is the limiting value of $f(x)$ as $x\to X_{k}^{\pm}$ from the left. Similar for $X_{k,R}^{\pm}$. Indeed, $u_2^\mathrm{B}(t,x)$ on $(X_k^-(t),X_k^+(t))$ is the continuation of $u_0^\mathrm{B}(t,x)$ from $x<X_k^-(t)$; in other words, $u_{0}^{\mathrm{B}}(t,X_{k,L}^{-})=u_0(y_0(t,X_{k,L}^{-}))$ is equal to $u_0(y_2(t,X_{k,R}^{-}))=u_2^{\mathrm{B}}(t,X_{k,R}^{-})$, where $y_j(t,x)$ is a solution of the characteristic equation~\eqref{eq:critical-points} that is left, right differentiable at $x=X_k^-(t)$ for $j=0,2$, respectively.  Similarly, $u_0^\mathrm{B}(t,x)$ on $(X_k^-(t),X_k^+(t))$ is the continuation of $u_0^\mathrm{B}(t,x)$ from $x>X_k^+(t)$ and both functions are equal to $u_0(y_0(t,x))$ at $x=X^+_k(t)$ where $y=y_0(t,x)$ is a smooth solution of~\eqref{eq:critical-points} in a neighborhood of $x=X_k^{+}(t)$. The functions $x\mapsto y_{j}(t,x)$, $j=0,2$, are both monotone increasing, $\partial_x[y_{0}(t,x)]\to +\infty$ as $x\to X_k^-(t)$ from the right, and $\partial_x[y_{2}(t,x)]\to +\infty$ as $x\to X_k^+(t)$ from the left.  Meanwhile the solution $u^\mathrm{B}_1(t,x)$ defined on $(X_k^-(t),X_k^+(t))$ is equal to $u_0(y_1(t,x))$ where $y_1(t,x)$ is a smooth and monotone decreasing solution of~\eqref{eq:critical-points} on this interval for which $\partial_x[y_1(t,x)]\to -\infty$ as $x\to X_k^\pm(t)$ from the interior of $(X_{k}^{-}(t),X_{k}^{+}(t))$.  Finally, we have $y_0(t,X_{k,L}^{-}(t))=y_2(t,X_{k,R}^{-}(t))$, $y_0(t,X_{k,R}^{-}(t))=y_1(t,X_{k,R}^{-}(t))$, and $y_1(t,X_{k,L}^{+}(t))=y_2(t,X_{k,L}^{+}(t))$.  

For $k=1,2,\ldots,K$, we write
\begin{align}
    \int_{X_{k}^{-}(t)}^{X_{k}^{+}(t)}u_2^\mathrm{B}(t,x)^2-u_1^\mathrm{B}(t,x)^2+u_0^\mathrm{B}(t,x)^2\,\dd x&=\sum_{j=0}^{2}(-1)^{j}\int_{X_{k}^{-}(t)}^{X_{k}^{+}(t)}u_0(y_j(t,x))^2\,\dd x,
\end{align}
then use the change of variables $w=y_j(t,x)$ to obtain
\begin{align}\label{eq:ujS1}
    \int_{X_{k}^{-}(t)}^{X_{k}^{+}(t)}u_2^\mathrm{B}(t,x)^2-u_1^\mathrm{B}(t,x)^2+u_0^\mathrm{B}(t,x)^2\,\dd x&=\int_{y_0(t,X_{k,L}^{-})}^{y_0(t,X_{k}^{+})}u_0(w)^2(1+2tu_0'(w))\dd w.
\end{align}
Similarly,
\begin{align}\label{eq:u0S0}
    \int_{X_{k-1}^{+}(t)}^{X_{k}^{-}(t)}u_0^\mathrm{B}(t,x)^2\,\dd x&=\int_{y_0(t,X_{k-1}^{+})}^{y_0(t,X_{k,L}^{-})}u_0(w)^2(1+2tu_0'(w))\dd w,
\end{align}
for $k=1,\ldots,K+1$, where we understand that $X_{0}^{+}(t)=-\infty$ and $X_{K+1}^{-}(t)=+\infty$.  It follows from~\eqref{eq:ujS1}, \eqref{eq:u0S0} that we can rewrite $I(t)$ as
\begin{equation}
    I(t)=\int_\mathbb{R} u_0(w)^2(1+2tu_0'(w))\,\dd w = \|u_0\|_{L^2(\mathbb{R})}^2 + \frac{2}{3}t\int_\mathbb{R}\frac{\dd}{\dd w}u_0(w)^3\,\dd w.
\end{equation}
Since $u_0$ is a rational function in $L^2(\mathbb{R})$, we have $u_0(w)\to 0$ as $w\to\pm\infty$ and hence the second term vanishes by the fundamental theorem of calculus.
\end{proof}

Applying Lemma~\ref{lem:I} to~\eqref{eq:normuZDI} completes the proof of the following proposition.

\begin{proposition}\label{lem:norm-uapp}
If $t\ge 0$ is such that  $J(t,x)\leq 1$ for every generic $x$ where $J(t,x)$ is well-defined, then
\[
\|u^\mathrm{ZD}(t,\diamond;\epsilon)\|_{L^2(\mathbb{R})}\to \|u_0\|_{L^2(\mathbb{R})} \quad \text{as }\epsilon\to0.
\]
\end{proposition}

Returning to~\eqref{eq:normu-uZD}, we see that the proof of Corollary~\ref{cor:CV-L2} will be complete after proving the following Proposition.

\begin{proposition}\label{prop:uuzdipto0}
If $t\ge 0$ is such that  $J(t,x)\leq 1$ for every generic $x$ where $J(t,x)$ is well-defined, then
\[
\left\langle u(t,\diamond;\epsilon)-u^\mathrm{ZD}(t,\diamond;\epsilon), u^\mathrm{ZD}(t,\diamond;\epsilon)\right\rangle\to 0 \quad \text{as }\epsilon\to0.
\]
\end{proposition}

\begin{proof}
Fix $t\ge 0$, $R>0$ large, $\delta>0$ small, and a rational initial condition $u_0$. Let us consider the set $S_\delta(t)$ from~\eqref{eq:Sdelta}, from which we now exclude also small intervals of width $2\delta$ centered at points of the discriminant locus $\mathcal{D}$ at the fixed time $t\ge 0$ that are not points of the caustic locus $\mathcal{C}$.  We denote the resulting smaller set by $\widetilde{S}_\delta(t)$ and denote its complement in $\mathbb{R}$ by $\widetilde{S}_\delta^\mathrm{c}(t)$.  The latter is just the union of intervals of width $2\delta$ centered at \emph{all} points (at most $4N$ in number) of $\mathcal{D}$ for fixed $t\ge 0$.
We decompose the scalar product 
\begin{equation}
\left\langle u(t,\diamond;\epsilon)-u^\mathrm{ZD}(t,\diamond;\epsilon), u^\mathrm{ZD}(t,\diamond;\epsilon)\right\rangle
	=\int_{\mathbb{R}}(u(t,x;\epsilon)-u^\mathrm{ZD}(t,x;\epsilon)) u^\mathrm{ZD}(t,x;\epsilon)\,\dd x
\end{equation}
into three parts by the disjoint union $\mathbb{R}=\{|x|>R\}\sqcup \widetilde{S}_\delta^\mathrm{c}(t)\sqcup ([-R,R]\cap \widetilde{S}_\delta(t))$ and estimate the three terms separately.

First, using the Cauchy-Schwarz and Minkowski inequalities,
\begin{multline}
\left(\int_{|x|>R}(u(t,x;\epsilon)-u^\mathrm{ZD}(t,x;\epsilon)) u^\mathrm{ZD}(t,x;\epsilon)\,\dd x\right)^2\\
\begin{aligned}
	&\leq \int_{|x|>R}(u(t,x;\epsilon)-u^\mathrm{ZD}(t,x;\epsilon))^2\,\dd x
	\cdot\int_{|x|>R}u^\mathrm{ZD}(t,x;\epsilon)^2\,\dd x\\
 &\leq \|u(t,\diamond;\epsilon)-u^\mathrm{ZD}(t,\diamond;\epsilon)\|_{L^2(\mathbb{R})}^2\int_{|x|>R}u^\mathrm{ZD}(t,x;\epsilon)^2\,\dd x\\
 &\le \left[\|u(t,\diamond;\epsilon)\|_{L^2(\mathbb{R})} + \|u^\mathrm{ZD}(t,\diamond;\epsilon)\|_{L^2(\mathbb{R})}\right]^2\int_{|x|>R}u^\mathrm{ZD}(t,x;\epsilon)^2\,\dd x.
 \end{aligned}
\end{multline}
By Proposition~\ref{lem:norm-uapp} and conservation of the $L^2(\mathbb{R})$ norm under Benjamin-Ono, the square of the sum of norms remains bounded as $\epsilon\to 0$. Moreover, for $R$ large enough, we have $J(t,x)=0$, and hence $u^\mathrm{ZD}(t,x;\epsilon)=u^\mathrm{B}_0(t,x)$ is the single branch of the solution of the inviscid Burgers equation with initial data $u_0$.  Since $u_0$ is rational, it is easy to see that $u^\mathrm{B}_0(t,x)=\mathcal{O}(x^{-1})$ as $x\to\pm\infty$. Therefore,
\begin{equation}
\left(\int_{|x|>R}(u(t,x;\epsilon)-u^\mathrm{ZD}(t,x;\epsilon)) u^\mathrm{ZD}(t,x;\epsilon)\,\dd x\right)^2	\leq \frac{C_1}{R}
\label{eq:C1}
\end{equation}
for some $C_1>0$ independent of $\epsilon$ (but possibly dependent on $t$).

Next, we use that $u^\mathrm{ZD}(t,\diamond;\epsilon)$ is bounded in $L^\infty(\mathbb{R})$ thanks to Proposition~\ref{lem:uapp-bound}, while as noted above $u(t,\diamond;\epsilon)-u^\mathrm{ZD}(t,\diamond;\epsilon)$ is bounded in $L^2(\mathbb{R})$, both uniformly as $\epsilon\to 0$, so the Cauchy-Schwarz inequality leads to
\begin{multline}
\Bigg(\int_{\widetilde{S}_\delta^\mathrm{c}(t)}(u(t,x;\epsilon)-u^\mathrm{ZD}(t,x;\epsilon)) 	u^\mathrm{ZD}(t,x;\epsilon)\,\dd x\Bigg)^2\\
\begin{aligned}
	&\leq \int_{\widetilde{S}_\delta^\mathrm{c}(t)}(u(t,x;\epsilon)-u^\mathrm{ZD}(t,x;\epsilon))^2\,\dd x\cdot
	\int_{\widetilde{S}_\delta^\mathrm{c}(t)}u^\mathrm{ZD}(t,x;\epsilon)^2\,\dd x\\
	&\leq C_2 \delta.
 \end{aligned}
 \label{eq:C2}
\end{multline}
because the first factor is bounded independently of $\epsilon$ and $\delta$, while for the second factor we use the fact that for rational initial data, 
$\operatorname{Leb}(\widetilde{S}_\delta^\mathrm{c}(t))\leq C_0 \delta$ for some $C_0>0$ (we can take $C_0=8N$) since there are only finitely many points of $\mathcal{D}$ at time $t$.  The constant $C_2$ is independent of $\epsilon$.

Finally, we use Theorem~\ref{thm:u-app}, asserting that $u(t,\diamond;\epsilon)-u^\mathrm{ZD}(t,\diamond;\epsilon)$ is uniformly convergent to $0$ on the compact set $\{|x|\leq R\}\cap \widetilde{S}_\delta(t)$.  Since also $u^\mathrm{ZD}(t,\diamond;\epsilon)$ is bounded in $L^\infty(\mathbb{R})$ uniformly for small $\epsilon$, we deduce that
\begin{equation}
\lim_{\varepsilon\to 0} \int_{|x|\leq R,\;x\in \widetilde{S}_\delta(t)}(u(t,x;\epsilon)-u^\mathrm{ZD}(t,x;\epsilon)) u^\mathrm{ZD}(t,x;\epsilon)\,\dd x
	=0.
 \label{eq:limit}
\end{equation} 

Combining~\eqref{eq:C1},~\eqref{eq:C2}, and~\eqref{eq:limit}, we have proven that
\begin{equation}
\limsup_{\varepsilon\to 0} \left|\int_{\mathbb{R}}(u(t,x;\epsilon)-u^\mathrm{ZD}(t,x;\epsilon)) u^\mathrm{ZD}(t,x;\epsilon)\,\dd x\right|
	\leq \sqrt{\frac{C_1}{R}}+\sqrt{C_2\delta}.
\end{equation}
Letting $R\to \infty$ and $\delta\to 0$ then completes the proof.
\end{proof}

\section{Numerical verification of Theorem \ref{thm:u-app}}\label{sec:numerics}
\begin{figure}
    \centering
    \includegraphics[width=0.6\linewidth]{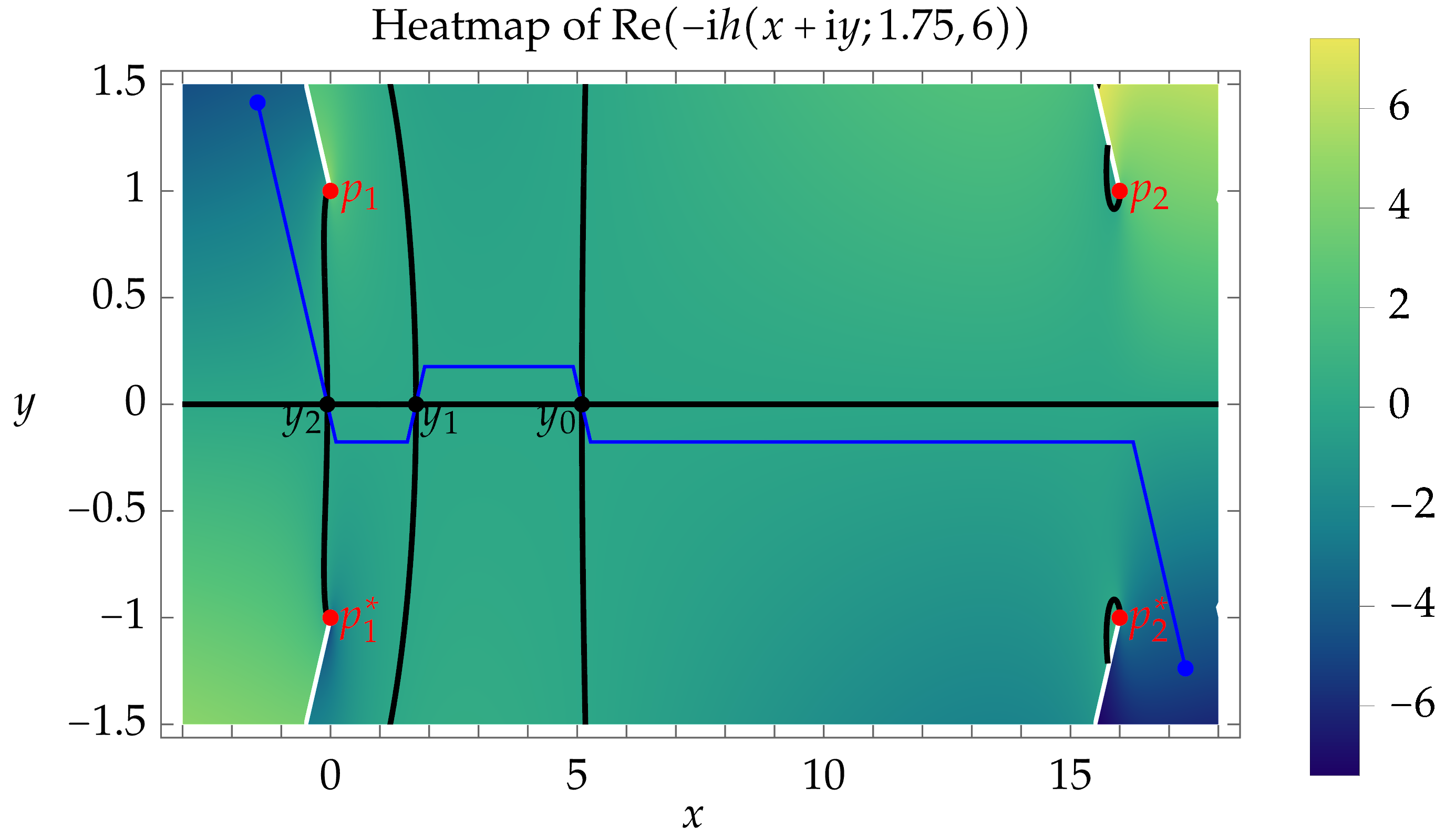}
    \caption{A heatmap of $\mathrm{Re}(-\ii h(z;1.75,6))$.  The black contours are the level set $\mathrm{Re}(-\ii h(z;1.75,6))=0$, the white lines are the branch cuts of $h(z;t,x)$, and the blue path is $C_{0,1}$, a truncation of $C_0$, which passes through the relevant saddle points $y_k$, $k=0,1,2$.}
    \label{fig:-ihheatmap}
\end{figure}
\begin{figure}
\begin{center}
    \includegraphics[width=0.48\linewidth]{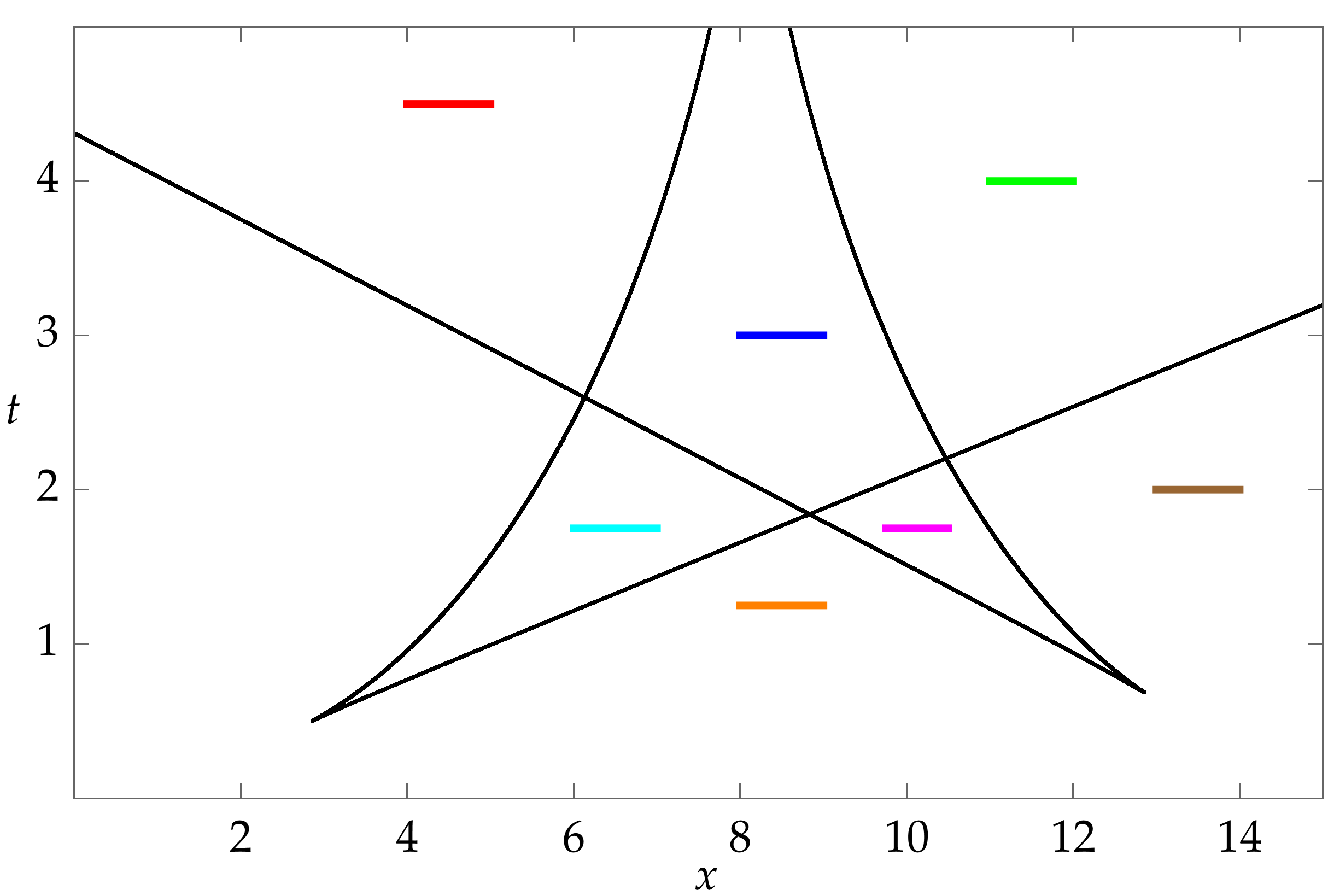}\hfill%
    \includegraphics[width=0.48\linewidth]{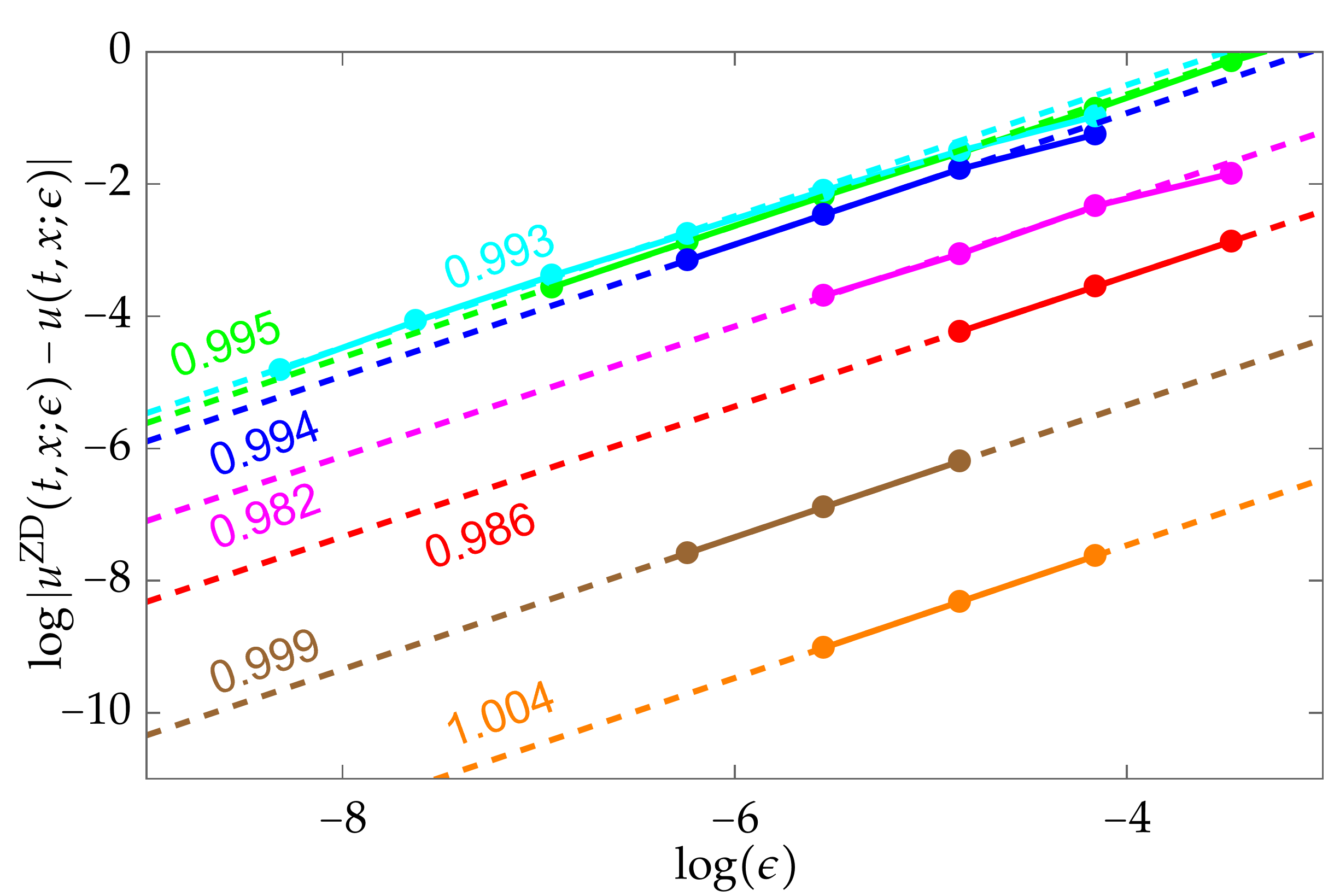}
\end{center}
\caption{The sup-norm of $u^{\mathrm{ZD}}(t,x;\epsilon)-u(t,x;\epsilon)$ is taken along the colored intervals in the left pane for various values of $\epsilon=2^{-k}$, with $u_0(x)$ being the same as in Figure \ref{fig: Burgers}.  The corresponding colored-dashed lines in the right pane are lines of best fit for the data points $(\log|u^{\mathrm{ZD}}(t,x;2^{-k})-u(t,x;2^{-k})|,\log(2^{-k}))$. The colored number is slope of the corresponding line of best fit.}
\label{fig: loglog}
\end{figure}

To demonstrate the validity of Theorem \ref{thm:u-app} numerically, we compute the sup-norm of $u^{\mathrm{ZD}}(t,x;\epsilon)-u(t,x;\epsilon)$ over various intervals in the $(x,t)$-plane (avoiding the caustic curves, see Figure \ref{fig: loglog}, left pane) with the initial condition $u_0(x)$ as in Figure \ref{fig: Burgers}.  Evaluating $u(t,x;\epsilon)$ boils down to computing integrals of the form
\begin{align}\label{eq:typical u int}
    \int_{W_k}f(z)\ee^{-\ii h(z;t,x)/\epsilon}\,\dd z, ~~~ k=0,1,\ldots,N,
\end{align}
see \eqref{eq:Amatrix}, \eqref{eq:Bmatrix}, Theorem \ref{thm:inversion-formula}, and Remark \ref{rem:change-contours}.  To accurately evaluate such integrals numerically, we deform the contours $W_k$, $k=0,\ldots,N$, to pass through the relevant saddle points such that $\mathrm{Re}(-\ii h(z;t,x))$ is maximized only when $z=z_s$ is a saddle point, see Section \ref{section:StokesGraphs}.  Then, we write \eqref{eq:typical u int} as
\begin{align}
    \int_{W_k}f(z)\ee^{-\ii h(z;t,x)/\epsilon}\,\dd z=\int_{W_{k,1}}f(z)\ee^{-\ii h(z;t,x)/\epsilon}\,\dd z+\int_{W_{k,0}}f(z)\ee^{-\ii h(z;t,x)/\epsilon}\,\dd z,
\end{align}
where the path $W_{k,1}$ contains all the relevant saddle points and $W_{k,0}$ is such that
\begin{align}
    \left|\frac{f(z_s)\ee^{-\ii h(z_s;t,x)/\epsilon}}{f(z)\ee^{-\ii h(z;t,x)/\epsilon}}\right|\leq10^{-16},
\end{align}
for all $z\in W_{k,0}$, where $z_s$ is any relevant saddle point.  See Figure \ref{fig:-ihheatmap}, blue contour, for an example of the path $W_{0,1}$.  Thus, we evaluate $u(t,x;\epsilon)$ numerically by replacing the paths $C_k$ with $W_{k,1}$ in~\eqref{eq:lambda-formula}.  Now we can approximate the sup-norm
\begin{align}\label{eq:supnorm}
    \sup_{x\in[a,b]}|u(t,x;\epsilon)-u^{\mathrm{ZD}}(t,x;\epsilon)|
\end{align}
for fixed $t$ by computing
\begin{align}\label{eq:maxnorm}
    \max_{k=0,1,\ldots,m}|u(t,x_k;\epsilon)-u^{\mathrm{ZD}}(t,x_k;\epsilon)|, ~~~ x_k=a+\frac{k}{m}(b-a).
\end{align}
The interval $[a,b]$ in \eqref{eq:supnorm} will be chosen to be the colored intervals in Figure \ref{fig: loglog}, left panel.  In the $J=0$ regions, $u(t,x;\epsilon)\simeq u^{\mathrm{B}}_0(t,x)$ is a slowly varying function, thus $m$ need not be large in \eqref{eq:maxnorm} to capture a good approximation of the sup-norm \eqref{eq:supnorm}.  For $J\geq1$, $u(t,x;\epsilon)$ becomes more oscillatory the smaller $\epsilon$ becomes so larger values of $m$ need to be taken in these regions.  For example, taking $t=4.5$ and the interval $[a,b]=[4,5]$ (the red interval in Figure \ref{fig: loglog}), we obtain the following data
\[
\begin{tabular}{|c|c|c|}
\hline
$\epsilon$ & $m$ & max \\
\hline
$2^{-5}$ & $1000$ & $0.056871$\\
$2^{-6}$ & $10000$ & $0.028720$\\
$2^{-7}$ & $10000$ & $0.014503$\\
\hline
\end{tabular}
\]
where `max' in the table refers to the max of \eqref{eq:maxnorm}.  Repeating this process for each of the colored intervals in Figure \ref{fig: loglog}, left panel, we create the log-log plot in the right panel of Figure \ref{fig: loglog}.  The lines of best fit are created by using the data points from the three smallest values of $\epsilon$ for each corresponding color.  Each line of best fit has slope nearly 1, which is consistent with Theorem \ref{thm:u-app} and demonstrates that the error term $\mathcal{O}(\epsilon)$ is sharp.

\begin{figure}
\begin{center}
    \includegraphics[width=0.48\linewidth]{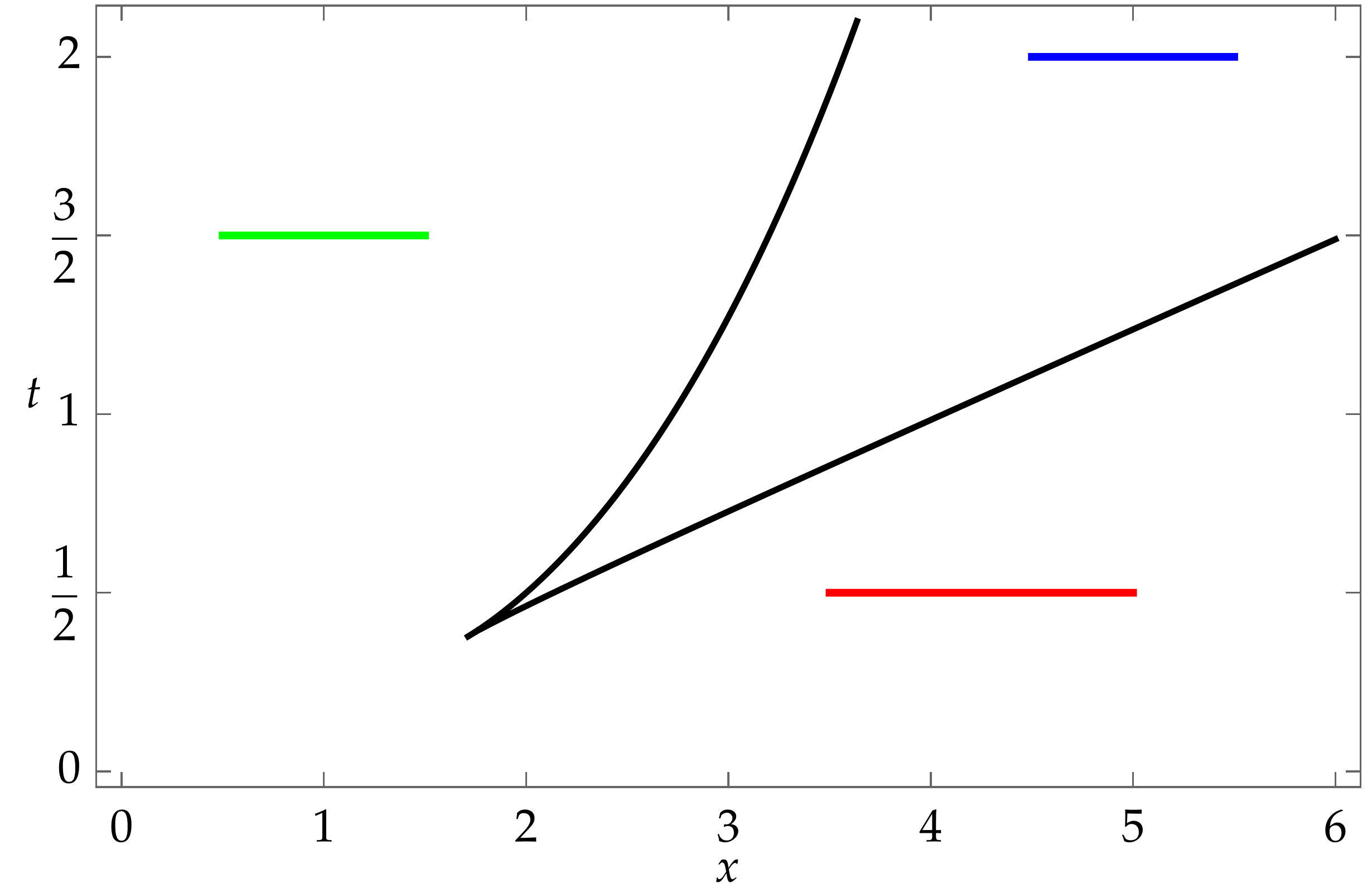}\hfill%
    \includegraphics[width=0.48\linewidth]{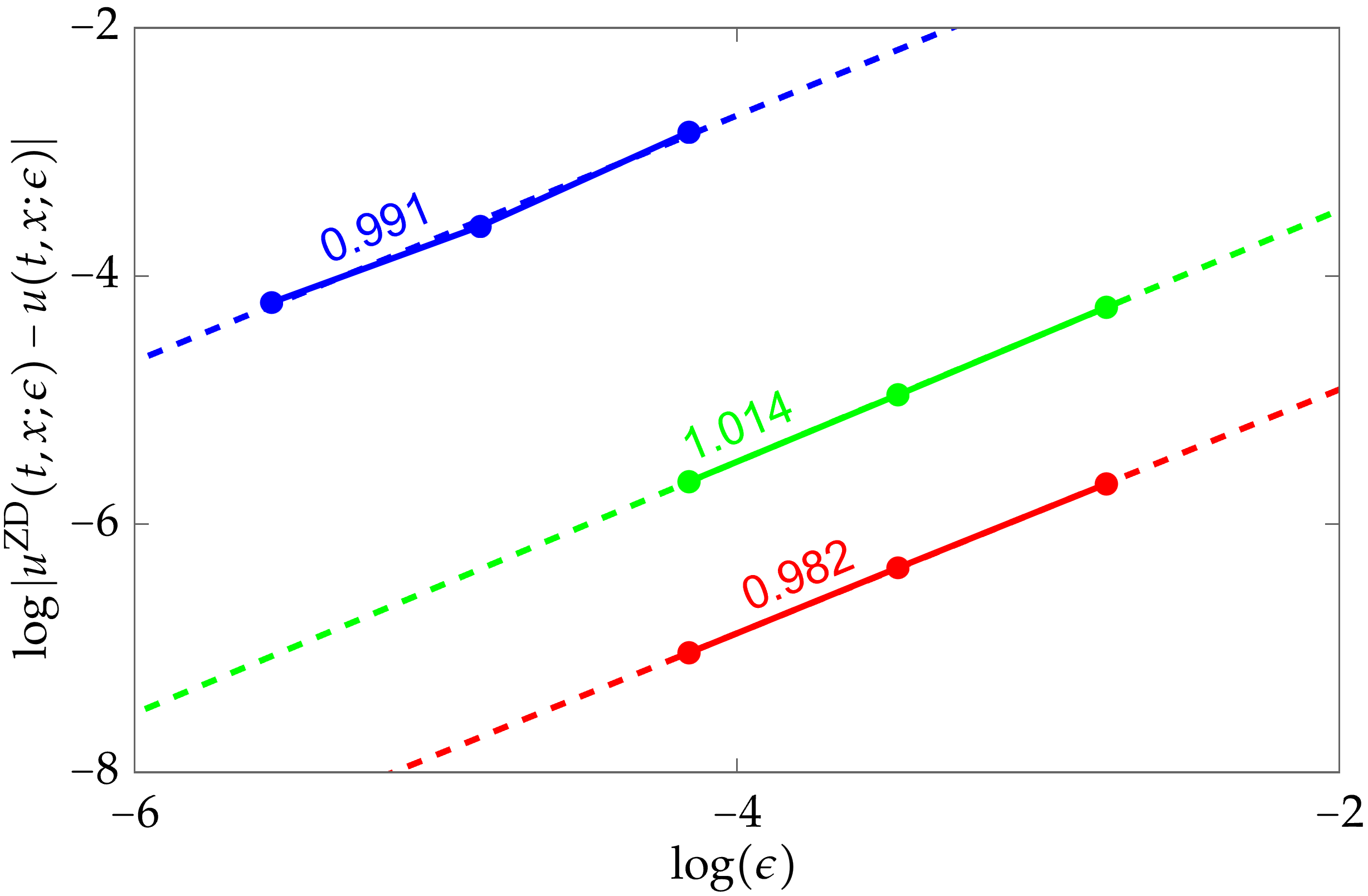}
\end{center}
\caption{As in Figure~\ref{fig: loglog}, but for the initial condition $u_0(x)=2/(1+x^2)$.}
\label{fig: loglog2}
\end{figure}

We repeat this process for the initial condition
\begin{align}
    u_0(x)=\frac{2}{1+x^2}=\frac{-\ii}{x-\ii}+\frac{\ii}{x+\ii}.
\end{align}
This example allows for a different method of computing $u(t,x;\epsilon)$ exactly without contour integration.  Indeed,
for $\epsilon=N^{-1}$ with $N\in\mathbb{N}$, $u(t,x;\epsilon)$ can be expressed explicitly via Matsuno's $N$-soliton formula \cite{Matsuno79} as
\begin{align}\label{u Matsuno}
    u(t,x;\epsilon)=2\epsilon\frac{\partial}{\partial x}\mathrm{Im}\left(\log\left(\det\left(\mathbb{I}+\frac{\ii}{\epsilon}\widetilde{\mathbf{A}}(t,x;\epsilon)\right)\right)\right),
\end{align}
where $\widetilde{\mathbf{A}}(t,x;\epsilon)$ is the $N\times N$ Hermitian matrix with elements
\begin{align}\label{eq: def A}
    A_{jk}(t,x;\epsilon)
    =\begin{cases}
    \displaystyle \frac{2\ii\epsilon\sqrt{\lambda_j\lambda_k}}{\lambda_j-\lambda_k}, &j\neq k, \\
    -2\lambda_j(x+2\lambda_jt), &j=k,
    \end{cases}
\end{align}
and $\lambda_j$, $j=1,\ldots,N$, are zeros of the $N^{\text{th}}$ degree Laguerre polynomial $L_N(-2\lambda/\epsilon)$, as was first discovered in \cite{KodamaAblowitzSatsuma82}.  We have created Figure~\ref{fig: loglog2} in the same fashion as Figure~\ref{fig: loglog} but we have used~\eqref{u Matsuno} for $u(t,x;\epsilon)$ instead of~\eqref{eq:lambda-formula}.  This experiment therefore demonstrates the validity of both the new explicit formula~\eqref{eq:lambda-formula} of $u(t,x;\epsilon)$ and the approximation by the zero-dispersion profile $u^{\mathrm{ZD}}(t,x;\epsilon)$ claimed in Theorem~\ref{thm:u-app}.

\bibliographystyle{siamplain}
\bibliography{references}

\end{document}